\documentclass[graybox]{svmult}

\usepackage{type1cm}        
\usepackage{makeidx}         
\usepackage{graphicx}        
\usepackage{multicol}        
\usepackage[bottom]{footmisc}

\usepackage{newtxtext}       %
\usepackage[varvw]{newtxmath}       
\usepackage{amsmath}
\usepackage{fancyhdr}

\makeindex             

\let\Bbb\mathbb
\let\goth\mathfrak

\def\og{\leavevmode\raise.3ex\hbox{$\scriptscriptstyle\langle\!\langle\,$}}

\def\fgf{\/\leavevmode\raise.3ex\hbox{$\scriptscriptstyle\,\rangle\!\rangle$}}

\def\fg{\fgf\ }

\newcommand{\carre}{\qed}

\let\stz\ss

\def\Demd#1|{\parindent=0pt\par{\sl D\'emonstration d#1}\pointir\parindent=20pt}

\let\petcap\sc

\def\finc{\vskip12pt}

\def\Dem{\parindent=0pt\par{\sl D\'emonstration}\pointir\parindent=20pt}

\def\Demo#1|{\parindent=0pt\par{\sl D\'emonstration #1}\pointir\parindent=20pt}

\makeatletter
\newdimen\indentTh\indentTh=0pt
\def\p@int{{\rm .}}
\def\p@intir{\discretionary{\rm .}{}{\rm .\kern.35em---\kern.7em}}
\def\pointir{\afterassignment\pointir@\let\next=}
\def\pointir@{\ifx\next\par\p@int\else\p@intir\fi\next}
\long\def\Thc#1|#2\finc{\Th{}{#1}{\pointir}{#2}}
\long\def\Th#1#2#3#4{\parindent=\indentTh\par\vskip5pt
{#1}{\petcap #2}{\sl #3}\parindent=20pt{\sl #4\par}\vskip 5pt\parindent=20pt}

\newdimen\indentssec\indentssec=20pt
\newdimen\indentrem\indentssec=0pt

\def\Demdsp#1|{\parindent=0pt\par{\sl D\'emonstration d#1.}\parindent=20pt}

\def\build#1_#2^#3{\mathrel{\mathop{\kern 0pt#1}\limits_{#2}^{#3}}}

\def\hdfl#1#2{\smash{\mathop{\hbox to 12mm{\rightarrowfill}}
\limits^{\scriptstyle#1}_{\scriptstyle#2}}}

\def\hdhfl#1#2{\smash{\mathop{\hbox to 12mm{\hookrightarrowfill}}
\limits^{\scriptstyle#1}_{\scriptstyle#2}}}

\def\hgfl#1#2{\smash{\mathop{\hbox to 12mm{\leftarrowfill}}
\limits^{\scriptstyle#1}_{\scriptstyle#2}}}

\def\hghfl#1#2{\smash{\mathop{\hbox to 12mm{\hookleftarrowfill}}
\limits^{\scriptstyle#1}_{\scriptstyle#2}}}

\newcount\secno \secno=0
\newcount\ssecno \ssecno=0
\newcount\sssecno \sssecno=0
\newcount\chapno \chapno=0
\newcount\notenumber \notenumber=1
\newcount\exino \exino=0
\newcount\expno \expno=0
\newdimen\indentsec\indentsec=20pt
\newdimen\indentssec\indentssec=20pt
\newdimen\indentsssec\indentsssec=20pt
\newdimen\indentrem\indentssec=0pt

\newdimen\indentTh\indentTh=0pt

\newdimen\indentth\indentssec=0pt
\def\sectiongen#1#2#3{\parindent=\indentsec\par\vskip .3cm
\vskip 0mm plus -20mm minus 1,5mm\penalty-50
{\bf #1}{\bf #2}{#3}\nobreak\parindent=20pt}%
\def\secc#1|{\sectiongen{}{#1}{\pointir}}
\def\nsecc#1|%
{\global\advance\secno by 1\global\ssecno=0\global\sssecno=0
\sectiongen{\the\secno\ }{#1}{\pointir}}
\def\secp#1|{\sectiongen{}{#1}{}\par}
\def\nsecp#1|%
{\global\advance\secno by 1\global\ssecno=0\global\sssecno=0
\sectiongen{\the\secno\ }{#1}{}\par}
\def\ssectiongen#1#2#3#4{\parindent=\indentssec\par\vskip .2cm
\vskip 0mm plus -20mm minus 1,5mm\penalty-50
{\bf #1}{\sl #2}{\sl #3}{#4}\nobreak\medskip\parindent=20pt}%
\def\ssecc#1|#2{\ssectiongen{}{#1}{\pointir}{#2}}
\def\nssecc#1|#2{\global\advance\ssecno by 1\global\sssecno=0
\ssectiongen{\the\secno.\the\ssecno\ }{#1}{\pointir}{#2}}
\def\ssecp#1|{\ssectiongen{}{#1}{}{}\par}
\def\nssecp#1|{\global\advance\ssecno by 1\global\sssecno=0
\ssectiongen{\the\secno.\the\ssecno\ }{#1}{}{}\par}
\def\sssectiongen#1#2#3#4{\parindent=\indentsssec\par\vskip .2cm
\vskip 0mm plus -20mm minus 1,5mm\penalty-50
{\bf #1}{\sl #2}{\sl #3}{#4}\nobreak\medskip\parindent=20pt}%
\def\sssecc#1|#2{\sssectiongen{}{#1}{\pointir}{#2}}
\def\nsssecc#1|#2{\global\advance\sssecno by 1
\ssectiongen{\the\secno.\the\ssecno.\the\sssecno\ }{#1}{\pointir}{#2}}
\def\sssecp#1|{\sssectiongen{}{#1}{}{}\par}
\def\nsssecp#1|{\global\advance\sssecno by 1
\sssectiongen{\the\secno.\the\ssecno.\the\sssecno\ }{#1}{}{}\par}
\long\def\Th#1#2#3#4{\parindent=\indentTh\par\vskip5pt
{#1}{\petcap #2}{\sl #3}\parindent=20pt{\sl #4\par}\vskip 5pt\parindent=20pt}
\long\def\pTh#1#2#3#4{\parindent=\indentth\par\vskip5pt
{#1}{\small \bf #2}{\small \sl #3}\parindent=20pt{\small \sl #4\par}\vskip 5pt\parindent=20pt}
\long\def\remarque#1#2#3#4{\parindent=\indentrem\par\vskip5pt
{#1}{\small \sl #2}{\small\sl #3}\parindent=20pt{\small#4\par}
\vskip 5pt\parindent=20pt}
\long\def\remarques#1#2#3#4{\parindent=\indentrem\par\vskip5pt
{#1}{\small \sl #2}{\small\sl #3}{\small#4\par}
\vskip 5pt\parindent=20pt}

\long\def\remarquesa#1#2#3#4#5{\parindent=\indentrem\par\vskip5pt
{#1}{\small \sl #2}{\small\sl #3}{\small#4}
\parindent=20pt{\small#5\par}
\vskip 5pt\parindent=20pt}

\long\def\Remarque#1#2#3#4{\parindent=\indentrem\par\vskip5pt
{#1}{ \sl #2}{\sl #3}\parindent=20pt{#4\par}
\vskip 5pt\parindent=20pt}
\long\def\remarquesn#1#2#3#4{\parindent=\indentrem\par\vskip5pt
{\small \sl #1}{\small #2}{\small\sl #3}{\small#4\par}
\vskip 5pt\parindent=20pt}

\long\def\Remarquen#1#2#3#4{\parindent=\indentrem\par\vskip5pt
{ \sl #1}{#2}{\sl #3}\parindent=20pt{#4\par}
\vskip 5pt\parindent=20pt}

\long\def\Thc#1|#2\finc{\Th{}{#1}{\pointir}{#2}}
\long\def\Thnc#1|#2|#3\finnc{\Th{#1}{#2}{\pointir}{#3}}
\long\def\Exic#1|#2\finc{{\global\advance\exino by 1}\Remarquen%
{Exercice }{\the\exino}{ #1\pointir}{#2}}
\long\def\Expc#1|#2\finc{{\global\advance\expno by 1}\Remarquen%
{Exemple }{\the\expno}{ #1\pointir}{#2}}
\long\def\exic#1|#2\finc{{\global\advance\exino by 1}\remarquesn%
{\bf Exercice }{\the\exino}{ #1\pointir}{#2}}
\long\def\expc#1|#2\finc{{\global\advance\expno by 1}\remarquesn%
{Exemple }{\the\expno}{ #1\pointir}{#2}}

\long\def\Ec#1\finc{\Th{}{}{}{#1}}
\long\def\thc#1|#2\finc{\pTh{}{#1}{\pointir}{#2}}%
\long\def\thsnc#1\finc{\pTh{}{}{}{#1}}

\long\def\Thp#1|#2\finp{\Th{}{#1}{\par}{#2}}
\long\def\thp#1|#2\finp{\pTh{}{#1}{\par}{#2}}
\long\def\rmc#1|#2\finc{\remarque{}{#1}{\pointir}{#2}}
\long\def\Rmc#1|#2\finc{\Remarque{}{#1}{\pointir}{#2}}

\long\def\rmp#1|#2\finp{\remarque{}{#1}{\par}{#2}}
\long\def\Rmp#1|#2\finp{\Remarque{}{#1}{\par}{#2}}

\long\def\parc#1\finc{\remarque{}{}{}{#1}}
\long\def\parcs#1\fincs{\remarques{}{}{}{#1}}

\long\def\parcsa#1\fins#2\fincsa{\remarquesa{}{}{}{#1}{#2}}

\def\Rm#1|{\parindent=0pt\par\vskip5pt{\sl #1}\pointir\parindent=20pt}
\def\preuved#1|{\parindent=0pt\par{\sl Preuve d#1}\pointir\parindent=20pt}
\def\qed{\quad\hbox{\hskip 1pt\vrule width 4pt height 6pt
         depth 1.5pt\hskip 1pt}}

\def\ffindem{\hfill\nobreak\penalty 1000 \hbox{\qed}\allowbreak\par\vskip 3pt}
\def\ffincdem{
			\hfill\nobreak\penalty 500
			\hbox{\carre\nobreak\qed}\allowbreak\par\vskip 3pt}
\def\fcarre{\hfill
	\hbox{\font\ppreu=cmsy10\ppreu\char'164\hskip-6.66666pt\char'165}}

\def\limproj{\mathop{\oalign{lim\cr\hidewidth$\longleftarrow$\hidewidth\cr}}}
\def\Indf{\parindent=0pt\par{\sl Induction finale}\pointir\parindent=20pt}
\def\ieme{\raise 1ex\hbox{\pc{}i\`eme|}\ }
\def\iemes{\raise 1ex\hbox{\pc{}i\`emes|}\ }
\input boxedeps.tex
\SetDVIWindowEPSFSpecial
\SetEPSFDirectory{}
\SetDefaultEPSFScale{250}
\EmulateRokicki

\global\indentsec=0pt\global\ssecno=0\global\secno=0
\global\ssecno=0\global\secno=0
\def\page {\leaders\hbox to 2mm{\hfil.\hfil}\hfill
}

\def\page#1#2{\leaders\hbox to 1mm{\hfil.\hfil}\hfill
\rlap{\hbox to 10mm{\hfill #1}{\ \  p.#2}}\par
}
\def\tpage #1{\leaders\hbox to 2mm{\hfil.\hfil}\hfill
\rlap{\hbox to 15mm{\hfill #1}}\par
}

\begin{document}
\fancyhf{}
\pagestyle{fancy}

\title*{Uniformisation des surfaces de Riemann}
\vskip5mm
\author{{\it par\/}  Alexis Marin, illustrations de Dorothea Vienne-Pollak}
\maketitle
\vskip10mm
\centerline{\Large  Mythe}
\centerline{\small (qui fut)}
\vskip2mm
\centerline{\bf contemporain}
\vskip50mm
\centerline{\sl d'apr\`es une id\'ee originale de}
\vskip2mm
\centerline{ J.-P.-A. Douaill\^y Jr.}
\vfill\eject
\thispagestyle{empty}
\null\parc
Classification de l'A.M.S. : 30-00, 30F10, 30C20, 57N05.
\vskip5mm
Classification de la librairie du congr\`es :

\quad Mythologie profane.

\quad Invertissement analytico-alg\`ebrique.

\quad Divertissement holomorphe.

\quad Exercice de style g\'eom\'etrique.

\vskip20mm

\centerline{R\'esum\'e}
\vskip2mm
De moulte lemmes, corollaires, cocorollaires\footnote%
{{\it c.a.d.\/} corollaire d'un corollaire.},
defemmes, decofaires\footnote%
{{\it c.a.d.\/} \'enonc\'es int\'egrant lemme ou corollaire \`a une d\'efinition le n\'ecessitant.}
subreptissement sorit\'es dans les appendices, les l\`evres humides et l'esprit clarifi\'e,
l'uniformisation,
 d\'elaissant
la savante analyse, sa Th\'eorie du potentiel et ses  Equations aux d\'eriv\'ees partielles,
c\`ede aux charmes \og \'el\'ementaires\fg de la na{\"{\i}}ve, mais efficace, r\`egle des signes :

$$\hbox{\rm Pour \/} m\ \hbox{\rm et \/} n\
\hbox{\rm entiers \/} (-1)^m\,(-1)^n=(-1)^{m+n}\/$$
\finc
\vskip15mm
{\leftskip+23mm\rightskip+25mm\parindent=0pt
{\small\it E\'ennodrooque, H\'et\'erolution, ainsi que la plupart des
termes du glossaire sont des barbarismes d\'epos\'es dont
l'emploi est soumis
\`a l'obtention d'une licence.

En attendant que les conditions d'abonnement et les
certificats correspondants soient disponibles dans les bureaux de
tabac, une tol\'erance est accord\'ee pour usage exclusif de compliment et/ou
quolibet de cour de r\'ecr\'eation.

Pour tout autre usage, public ou priv\'e, m\^eme licencieux ou domestique,
les demandes d'abonnement s'obtiennaient

\centerline{(jusqu'au Mon, 1 Mar 2004 19:06:59 +0100)}
par courier \'electronique
\`a :
\vskip2mm
\centerline{"Annales de l'Institut Fourier (Monique Vitter et Nathalie Catrain)"}
\centerline{ annalif@ujf-grenoble.fr}}
\par}

\vfill\eject
%
\pagenumbering{roman}
\setcounter{page}{1} 
\fancyhead[LO]{A. M. \quad Uniformisation des surfaces de Riemann}
  \fancyhead[RO]{\thepage}
  \fancyhead[RE]{\hfill\small  Aper\c cu sur quatre paragraphes et autant d'appendices \hfill}
  \fancyhead[LE]{\thepage}

\null\vskip1mm
\centerline{\Large Aper\c cu sur quatre paragraphes et autant d'appendices}
\vskip6mm
\centerline{r\'eduisant}
\vskip4mm
\centerline{\bf le th\'eor\`eme de Koebe \`a celui de Riemann dans le plan\/}
\vskip20mm
Une {\it surface de Riemann\/} ${\goth S\/}\/$,
{\it i.e.\/} vari\'et\'e holomorphe s\'epar\'ee
de dimension com\-plexe $1\/$, 
est {\it simplement connexe\/}
si tous ses rev\^etements
sont triviaux.

 Plan complexe,
disque unit\'e
et sph\`ere de
Riemann
sont des surfaces de  Riemann connexes et simplement connexes.
Le th\'eor\`eme
d'uniformisation affirme :

{\sl A {\it isomorphisme\/},\/
{\rm {\it i.e.\/} hom\'eomorphisme holomorphe\/}, pr\`es ce sont les seules.\/}

Ce r\'esultat \'enonc\'e par Riemann au \S XXI de sa dissertation de 1851,
mais avec une d\'emonstration qui n'a \'et\'e  pleinement justifi\'ee
par Hilbert qu'en 1909, est usuellement
 attribu\'e ind\'ependamment
\`a Poincar\'e et \`a Koebe
dans deux articles de 1907 et
a fait couler beaucoup d'encre$^{0}\!\!\/$.
Sans pr\'etendre, comme le dit
Hermann Weyl de Koebe, avoir, la vie enti\`ere
\og darauf verwendet, das Problem der
Uniformisierung nach allen Richtungen und mit den verschiedensten
Methoden durchzudenken\fgf$^{1}\/$
le pr\'esent texte d\'eduit le th\'eor\`eme d'uniformisation
du \og th\'eor\`eme de Riemann dans le plan\fg de
tous les manuels d'analyse complexe$^{2}\/$ :
\Thc Th\'eor\`eme 0| Un ouvert connexe strict du plan complexe, sur lequel
toute fonction  holomorphe ne prenant pas la valeur z\'ero
a une racine carr\'ee holomorphe, est isomorphe au disque unit\'e.\finc
Un \'enonc\'e de Jordan$^{3}\/$
caract\'erise donc les {\it domaines de Jordan\/}, {\it i.e.\/} int\'erieurs de sous-vari\'et\'e
topologique du plan
de bord connexe non vide :
\Thc Corollaire 0| Un domaine de Jordan est isomorphe au disque unit\'e.
\finc

La r\'eduction ici propos\'ee est \og \'el\'ementaire\fgf, en le sens
que d'une part elle ne n\'ecessite pas d'outil homologique, ni de triangulation
(elle donne en particulier la paracompacit\'e$^{4}\/$
des surfaces de Riemann)
et d'autre part, contrairement aux diff\'erents trait\'es modernes
sur les surfaces de Riemann$^{5}\/$,
la seule analyse qu'elle utilise sont les r\'esultats de base
sur les fonctions holomorphes.
\vfill\eject
\fancyhead[LO]{A. M. \quad Uniformisation des surfaces de Riemann}
  \fancyhead[RO]{\thepage}
  \fancyhead[RE]{\S 1 Le fil des \'enonc\'es}
  \fancyhead[LE]{\thepage}
\null\vskip5mm
\centerline{\bf La preuve et ses sept \'enonc\'es}
\vskip5mm
\centerline{\small\bf Le cas compact\/}
\vskip10mm
Le \S 1 \'etablit une premi\`ere caract\'erisation de la sph\`ere de Riemann :
%
\Thc Lemme A|  Une surface de Riemann compacte ${\goth S\/}\/$
union de deux ouverts isomorphes  au disque unit\'e 
 est isomorphe \`a la sph\`ere de Riemann.
\finc
%
Ici la seule
hypoth\`ese de nature topologique sur la surface de Riemann ${\goth S\/}\/$ est la compacit\'e
et il n'y a pas d'hy\-po\-th\`ese de r\'egularit\'e de la fronti\`ere.
Il n'en sera pas de m\^eme pour caract\'eriser dans 
les surfaces de Riemann les ouverts isomorphes au disque, ici d\'enomm\'es {\it \'el\'ementaires\/}, ou ceux,
dits {\it standards\/},
dont les composantes connexes
sont isomorphes au disque ou \`a la sph\`ere de Riemann :

Un {\it cycle analytique\/} (ou {\it $\omega\/$-cycle\/})
de ${\goth S\/}\/$ est un ferm\'e $\Gamma\/$ de ${\goth S\/}\/$,
union localement finie d'une famille $(\beta_\lambda)_{\lambda\in \Lambda}\/$
d'arcs analytiques r\'eels  tel que tout point $p\in \Gamma\/$ est extr\'emit\'e
d'un nombre pair des $\beta_\lambda\/$.
Un {\it polygone analytique (ouvert) (\/{\rm ou\/} $\omega\/$-polyone (ouvert))\/} 
est un ouvert de ${\goth S\/}\/$, int\'erieur de sa fermeture, et de fronti\`ere un $\omega\/$-cycle.

Une surface de Riemann est {\it planaire\/} si tout $\omega\/$-cycle compact la s\'epare.

La construction, dans l'appendice 1, du rev\^etement double associ\'e \`a un $\omega\/$-cycle
  \'etablit qu'une surface de Riemann simplement connexe est planaire.

Le \S 2 d\'emontre l'analogue du  Lemme A pour les $\omega\/$-polygones planaires :
%
\Thc Lemme C| Un polygone analytique relativement compact $\Omega\/$ dans une surface de
Riemann planaire ${\goth P\/}\/$, 
de fermeture $\overline{\Omega}\/$ recouverte par deux ouverts standards
$U_1\/$ et $U_2\/$, et d'int\'erieur du compl\'ementaire
${\goth P\/}\!\setminus\!\overline{\Omega}\/$ connexe, est standard.
\finc
%
{\parindent=0pt\par\vskip -.1cm
\vskip 0mm plus -20mm minus 1,5mm\penalty-50
qui, au \S 4, implique une seconde caract\'erisation de la
sph\`ere de Riemann :
\nobreak\parindent=20pt}%
%
\Thc Th\'eor\`eme 1| Une surface de Riemann compacte connexe planaire ${\goth P\/}\/$
 est isomorphe
\`a la sph\`ere de Riemann.
\finc
%
{\parindent=0pt\par\vskip -.1cm
\vskip 0mm plus -20mm minus 1,5mm\penalty-50 
d'o\`u le cas compact du th\'eor\`eme d'uniformisation. 
\nobreak\parindent=20pt}%
\vfill\eject
\centerline{\bf Arrondissement des brisures du bord et le cas g\'en\'eral\/}
\vskip3mm
Le \S 2 utilise l'extension,
 aux $\omega\/$-polygones,
du {\it double de Klein\/} d'une surface de Riemann \`a bord ${\bf R\/}\/$,
une surface de
Riemann ${\goth D\/}\,{\bf R\/}\supset {\bf R\/}\/$ la contenant et munie d'une
{\it h\'et\'erolution\/}, {\it i.e.\/} involution conforme nulle part holomorphe,   $\sigma_{\bf R\/}\/$
 d'ensemble des points fixes ${\rm Fix\/}\,\sigma_{\bf R\/}=\partial{\bf R\/}\/$ le bord  de ${\bf R\/}\/$ et
ayant ${\bf R\/}\/$ pour domaine fondamental. Le cas usuel rappel\'e dans l'appendice 2, cette extension
occupe l'appendice 3 : 

Un {\it bord\/} $\beta\/$ d'un $\omega\/$-polygone $U\/$
est une courbe param\'etr\'ee $\beta:B\rightarrow{\rm Fr\/}\,U\/$ de sa fronti\`ere, dont
tout param\`etre $s\in B\/$ a un voisinage d'image par $\beta\/$
l'union de deux arcs non s\'epar\'ees pr\`es de $\beta\,(s)\/$ par le compl\'ementaire de $U\/$ et analytique
injective hors de $\beta^{-1}\,(\Delta_U)\/$, ses {\it coins\/},  pr\'eimage
du discret ensemble $\Delta_U\/$ des {\it sommets\/} de  $U\/$.

Il y a un bord d'image ${\rm Fr\/}\,U\/$, le {\it bord $\partial U\/$ de\/} $U\/$, et
 pour tout bord $\beta\/$ de $U\/$, 
une surface de Riemann  ${\goth D\/}_\beta\,U\/$, le {\it double  de \/{\rm l'$\omega\/$-plolygone\/}
$U\/$ sur\/} le bord $\beta\/$, munie
de $\sigma_{\beta}\/$ et $i_\beta : U_\beta\rightarrow {\bf U\/}_{^a\beta}\/$, h\'et\'ero\-lution  et
hom\'eomorphisme, holomorphe sur $U\/$, 
de
$$U_\beta=(U\cup(\beta\,(B)\!\setminus \Delta_U))\amalg B/%
\{\beta\,(s)\sim s\ \hbox{\rm si \/} s\in B\ \hbox{\rm et \/}\beta\,(s)\notin\Delta_U\}\/$$
sur une sous-surface de Riemann \`a bord ${\bf U\/}_{^a\beta}\/$ de ${\goth D\/}_\beta\,U\/$, dite
{\it arrondie de $U\/$ sur $\beta\/$\/}, et dont
 $({\goth D\/}_\beta\,U,{\bf U\/}_{^a\beta},\, \sigma_{\beta})\/$ est double de Klein usuel.

Le double d'un $\omega\/$-polygone \'el\'ementaire sur
un arc non dense de sa fronti\`ere \'etant \'el\'e\-men\-taire,
le Lemme A donne
 une caract\'erisation du disque unit\'e,
le Lemme B, d'\'enonc\'e utilisant le vocabulaire ouvrant le \S 2, et
dont l'application r\'ecurrente \'etablira le Lemme C et, au \S 3, 
la structure des composantes compactes de bord permettant de  \og boucher\fg les surfaces de Riemann
compactes \`a bord :
\Thc Lemme D|
Une
composante compacte $X\/$ du bord $\partial\, {\bf R\/}\/$ d'une surface de Riemann \`a bord ${\bf R\/}\/$
a un voisinage dans ${\bf R\/}\/$ isomorphe \`a un voisinage dans le disque unit\'e ferm\'e, de son bord,
le cercle unit\'e.
\finc
\vskip-2mm
\Thc Corollaire D| Une surface de Riemann \`a bord compacte ${\bf R\/}\/$
 est isomorphe \`a une sous-surface ${\bf S\/}\/$ d'une surface de Riemann compacte
${\goth S\/}\/$ avec ${\goth S\/}\!\setminus\!{\bf S\/}\/$ standard.
\finc
Tout compact d'une surface de Riemann ${\goth Q\/}\/$,
\'etant inclus dans
un $\omega\/$-polygone relativement compact, a donc un voisinage isomorphe \`a un ouvert $U\/$ d'une surface de
Riemann compacte ${\goth S\/}\/$,  planaire si ${\goth Q\/}\/$ l'est. Au \S 4
le
crit\`ere de Montel, expos\'e \`a l'appendice 4 (et utilis\'e aussi aux
\S 1 et 3), justifiera  la terminologie \og planaire\fg :
\Thc Th\'eor\`eme 2| Une surface de Riemann connexe planaire ${\goth Q\/}\/$ est
isomorphe \`a un ouvert de la sph\`ere de Riemann.
\finc
{\parindent=0pt\par\vskip -.2cm
\vskip 0mm plus -20mm minus 1,5mm\penalty-50
qui ram\`ene le th\'eor\`eme d'uniformisation au th\'eor\`eme de Riemann dans le plan.
\nobreak\parindent=20pt}%
\centerline{\hfill$\clubsuit\ \diamondsuit\ \heartsuit\ \spadesuit\/$\hfill}

Les appendices soritteront alors
vers  commentaires, r\'ef\'erences, notes, index,  table des mati\`eres et 
g\'en\'erique de fin : \'epilogue d'une histoire dont, avant de la plus habiller, 
il convient de pr\'esenter les personnages  :
\vfill\eject
%
%
\fancyhead[LO]{A. M. \quad Uniformisation des surfaces de Riemann}
  \fancyhead[RO]{\thepage}
  \fancyhead[RE]{\quad\quad Pr\'elude\hfill  Notations de quatres paragraphes et autant d'appendices \hfill}
  \fancyhead[LE]{\thepage}

{\small
{\parindent=0pt\par
\vskip 0mm plus -20mm minus 1,5mm\penalty-50
\centerline{\bf Objets et morphismes usuels dans le plan et la sph\`ere}
\vskip5mm
\centerline{\sl Les mod\`eles et quelques morphismes et isomorphismes tant alg\'ebriques que conformes.\/}%
\nobreak\parindent=20pt}%
${\Bbb N\/},\ {\Bbb Z\/}\/$ , {\it Entiers naturels\/} et {\it relatifs\/}, (mo-,  an-)neau
des deux premi\`eres
op\'erations $+\/$ et $\times\/$.

${\Bbb F\/}_2=\{0,\,1\}\/$, {\it Corps  \`a deux
\'el\'ements\/} isomorphe au {\it corps ${\Bbb Z\/}/2\,{\Bbb Z\/}\/$ des entiers modulo $2\/$\/}.

${\Bbb C\/}=\{x+i\,y\,;\, x\ \hbox{\rm et \/} y\ \hbox{\rm r\'eels\/}\}\/$, \/{\it Plan complexe\/},
un corps \'etendant celui des r\'eels o\`u $i^2=-1\/$.

${\Bbb R\/}=\{x+i\,0\,;\, x\ \hbox{\rm r\'eel\/}\}\subset{\Bbb C\/}\/$, \/{\it Axe r\'eel\/}, un sous-corps de ${\Bbb C\/}\/$.

${\Bbb C\/}^{\ast}={\Bbb C\/}\!\setminus\!\{0\}\/$, \/{\it Plan complexe \'epoint\'e\/}, groupe multiplicatif du 
corps ${\Bbb C\/}\/$.

${\Bbb R\/}^{\ast}={\Bbb R\/}\!\setminus\!\{0\}\/$, \/{\it Axe r\'eel \'epoint\'e\/}, groupe multiplicatif du 
sous-corps ${\Bbb R\/}\/$ de ${\Bbb C\/}\/$.

$\mu_2\!=\!\{\pm\}\!=\!\{\pm1\}\!\subset\!{\Bbb R\/}^\ast\!\subset\!{\Bbb C\/}^\ast\/$, {\it Doublon des signes\/}, 
sous-groupe des racines carr\'ees de l'unit\'e.

${\Bbb R\/}_\epsilon=\{x+i\,0\,;\, \epsilon\,x\geq 0\}$, pour $\epsilon\in\mu_2\/$,
\/{\it Demi-axes positif\/} si $\epsilon=+\/$, \/{\it n\'egatif\/} si $\epsilon=-\/$.

$\Re,\, \Im :{\Bbb C\/}\rightarrow{\Bbb R\/}\/$, \/{\it Parties r\'eelle\/} 
et {\it imaginaire\/} associant \`a $z=x+i\,y,\  \Re z=x\/$ et $\Im z=y\/$.

${\rm conj\/}\/$, \/{\it Conjugaison complexe\/} : isomorphisme de corps 
$z=\Re z+i\,\Im z\mapsto {\bar z}=\Re z-i\,\Im z\/$.

$|\ | :{\Bbb C\/}\rightarrow{\Bbb R\/}_+\/$, \/{\it Module\/}, $z\mapsto \sqrt{z\,{\bar z}}=\sqrt{(\Re\,z)^2+(\Im z)^2}\/$,
il est multiplicatif : $|z\,z'|=|z|\,|z'|\/$.

$D=\{z\in{\Bbb C\/}\,;\, |z|<1\}\/$, \/{\it Disque (unit\'e (centr\'e en l'origine $0\/$),\/ {\rm ou\/} central)\/}.

$D^{\ast}=D\!\setminus\!\{0\}\/$, \/{\it Disque (unit\'e) \'epoint\'e\/}.

$D_+=\{z\in{\Bbb C\/}\,;\, |z|<1,\ \Im z\geq 0\}\/$,
\/{\it Demi-disque  (Nord unit\'e)\/}, de bord l'intervalle $]-1,\,1[\/$.

$S^1=\{z\in{\Bbb C\/}\,;\,|z|=1\}\/$, \/{\it Cercle unit\'e\/}. Sous- groupe compact
de ${\Bbb C\/}^\ast\!\!\/$, y ayant pour voisinage  :

$A_r=\{z\in{\Bbb C\/}\,;\,r<|z|<r^{-1}\}\/$, o\`u $0<r<1\/$, \/{\it Anneau invariant (de rayon int\'erieur $r\/$)\/}.

$\overline{D}=\{z\in{\Bbb C\/}\,;\,|z|\leq 1\},\ (\overline{D}^\ast=\overline{D}\setminus\{0\})\/$,
\/{\it Disque unit\'e ferm\'e (\'epoint\'e)\/}, de bord $S^1\/$.

$H=\{z\in{\Bbb C\/}\,;\,\Im z>0\}\/$
,  {\it Demi-plan 
de Poincar\'e
\/}.

$D_{z,\,r}=z+r\,D\/$, pour $z\in{\Bbb C\/},\, r>0\/$, {\it Disque de centre $z\/$ et rayon $r\/$\/},
sans centre $D_{r}=D_{0,\,r}\/$.

$\overline{D}_{z,\,r}=z+r\,\overline{D}\/$, pour $z\in{\Bbb C\/},\, r>0\/$,
{\it Disque ferm\'e de centre $z\/$ et rayon $r\/$\/},
$\overline{D}_{r}=\overline{D}_{0,\,r}\/$.

$S_{z,\,r}=z+r\,S^1\/$, pour $z\in{\Bbb C\/},\, r>0\/$, {\it Cercle de centre $z\/$ et rayon $r\/$\/},
sans centre $S_{r}=S_{0,\,r}\/$.

$P_1\,({\Bbb C\/})={\Bbb C\/}\cup\{\infty\}\/$, \/{\it Sph\`ere de Riemann\/}, contenant  
les espaces pr\'ec\'edents et o\`u,
selon la coutume [si $u\in{\Bbb C\/}, v\in{\Bbb C\/}^\ast\/$ alors
$\infty +u=\infty,\, v\,\infty=\infty\,v=\infty,\, \frac{v}{0}=\infty, \frac{v}{\infty}=0,\,
\overline{\infty}=\infty\/$],
agissent 
{\it homo\/} $z\mapsto\frac{a\,z+b}{c\,z+d}\/$ et {\it h\'et\'ero\/}
$z\mapsto\frac{a\,{\bar z}+b}{c\,{\bar z}+d}\/$ {\it graphies\/}
[$\ a, b, c, d\in{\Bbb C\/},\ a\,d-b\,c\!\ne\!0\/$]. Une graphie commutte \`a ${\rm conj\/}\/$ si et
seulement si elle est r\'eelle, {\it i.e.\/} a un repr\'esentant avec $a,\, b,\, c,\, d\in {\Bbb R\/}\/$.
\vskip1mm

$L_{\epsilon}=\{z\in{\Bbb C\/} ;\, \epsilon\,\Re\,z>0\}, \epsilon\in\mu_2\/$, 
\/{\it Demi-plan lat\'eraux\/}. 
{\it Sym\'etrie polaire\/} $\imath : z\mapsto -z^{-1}\/$
(centr\'ee aux \/{\it p\^oles \/} $\{\pm i\}\/$)
\'echange demi-plan \/{\it Est\/} $L_+\/$ et \/{\it Ouest\/} $L_-\/$, chacun
iso\-mor\-phe au disque central $D\/$ par
l'\/{\it $\epsilon\/$-r\'etrograde en $\epsilon\,i\/$ quart de tour polaire,\/} 
$\vartheta_\epsilon : L_\epsilon\rightarrow D,\, 
\vartheta_\epsilon\,(z)= (\frac{z-1}{z+1})^\epsilon\/$.

$I=[0,\,1]\subset{\Bbb R\/}\subset{\Bbb C\/}\/$, {\it Intervalle unit\'e\/}, d'\/{\it int\'erieur formel\/}
${\breve I}=I\!\setminus\!\{0,\,1\}\/$ ayant pour voisinages

$P_{N,\,h}\!=\!\{z\in{\Bbb C\/};\, h\,|\Im\,z|\!<\!(\Re z-|z|^2)^{N+1}\}\/$, pour $h\!>\!0\/$ et $N\/$ pair,
{\it Perle\/}  $(N,\,h)\/$-fine sur $I\/$.
\vskip-1mm

$\exp:{\Bbb C\/}\!\rightarrow\!{\Bbb C\/}^\ast,\
z\mapsto\build{\sum}_{n=0}^{\infty}\!\frac{x^n}{n!}\/$,
{\it Exponentielle\/}, morphisme $2\,\pi\,i\/$-p\'eriodique ${\Bbb C\/}\!\setminus\!{\Bbb R\/}_-\/$-scind\'e  par :

$\log : {\Bbb C\/}\!\setminus\!{\Bbb R\/}_-\rightarrow\{z\in{\Bbb C\/}\,;\, |\Im z|<\pi\}\/$,
isomorphisme \/{\it D\'etermination principale du logarithme\/}.

$\epsilon : {\Bbb Z\/}\rightarrow\mu_2\/$, {\it Parisigne\/}, morphisme $n\mapsto (-1)^n\/$,
induit l'isomorphisme $\underline{\epsilon\/} : ({\Bbb F\/}_2,\,+)\rightarrow \mu_2\/$.

\vfill\eject%

{\parindent=0pt\par\vskip 0cm
\vskip 0mm plus -20mm minus 1,5mm\penalty-50
\centerline{\sl D'autres homo (et h\'et\'ero) graphies utiles.\/}%
\nobreak\parindent=20pt}%
\vskip1mm

${\rm Id\/}\/$, \/{\it Identit\'e\/} $z\mapsto z\/$ de $P_1\,({\Bbb C\/})\/$.

$\imath_\lambda\/$, pour $\lambda\in{\Bbb C\/}^\ast$, est l'involution holomorphe 
$z\mapsto -\frac{\lambda}{z}\/$, donc $\imath_1=\imath\/$ est {\it Sym\'etrie polaire\/}.

$\sigma_t\/$, pour $t\in{\Bbb R\/}^\ast\/$, est l'h\'et\'erolution
$\sigma_t\!=\!{\rm conj\/}\circ\imath_{-t}\/$, libre si et seulement si $t<0\/$, pamis elles :

$\sigma_{-1}\/$, {\it Antipode\/}.

$\sigma_{r^2}\/$, \/{\it Inversion fixant (le cercle) $S_r\/$\/}, sans pr\'eciser, \/{\it Inversion \/}
est $\sigma_1\/$.

$\varphi_\alpha : D\rightarrow D\/$, pour $\alpha\in D\/$, est $ z\mapsto\frac{z-\alpha}{1-{\bar \alpha}\,z}\/$,
elle est inversible d'inverse $\varphi_{-\alpha}\/$.

$h_\lambda\/$, {\it Similitude de rapport complexe\/}  $\lambda\in{\Bbb C\/}^\ast\/$
 est $z\mapsto \lambda\,z$.

$h_t\/$, {\it Homoth\'etie de rapport\/} $t\/$ si $t\in{\Bbb R\/}^\ast\/$.

$\tau_v\/$, {\it Translation de vecteur \/} $v\in{\Bbb C\/}\/$  est $z\mapsto z+v\/$.

{\parindent=0pt\par\vskip 0mm
\vskip 0mm plus -20mm minus 1,5mm\penalty-50
\vskip3mm
\centerline{\sl Et les triviales\/}%
\nobreak\parindent=20pt}%
\vskip1mm
$c_k\/$,  Applications constantes sur $P_1\,({\Bbb C\/})\/$ d'image 
le singleton $\{k\}\/$.

\vskip10mm
{\parindent=0pt\par\vskip 0cm
\vskip 0mm plus -20mm minus 1,5mm\penalty-50
\centerline{\bf Op\'erations et invariants topologiques}
\nobreak\parindent=20pt}%
\vskip5mm
$X\amalg Y\/$, {\it Espace somme\/} (disjointe) des espaces topologiques $X\/$ et $Y\/$.

$\build{\amalg}_{\lambda\in\Lambda}^{}\,X_\lambda\!=\!\{(x,\,\lambda)\in X\times\Lambda; x\in X_\lambda\}\/$ somme disjointe
d'une famille 
de sous-espaces de $X\/$.
\vskip1mm
$\overline{A}\/$, {\it Fermeture\/} d'une partie $A\subset X\/$ d'un espace topologique $X\/$.

${\rm Int\/}\,A=\build{A}_{}^{\circ}\/$, {\it Int\'erieur\/} d'une partie $A\subset X\/$ d'un espace topologique $X\/$.

${\rm Int\/}_Y\,A=\build{A}_{}^{\circ}\cap Y\/$, {\it Int\'erieur relatif\/} dans un sous-espace
$Y\/$ d'une partie $A\/$ de l'espace $X\/$.

${\cal C\/}\,(X,\,Y)\/$, {\it Ensemble des applications continues\/} d'un espace $X\/$ dans un espace $Y\/$.

${\rm incl\/}_A^X\ ({\rm ou \/}\ A\hookrightarrow X\,)\in{\cal C\/}\,(A,\,X)\/$, {\it Inclusion\/} du sous-espace $A\/$ dans l'espace $X\/$.

$c\,X={\rm card\/}\,{\goth c\/}\,X\/$, {\it Connexit\'e\/} de $X\/$, cardinal   des
{\it composantes connexes\/}, sous-espace dense de 

${\bar {\goth c\/}}\,X={\goth X\/}\,A\,(X)\/$, {\it Compact des caract\`eres\/} de $A\,(X)\/$,
{\it c.a.d.\/}  morphismes
d'anneau sur ${\Bbb F\/}_2\/$   o\`u

 $A\,(X)={\cal C\/}\,(X,\,{\Bbb F\/}_2)\/$, {\it Alg\`ebre des ${\Bbb F\/}_2\/$-constantes locales\/} de $X\/$,
de groupe additif
isomorphe  \`a

${\goth E\/}\,X={\cal C\/}\,(X,\,\mu_2)\/$,
{\it  Groupe
des signes locaux\/},
 fini si et seulement si $c\,X\/$ l'est, en ce cas
$$c\,X={\rm dim\/}_{{\Bbb F\/}_2}\,{\goth E\/}\,X$$

}
%
Lemmes A \`a D, Th\'eor\`emes 0 \`a 2, Corollaires 0 et D, d\'ej\`a \'evoqu\'es,
sont d\'esormais ainsi d\'enomm\'es, mais
le $m^{\hbox{\scriptsize i\`eme}}$ \'enonc\'e de l'appendice $n\/$ sera cit\'e, \'eventuellement pr\'ec\'ed\'e d'un
\og nom d'auteur\fgf, par {\bf n.m\/} en caract\`eres gras.

Le chiffre arabe ${n}\/$ en exposant d'un mot$^{n}\/$, 
comme les $^{0}\/$ \`a $^{5}\/$ apparus page i,
indique au lecteur que le fait ou l'argument \'evoqu\'e par ce mot est justifi\'e
par une Blitzbeweis, ou une r\'ef\'erence,  dans la $n^{\hbox{\scriptsize i\`eme}}$ des notes,  en bout de texte,
 p. {\uppercase\expandafter{\romannumeral 3}} \`a {\uppercase\expandafter{\romannumeral 6}}.

\vfill\eject
\null\vfill
\setcounter{page}{1}
\pagenumbering{arabic}
\TrimTop{-7pct}
\TrimBottom{-3pct}
\centerline{\BoxedEPSF{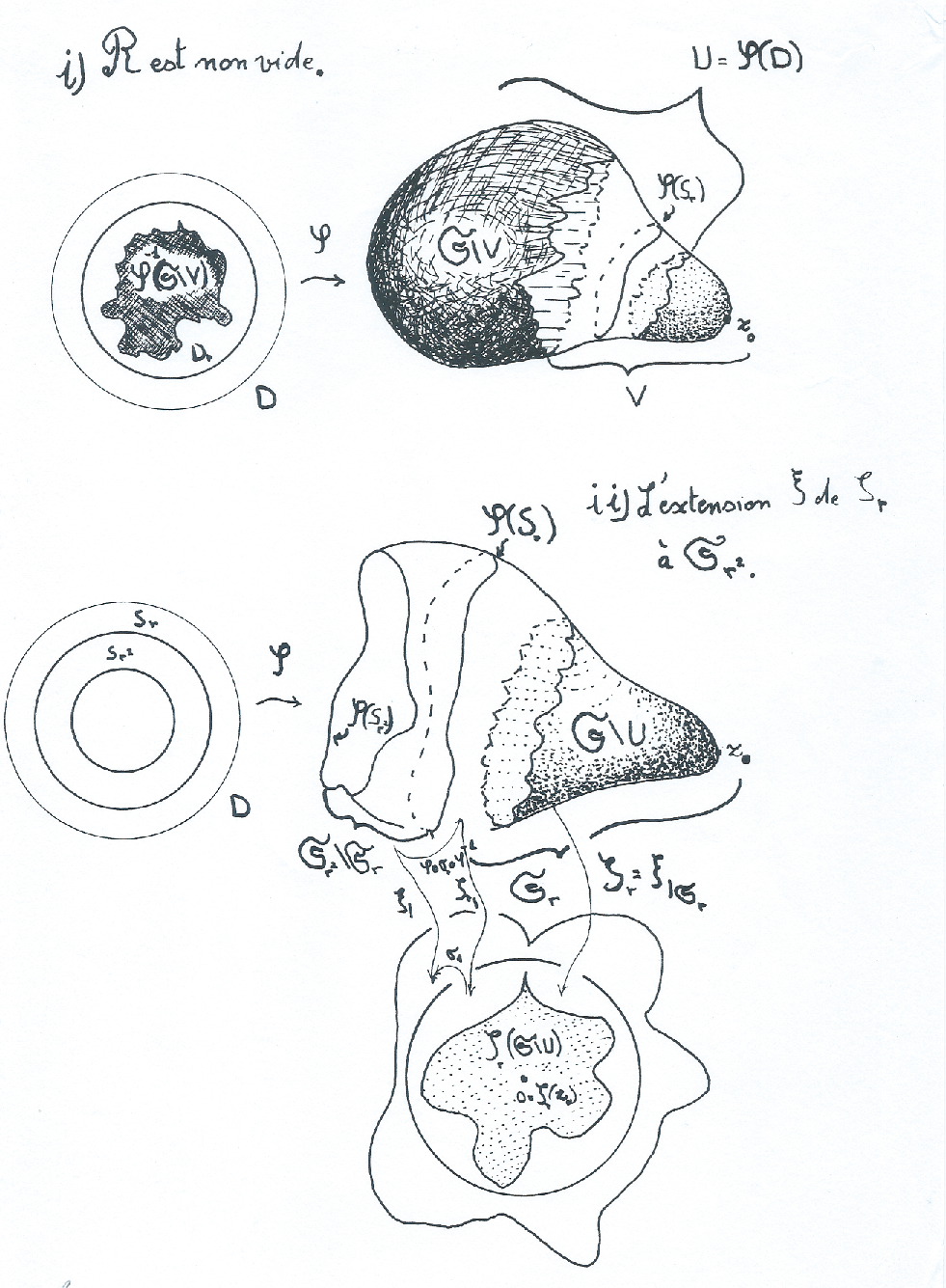 scaled 500} }
\centerline{Planche a}
\vskip2mm
\centerline{D\'emonstration du Lemme A}
\vfill
\null\vfill\eject
%

%
\null\vskip5mm
\fancyhead[LO]{A. M. \quad Uniformisation des surfaces de Riemann}
  \fancyhead[RO]{\thepage}
  \fancyhead[RE]{\hfill\small D\'emonstration du Lemme A\hfill\quad}
  \fancyhead[LE]{\thepage}
\centerline{\Large D\'eroulement de la preuve}
\vskip5mm
{\parindent=0pt\par\vskip .3cm
\vskip 0mm plus -20mm minus 1,5mm\penalty-50
{\bf 1\ \/}%
{\bf D\'emonstration du  Lemme A}{\bf \pointir}%
\nobreak\parindent=20pt}%
Soient $U\/$ et $V\/$ les deux ouverts
\'el\'emen\-taires qui
recouvrent ${\goth S\/}\/$ et $\varphi : D\rightarrow U\/$
un iso\-mor\-phisme. Puisque $U\/$ n'est pas compact, mais ${\goth S\/}\/$ l'est,
il y a  un point $z_0\/$ hors de $U\/$.

Pour $r\in ]0, 1[\/$, comme ${\goth S\/}\/$ est s\'epar\'ee, $\overline{D}_r\/$
compact et $\varphi\/$ continue, le
compl\'emen\-taire ${\goth S\/}_r={\goth S\/}\!\setminus\!\varphi\,(\overline{D}_r)\/$ de
$\varphi\,(\overline{D}_r)\/$ est ouvert dans ${\goth S\/}\/$.

Soit  ${\cal R\/}$ l'ensemble des $r\!\in]0, 1[\/$
tels qu'il y a un isomorphisme $\psi_r\/$ de ${\goth S\/}_r\/$ sur $D\/$.
Si $r\in{\cal R\/}$ l'isomorphisme $\zeta_r=\varphi_{\psi_r\,(z_0)}\circ \psi_r: {\goth S\/}_r\rightarrow D\/$
v\'erifie $\zeta_r\,(z_0)=0\/$.

{\parindent=0pt\par\vskip .3cm
\vskip 0mm plus -20mm minus 1,5mm\penalty-50
i)\  {\sl ${\cal R\/}\/$ est non vide\/}{\bf \pointir}%
\nobreak\parindent=20pt}%
 En effet le compact ${\goth S\/}\!\setminus\!V\/$
 de $U=\varphi\,(D)\/$ est inclus dans $\varphi\,(D_r)\/$
pour $r\/$ assez proche
de $1\/$. En ce cas, $\varphi\/$ \'etant$^{6}\/$ ouverte,
 $\overline{\goth S\/}_r={\goth S\/}\!\setminus\!\varphi\,(D_r)\/$
est une sous-vari\'et\'e compacte de $V\/$ d'int\'erieur ${\goth S\/}_r\/$
et de bord, \'egal \`a sa fronti\`ere 
${\rm Fr\/}\,{\goth S\/}_r={\rm Fr\/}\,\varphi\,(D_r)=\varphi\,(S_r)\/$, connexe.
 Ainsi l'$\omega\/$-polygone
${\goth S\/}_r\subset V\simeq D\subset{\Bbb C\/}\/$,
isomorphe \`a un domaine
de Jordan, est \'el\'ementaire par le
{\petcap Corollaire 0\/}.\fcarre

{\parindent=0pt\par\vskip .3cm
\vskip 0mm plus -20mm minus 1,5mm\penalty-50
ii)\ {\sl ${\cal R\/}\/$
est stable par $r\!\!\mapsto\!\!r^2\/$\/}{\bf \pointir}%
\nobreak\parindent=20pt}%
L'isomorphisme $\zeta_r : {\goth S\/}_r\rightarrow D\/$ a, par le principe de
sym\'etrie$^{7,\,7'}\/$ de Schwarz,
une extension holomorphe injective $\xi : {\goth S\/}_{r^2}\rightarrow {\Bbb C\/}\/$.

Son image $\xi\,({\goth S\/}_{r^2})\/$
est$^{6}\/$ ouverte,
distincte du plan complexe ${\Bbb C\/}$ puisque,
sinon, le th\'eor\`eme$^{8,\,8'}\/$ des singularit\'es inexistantes
\'etendrait
\`a la sph\`ere de Riemann $P_1\,({\Bbb C\/})={\Bbb C\/}\cup\{\infty\}\/$
l'inverse $\xi^{-1}:{\Bbb C\/}\rightarrow {\goth S\/}_{r^2}\subset {\goth S\/}\/$ de $\xi\/$
en une application
holomorphe, qui n'\'etant pas constante est d'image dans ${\goth S\/}\/$
\`a la fois compacte et$^{6}\/$ ouverte, contredisant que
$\overline{\goth S\/}_{r^2}={\goth S\/}\!\setminus\!\varphi\,(D_{r^2})\/$ n'est pas ouvert.

L'ouvert $\xi\,({\goth S\/}_{r^2})\/$, union de $\xi\,({\goth S\/}_{r^2}\!\setminus\!\overline{\goth S\/}_r)=
\xi\circ\varphi\,(D_r\!\setminus\!\overline{D}_{r^2})\/$, 
homom\'eomorphe \`a un anneau semi-ferm\'e, et du disque ferm\'e
$\overline{D}=\overline{\xi\,({\goth S\/}_r)}\/$,
satisfait aussi$^{9}\/$ \`a la derni\`ere
hypoth\`ese
du {\petcap Th\'eor\`eme 0\/},
ainsi $\xi\,({\goth S\/}_{r^2})\/$
est isomorphe \`a $D\/$ et $r^2\/$ est dans ${\cal R\/}$.\fcarre

Le
crit\`ere de Montel {\bf 4.1\/}
donne, si $r\/$ est dans
${\cal R\/}$, un isomorphisme
$$\psi : {\goth S\/}\!\setminus\!\varphi^{-1}\,(0)=\build{\cup}_{n=1}^{\infty}{\goth S\/}_{r^{2^n}}%
\rightarrow W\subset{\Bbb C\/}\/$$
 de
${\goth S\/}\!\setminus\!\varphi^{-1}\,(0)\/$
sur un ouvert $W\/$ de ${\Bbb C\/}\/$ qui, d'apr\`es le th\'eor\`eme$^{8,\,8"}\/$ des
singularit\'es inexistantes, s'\'etend en 
$\Psi : {\goth S\/}\rightarrow P_1\,({\Bbb C\/})\/$ holomorphe injective, donc
isomorphisme puis\-que ${\goth S\/}\/$ est compacte et la sph\`ere de Riemann
$P_1\,({\Bbb C\/})\/$ est s\'epar\'ee et connexe.

\ffincdem
\vfill\eject
\null\vfill

\TrimTop{-7pct}
\TrimBottom{-3pct}
\centerline{\BoxedEPSF{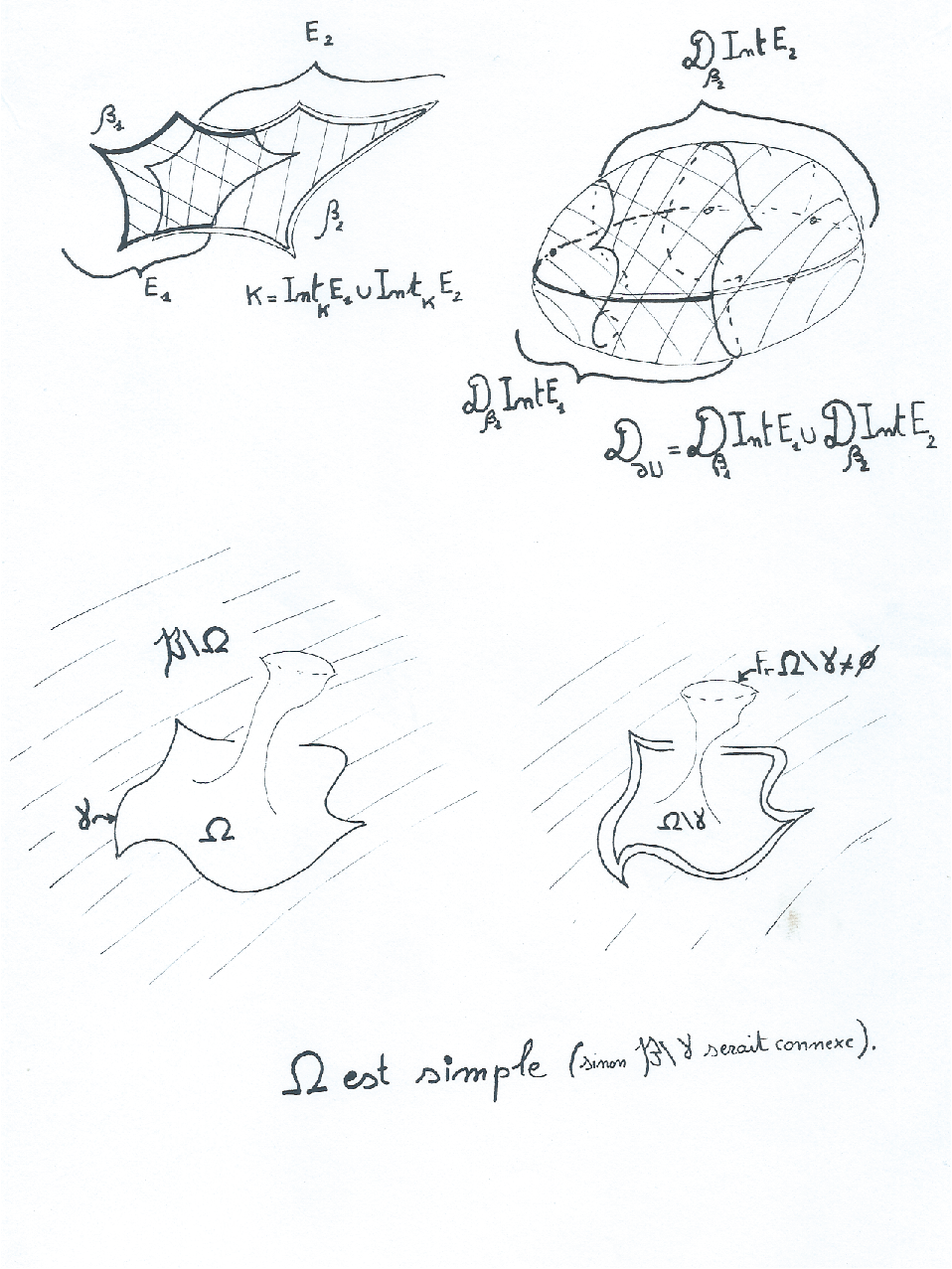 scaled 500} }
\centerline{Planche b}
\vskip2mm
\centerline{Le Lemme B et Affirmation 1}
\vfill
\vfill\eject
%
%
\fancyhead[LO]{A. M. \quad Uniformisation des surfaces de Riemann}
  \fancyhead[RO]{\thepage}
  \fancyhead[RE]{\S 2\hfill\small Le Lemme B et d\'emonstration du Lemme C\hfill\quad}
  \fancyhead[LE]{\thepage}
{\parindent=0pt\par\vskip .3cm
\vskip 0mm plus -20mm minus 1,5mm\penalty-50
{\bf 2\ \/}%
{\bf Le  Lemme B et d\'emonstration du  Lemme C\/}{\bf \pointir}%
\nobreak\parindent=20pt}%
Un {\it $\omega\/$-polygone ferm\'e\/} $F\/$
est la fermeture $F=\overline{U}=U\cup{\rm Fr\/}\,U\/$ d'un $\omega\/$-polygone ouvert $U\/$
(son {\it int\'erieur\/}).

Le {\it bord\/} $\partial F\/$
 d'un $\omega\/$-polygone ferm\'e $F\/$ est $\partial F=\partial\,{\rm Int\/}\,F\/$, celui de son int\'erieur.

Un $\omega\/$-polygone (ouvert ou ferm\'e) $X\/$ est {\it r\'egulier (\/{\rm resp. \/}simple)\/} si son bord est injectif
(resp. et d'image connexe). Le bord $\partial X\/$ s'identifie alors \`a son image ${\rm Fr\/}\,X\/$.

Une {\it cellule\/}
est un $\omega\/$-polygone compact simple d'int\'erieur \'el\'ementaire.

Une cellule $E\subset F\/$
est {\it  p\'eriph\'erique dans\/} un $\omega\/$-polygone ferm\'e $F\/$
la contenant si l'intersection
des images de leurs bords est un arc ferm\'e non r\'eduit \`a un point.
\Thc Lemme B| Un polygone analytique compact 
et simple $K\/$
union des int\'erieurs relatifs de deux
de ses cellules p\'eriph\'eriques $E_1\/$ et $E_2\/$ est une cellule.
\finc
\Dem Soit $U={\rm Int\/}\,K\/$ et 
$\beta_i=\partial {E_i}\!\setminus\!{\overline{U\cap\partial E_i}\/}\/$.
La surface de Riemann  ${\goth S\/}={\goth D\/}_{\partial U}\,U\/$ double de $U\/$ sur $\partial U\/$
(voir  l'appendice 3)
 est  compacte et
recouverte par les doubles
${\goth D\/}_{\beta_i}\,({\rm Int\/}\,E_i)\/$ des int\'erieurs
des cellules $E_i\/$ sur les arcs $\beta_i\/$, deux disques par
{\bf 3.6\/}. Ainsi ${\goth S\/}\/$ est,
d'apr\`es le {\petcap Lemme A\/}, isomorphe \`a 
$P_1\,({\Bbb C\/})\/$.

Toute h\'et\'erolution non libre de 
$P_1\,({\Bbb C\/})\/$ \'etant$^{10}\/$ conjugu\'ee
\`a l'inversion $\sigma_1\/$, l'int\'erieur
$U\/$ du polygone analytique compact $K\/$,
isomorphe \`a une des composantes du compl\'ementaire
${\goth D\/}_{\partial U}\,U\!\setminus\!{\rm Fix\/}\,\sigma_{{\partial U}}\/$,
est \'el\'ementaire, donc $K\/$ est une cellule.\ffindem
{\parindent=0pt\par{\sl D\'emonstration du {\petcap Lemme C\/}\/}\pointir\parindent=20pt}%
L'union du connexe ${\goth P\/}\!\setminus\!\Omega$ et
des composantes non ferm\'ees de $\Omega\/$, sauf au plus une, \'etant connexe \`a
compl\'ementaire connexe dans la composante de ${\goth P\/}\!\setminus\!\Omega\/$ dans ${\goth P\/}\/$, 
la preuve pour $\Omega\/$ connexe suffit. En ce cas :
\Thc Affirmation 1| Le polygone analytique relativement compact $\Omega\/$ est simple.
\finc
\Dem Sinon  sa fronti\`ere
contiendrait strictement une courbe sim\-ple ferm\'ee $\gamma\/$ de compl\'ementaire
${\goth P\/}\!\setminus\!\gamma=
(\Omega\!\cup\!({\rm Fr\/}\,\Omega\!\setminus\gamma))\cup
(({\rm Fr\/}\,\Omega\!\setminus\gamma)\!\cup\!({\goth P\/}\!\setminus\overline{\Omega}))\/$
con\-nexe, car union de deux connexes non disjoints, contredisant la planarit\'e de ${\goth P\/}\/$.
\ffindem%

Si l'une des composantes de l'un des ouverts $U_k\/$ rencontre $\Omega\/$ et est compacte,
elle serait composante connexe isomorphe \`a $P_1\,({\Bbb C\/})\/$ de ${\goth P\/}\/$ contenant 
le $\omega\/$-poly\-gone simple $\Omega\/$. Ce dernier serait donc, par le
{\petcap Corollaire 0\/},  standard.\fcarre

Sinon chaque composante des $U_k\/$ est
\'el\'ementaire, et :
\Thc Affirmation 2| Il y a des cellules
$E_{k, h},\ k=1,\, 2,\ h=1,\,\ldots,\,n_k\/$ telles que :

{\it (i)\/} Les int\'erieurs $U_{k,\,h}\/$ des cellules $E_{k,\,h}\/$ recouvrent $\overline{\Omega}\/$.

{\it (ii)\/} \`A $k\/$ fix\'e les cellules $E_{k,\,h}\/$ sont deux \`a deux disjointes.

{\it (iii)\/} Les composantes de $\partial E_{k, h}\cap\Omega\/$ 
 sont 
$\alpha_{k,\,l},\ l=1,\,\ldots,\, m_k\/$, en nombre fini,
et leurs fermetures $\beta_{k,\,l}=\overline{\alpha_{k,\,l}}\/$ sont deux \`a deux disjointes.

{\it (iv)\/} Aucune sous-famille stricte des $E_{k, h}\/$ n'a les propri\'et\'es
{\it (i)\/}, {\it (ii)\/} et {\it (iii)\/}.
\finc
\vfill\eject
\null\vfill

%
\TrimTop{-7pct}
\TrimBottom{-3pct}
\centerline{\BoxedEPSF{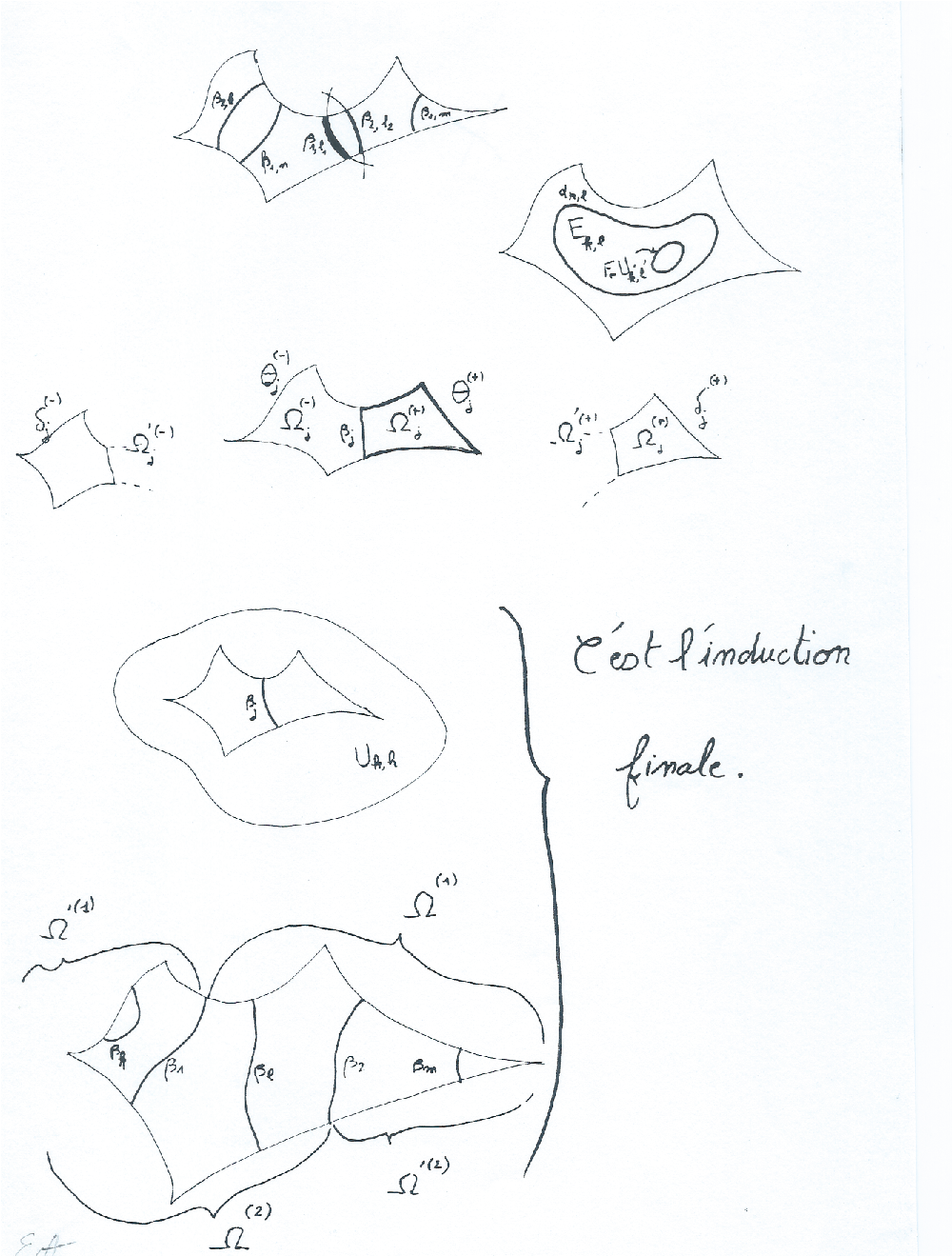 scaled 500} }
\centerline{Planche c}
\vskip2mm
\centerline{fin de d\'emonstration du Lemme C}
\vfill
\null\vfill\eject
\Dem Le compact $\overline{\Omega}\/$ est inclus dans l'union d'une famille finie,
dont aucune sous-famille ne le recouvre,
de composantes $U_{k, h}\/$  des $U_k\/$ munies
d'isomorphismes
$\varphi_{k, h} : D\rightarrow U_{k, h}\/$. Pour $r<1\/$
assez proche de $1\/$, les 
$\varphi_{k, h}\,({D_r})\/$ recouvrent $\overline{\Omega}\/$
et leurs bord $\varphi_{k, h}\,(S_r)\/$ ne passent par aucun sommet
de $\Omega\/$ et sont transverses aux arcs
analytiques dont ${\rm Fr\/}\,\Omega\/$ est l'union.

Les  cellules $E_{k,\,h}=\varphi_{k, h}\,(\overline{D}_r)\/$ satisfont donc \`a {\it (i)\/},
{\it (ii)\/} et {\it (iv)\/},
 les $\alpha_{k, l}\/$ sont en nombre fini et, \`a $k\/$ fix\'e,
leur fronti\`eres dans $\partial E_{k,\,h}=\partial U_{k,\,h}\/$ ainsi que les $\beta_{k,\, l}\/$
sont disjoints. Ainsi  {\it (iii)\/} est aussi v\'erifi\'ee, car $\beta_{k,\,l}\subset\overline{\Omega}\/$,
pour $k=1,\,2\/$ et tout $l\/$ mais
{\parskip=2pt plus 1pt minus 1pt
$$\beta_{1,\,l_1}\!\cap\beta_{2,\,l_2}\subset
\partial E_{1, h_1}\cap \partial E_{2, h_2}\subset {\goth P\/}\setminus\cup \varphi_{k, h}\,({D_r})
\subset {\goth P\/}\setminus\overline{\Omega}\/$$
\vskip-2mm\ffindem}
Si $\alpha_{k,\,l}\/$ est
 simple ferm\'ee, donc bord 
$\alpha_{k,\,l}=\partial E_{k, h}\/$ d'une des cellules, elle est incluse dans
un $U_{k'\!,\,h'}\/$ (o\`u $k'\!=\!3-k\/$). Or, par {\it (iv)\/}, $U_{k,\,h}\!\not\subset\!U_{k'\!,\,h'}\/$,
 et $U_{k,\,h}\/$
ren\-contre, donc contient, la fronti\`ere ${\rm Fr\/}\,U_{k'\!,\,h'}\/$. 
Ainsi $\Omega\/$, de fronti\`ere connexe, est inclus dans la composante compacte
de ${\goth P\/}\/$,
union des deux ouverts \'el\'emen\-taires $U_{k,\, h}\/$ et 
$U_{k'\!,\, h'}\/$. Par
 {\petcap Lemme A\/} et {\petcap Corollaire 0\/}, le $\omega\/$-polygone simple $\Omega\/$ est donc
 standard.\fcarre

Sinon chaque $\beta_{k,\,l}\/$ joint
deux points du bord $\partial \Omega\/$ de $\Omega\/$. Une induction sur le cardinal $m\/$
de la famille ${\cal B\/}=\{\beta_j\, ;\, j=1,\, \ldots,\, m\}\/$
 de ces arcs  conclura, gr\^ace \`a l' 
\Thc Affirmation 3| Soit $\theta_{j}^{(\eta)}\/$, pour
$\eta\in\mu_2\/$, les deux arcs  
joignant les extr\`emit\'es de $\beta_{j}\/$ dans le bord
$\partial \Omega\/$.
Chaque $\alpha_j\/$ s\'epare $\Omega\/$ en deux $\omega\/$-polygones
 connexe $\Omega_{j}^{(\eta)}\/$
 de bord les $\omega\/$-courbes simples ferm\'ees
$\delta_{j}^{(\eta)}\/$
union de l'arc $\beta_{j}\/$ et 
de $\theta_{j}^{(\eta)}\/$.  
\finc
\vskip-2mm
\Dem Selon  {\bf 1.3\/}{\it (ii)\/},
la $\omega\/$-courbe simple ferm\'ee $\delta\_{j}^{(\eta)}\!\!\!\!\/$,
disjointe du connexe
$X^{(\eta)}\!=\!({\goth P\/}\!\setminus\!\overline{\Omega})\cup{\breve{\theta}}_j^{(-\eta)}\!\!\!\!\!\!\!\/$,\ \ \
s\'epare le planaire ${\goth P\/}\/$ en deux $\omega\/$-polygones connexes :

L'un ${\Omega'}_{j}^{(\eta)}\supset X^{(\eta)}\/$
contient $X^{(\eta)}\/$ et l'autre $\Omega_{j}^{(\eta)}\subset\Omega\/$ est inclus dans $\Omega\/$.

Ainsi
$\Omega_{j}^{(\eta)}\subset{\Omega'}_{j}^{(-\eta)}\/$ et 
l'union 
$\Omega_{j}^{(-)}\cup{{\alpha}}_{j}\cup\Omega_{j}^{(+)}\!\!\!\!\/$,\quad
disjointe, ferm\'ee dans $\Omega\/$ et,
contenant un voisinage de ${{\alpha}}_{j}\/$,
aussi ouverte donc \'egale \`a $\Omega\/$ car $\Omega\/$ est connexe.\ffindem

\Indf Si $m\leq 1\/$ 
l'arc \'eventuel $\beta_{j}\/$ de 
${\cal B\/}\/$
est inclus dans un $U_{k,\,h}\/$ rencontrant $\Omega\/$,
mais 
$\partial U_{k,\,h}\/$ est disjoint de $\Omega\/$.
L'ouvert simple $\Omega\/$ est donc inclus dans $U_{k\,h}\/$
donc, par le {\petcap Corollaire 0\/}, \'el\'ementaire.\fcarre

Sinon ${\cal B\/}\/$ a deux arcs distincts $\beta_j,\ j=1,\,2\/$.
Soit $\Omega^{(j)}\/$
la composante de $\Omega\!\setminus\!\beta_j\/$ de fermeture
$F^{(j)}=\overline{\Omega}^{(j)}\/$
contenant $\beta_{3-j}\/$ et
$\Omega'^{(j)}\/$ l'autre composante.

Comme
${\goth P\/}\!\setminus\!\overline{\Omega}^{(j)}=
({\goth P\/}\!\setminus\!\overline{\Omega}\cup
(\partial\Omega\!\setminus\!\partial\Omega^{(j)}))\cup
((\partial\Omega\!\setminus\!\partial\Omega^{(j)})\cup\Omega'^{(j)})\/$ est connexe, et 
{\it (i)\/}
\`a {\it (iv)\/} de l'{\petcap Affirmation 2\/} pour $\Omega^{(j)}\/$
sont satisfaites par une partie de la famille des
$E_{k,\, h}\/$ avec un nombre d'arcs $m^{(j)}<m\/$
inf\'erieur,
 $F^{(j)}\/$ est une cellule.

Les int\'erieurs relatifs de deux cellules $F^{(j)}\/$ recouvrant
le compact r\'egulier
$\overline{\Omega}\/$ et y \'etant p\'eriph\'eriques,
$\overline{\Omega}\/$ est, par le {\petcap Lemme B\/}, 
une cellule. Ainsi  $\Omega\/$
est \'el\'ementaire.

\ffincdem
\vfill\eject
\null\vfill

\TrimTop{-7pct}
\TrimBottom{-3pct}
\centerline{\BoxedEPSF{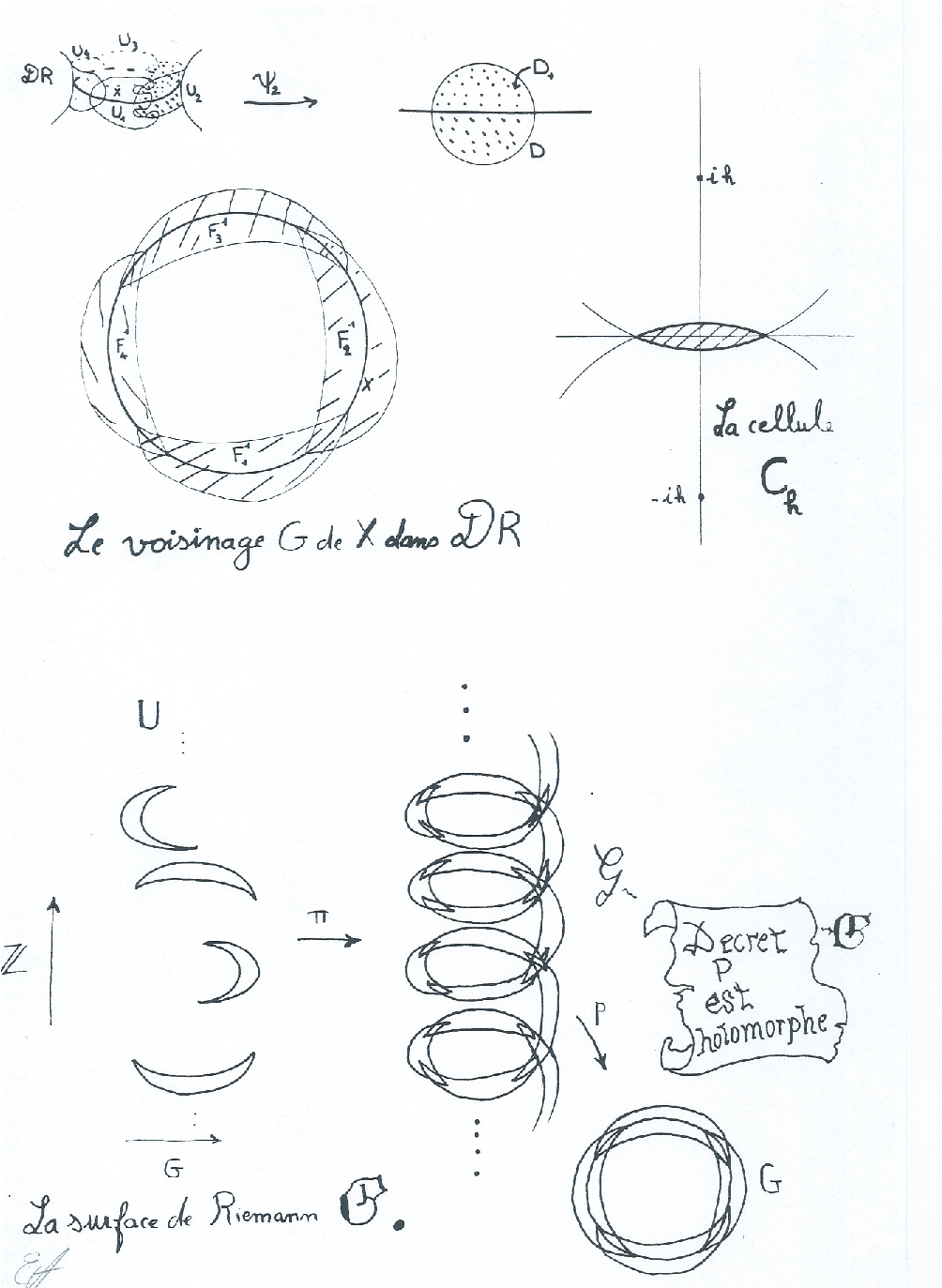 scaled 500} }
\centerline{Planche d}
\vskip2mm
\centerline{D\'emonstration du Lemme D}
\vfill
\null
\vfill\eject
%
\fancyhead[LO]{A. M. \quad Uniformisation des surfaces de Riemann}
  \fancyhead[RO]{\thepage}
  \fancyhead[RE]{\S 3\hfill\small D\'emonstration des Lemme D et Corollaire D\hfill\quad}
  \fancyhead[LE]{\thepage}
{\parindent=0pt\par\vskip .3cm
\vskip 0mm plus -20mm minus 1,5mm\penalty-50
{\bf 3\ \/}%
{\bf D\'emonstration des Lemme D et Corollaire D\/}{\bf \pointir}%
\nobreak\parindent=20pt}%
La composante compacte $X\/$ de $\partial{\bf R\/}\/$ est recouverte par les int\'erieurs $U_k\/$
de $\!N\!\/$ cellules $E_k$ du double $\!{\goth D\/}{\bf R\/}\/$, invariantes par $\sigma_{\bf R\/}\/$, avec de plus : \quad
$(1)\/$ $N-1$ des $U_k\/$ ne recouvrent pas $X\/$.

En les prenant suffisamment petites, on a en outre $N>2\/$ et :

$(2)\/$ Si deux de ces cellules s'intersectent, elles sont incluses dans une carte.

Pour $k=1,\,\ldots,\,N\/$
soit
$\psi_k:(\overline{D},\,D,\, D_+;\, {\rm conj\/})\longrightarrow
(E_k,\,U_k,\,U_k\cap {\bf R\/};\,{\sigma_{\bf R}}_|)\/$,
 iso\-mor\-phimes de triple conjuguant
${\rm conj\/}$ aux restrictions de $\sigma_{\bf R\/}\/$, et $\psi_{N+1}=\psi_1\/$.
Quitte \`a remplacer $E_k\/$ par $\psi_k(\overline{D}_r)\/$
[et $\psi_k\/$ par $\psi_k\!\circ\!h_r\/$], pour $r\!<\!1\/$ assez 
 proche de $1\/$,

$(3)$ Les $\psi_k\/$ ont des extensions holomorphes injectives d\'efinies sur des voisinages de $\overline{D}\/$ et les intervalles $I_k\!=\!E_k\!\cap\!X\!=\!\psi_k\,([-1, 1])\/$
sont, pour $k\!\leq\!N\/$ de bords disjoints.

  Alors, $X\/$ \'etant connexe, {\bf 2.1\/}
renum\'erote les $I_k\/$ de sorte qu'orient\'es par les $\psi_k\/$ et avec
$I_{N+1}=I_{1}\/$, les intersections $I_k\cap I_{k+1}\/$ sont,  pour $k\leq N\/$, des sous-intervalles ferm\'es
 stricts de $I_k\/$ et de
$I_{k+1}\/$, de bord l'extr\'emit\'e finale de $I_k\/$ et l'origine de $I_{k+1}\/$,
les autres intersections $I_k\cap I_l=\emptyset\/$ \'etant vides.

La cellule
$C_{h}=\tau_{-1}\circ h_2\,(\overline{P}_{0,\,h})=\{z\in\overline{D}\  ;
\,
|z+ \eta\,i\, h|^2\leq 1+h^2 \, \hbox{\small pour \/}\,  \eta\in\mu_2\}\/$, est un voisinage de $]-1,\,1[\/$
 de bord  deux arcs lisses et
transverses \`a ${\Bbb R\/}$,
mais tendant vers $[-1,\, 1]\/$ dans la topologie $C^1\/$ 
quand la finesse $h\/$
tend vers l'infini.

Par r\'ecurrence 
sur $k=1,\ldots,\, N\/$, il y a
$h_k>0\/$, tels que,
si $t\geq 1,\, k\leq N\/$, les  cellules
$F^t_k=\psi_k\,(C_{t\,h_k})\/$ et $F^t_{N+1}=F^t_1\/$ v\'erifient :

{\it (i)\/} L'intersection $\partial F^t_k\cap \partial F^t_{k+1}\/$ est transverse et r\'eduite \`a deux points.

{\it (ii)\/} Si $1<|k-l|<N-1\/$, les cellules $F^t_k\/$ et $F^t_l\/$ sont disjointes.

L'union $G=G_1\cup\cdots\cup G_N\/$, o\`u $G_k={\rm Int\/}\,F^1_k\/$ est un
ouvert, $\sigma\/$- ({\it i.e.\/} invariant par $\sigma={\sigma_{\bf R}}_{|G}\/$)-voisinage
de $X\/$ dans  de ${\goth D\/}\,{\bf R\/}\/$. Notant $\bar n\/$ le repr\'esentant
dans $\{1,\,\ldots,\,N\}$
de la classe de l'entier $n\in{\Bbb Z\/}\/$ modulo $N\/$,
soit dans $G\times{\Bbb Z\/}\/$ l'ouvert%
$$U=\{(x,\,n)\in G\times{\Bbb Z\/}\,;\, x\in G_{\bar n}\}$$
{et $\tau_U,\, \sigma_U :U\rightarrow U\/$ les restrictions \`a $U\/$
de $(x,\, n)\mapsto(x,\, n+N)\/$ et
$\sigma\times{\rm Id\/}_{\Bbb Z\/}\/$.}

Sur $U\/$ est d\'efinie l'\'equivalence $\sim\/$ 
 par
$(x,\, n)\sim(y,\, m)\/$ si et seu\-lement si%
$$x=y\ \hbox{\rm et \/} |n-m|\leq 1$$

La restriction $P\/$ \`a $U\/$ de 
${\rm pr\/}_1 : G\times{\Bbb Z\/}\rightarrow G\/$
se factorise $P=p\circ\pi\/$
en $\pi\/$, application quotient 
 de $U\/$ sur l'espace quotient
${\cal G\/}=U/_{\sim}\/$, et la {\it projection\/} $p\/$.

L'espace ${\cal G\/}\/$ est, car$^{12}$ $N\geq3\/$, s\'epar\'e.
La relation $\sim\/$ \'etant ouverte et la com\-po\-sante $U_n=G_{\bar n}\times\{n\}\/$ de $P^{-1}\,(G_{\bar n})\/$
hom\'eomorphe par $P\/$ \`a $G_{\bar n}\/$, les ouverts $U_n\/$ de $U\/$
sont hom\'eomorphes, par la restriction de $\pi\/$, aux
ouverts
${\cal G\/}_n=\pi\,(U_n)\/$ de ${\cal G\/}\/$.

En d\'ecr\'etant  $p\/$ holomorphe, cette surface
s\'epar\'ee ${\cal G\/}$ devient, une surface de Riemann ${\goth G\/}$,
munie de l'h\'et\'erolution $\sigma_{\goth G\/}\/$
quotient de $\sigma_U\/$ et
de l'action holomorphe propre libre et d'orbites les pr\'eimages de $p\/$,
du groupe cyclique infini engendr\'e par l'isomorphisme $\tau_{\goth G\/}\/$ induit par $\tau_U\/$. Ainsi $p\/$
descend en $\overline{p} : {\goth G\/}/_{<\tau_{\goth G\/}>}\rightarrow G\/$, isomor\-phisme \'equivariant
de la surface de Riemann quotient sur un $\sigma\/$-voisinage de $X\/$.
\vfill\eject
\null\vfill
\TrimTop{-7pct}
\TrimBottom{-3pct}
\centerline{\BoxedEPSF{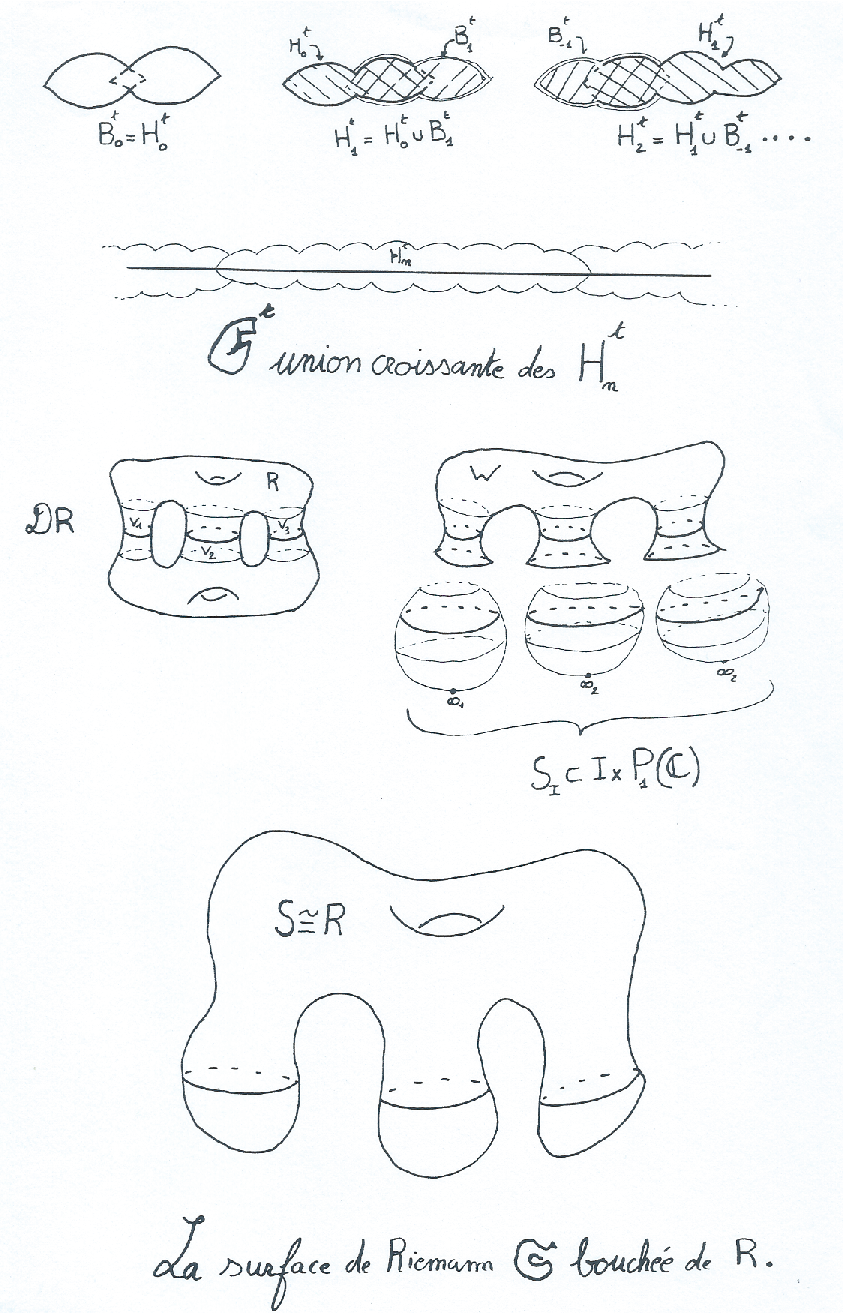 scaled 500} }
\centerline{Planche e}
\vskip2mm
\centerline{D\'emonstration du Corollaire D}
\vfill
\null\vfill\eject

Soit $q_n\/$ l'inverse de l'isomorphisme
de ${\goth G\/}_n\cup{\goth G\/}_{n+1}\/$
sur $G_{\bar n}\cup G_{\overline{n+1}}\/$ que $p\/$ induit.  
Pour $t>1\/$
l'union $K^t_k=F^t_k\cup F^t_{k+1}\/$, $\omega\/$-polygone ({\bf 3.1\/}{\it (iv)\/}) compact  simple de $G\/$,
est par $(2)\/$ et {\petcap Corollaire 0\/}, une cellule, d'image par $q_n\/$
 not\'ee $B^t_n=q_n\,(K^t_{\bar n})\/$.

Si $t\!>\!1\/$, la suite $H^t_{2\,m}=\build{\cup}_{j=-m}^{m}B^t_j,\
H^t_{2\,m+1}=\build{\cup}_{j=-m}^{m+1}B^t_j\/$,
de $\omega\/$-polygones ({\bf 3.1{\it (iv)\/}\/})
 compacts,  
d\'ebute par la cellule
$H^t_0\/$, puis $H^t_{p+1}\/$ est union des int\'erieurs
relatifs, de $H^t_{p}\/$ coupant $\partial H^t_{p+1}\/$ en
un arc, et de la cellule p\'eriph\'erique $B^t_{q}\/$ o\`u $q=(-1)^p\,[\frac{p+2}{2}]\/$.

Par {\petcap Lemme B\/}, les $H^t_{p}\/$ sont de proche en proche des cellules.
L'union de leurs int\'e\-rieurs ${\goth G\/}^t\/$,
est $\sigma\/$-in\-vari\-ante, et
\'el\'ementaire par Montel {\bf 4.1\/} et Liouville$^{11'}\/$.

L'image $G^t=p\,({\goth G\/}^t)\/$, 
voisinage $\sigma_{\bf R\/}\/$-invariant de $X\/$ dans ${\goth D\/}\,R\/$, 
est alors \'equivariament$^{13'}\/$ iso\-mor\-phe au quotient de $D\/$ par une$^{13}\/$ homographie r\'eelle $\theta\/$
libre de $D\/$, \'equivariament$^{13''}\/$
conjugu\'e
 \`a une l'homoth\'etie
$h_t\/$ de rapport $t>1\/$ de $L_+\/$.

Il suffit de remarquer que, si $r=e^{-\frac{\pi^2}{\log\,t}}\/$, le
quotient ($L_+/_{h_{t}},\, {\rm conj\/})\/$
est, par $\exp\circ h_{\frac{2\,\pi\,i}{\log t}}\circ\log\/$,
\'equivari\-ament isomorphe 
\`a l'anneau invariant $(A_r,\, \sigma_1)\/$ de rayon int\'erieur $r\/$,
un voisinage invariant de $S^1\/$ dans
$(P_1\,({\Bbb C\/}),\,{\rm conj\/})\simeq({\goth D\/}\,\overline{D},\,\sigma_{\overline{D}})\/$.
\ffindem

{\parindent=0pt\par{\sl D\'emonstration du {\petcap Corollaire D\/}\/}\pointir\parindent=20pt}%
 Soit
$X_i\/$, pour $i\/$ dans un ensemble fini $I\/$,
 les composantes de $\partial {\bf R\/}\/$. Le {\petcap Lemme D\/} donne,
sur des voisinages
$V_i\/$ dans 
${\goth D\/}\,{\bf R\/}\/$, deux \`a deux disjoints
des $X_i\/$, des
 isomorphisme $\phi_i : V_i\rightarrow A_{r_i}\/$,
avec $\phi_i\,(V_i\cap {\bf R\/})=A_{r_i}\cap\overline{D}\/$ et \'equivariants.
Soit $V=\cup_{i\in I}V_i\/$, alors  ${\bf R\/}\/$ a
$W=V\cup {\bf R\/}\/$ pour voisinage dans ${\goth D\/}{\bf R\/}\/$.
Dans $I\times P_1\,({\Bbb C\/})\/$, soit l'ouvert $S_I=\{(i,\,z)\in I\times P_1\,({\Bbb C\/})\,;|z|>r_i\}\/$.

L'application quotient de l'\'equivalence $\sim\/$ sur $S=W\amalg S_I$
d\'efinie par :
$$W\supset V_i\ni p\sim (i,\,\phi_i\,(p))\in \{i\}\times A_{r_i}\subset S_I\/$$
est not\'ee $\pi :  S\rightarrow{\cal S\/}\/$. Comme $\sim\/$ est triviale sur $W\/$ et sur $S_I\/$,
et les $\phi_i\/$ holomorphes, le quotient est une vari\'et\'e holomorphe. Elle est s\'epar\'ee car un point
de $\pi\,(W)\/$ et un point de $\pi\,(S_I)\/$, non tous deux dans l'un  de
ces deux ouverts s\'epar\'es, sont respec\-tivement dans $\pi\,({\bf R\/}\!\setminus\!\partial{\bf R\/})\/$
et $\pi\,(I\times(P_1\,({\Bbb C\/})\!\setminus\!\overline{D}))\/$, deux ouverts disjoints.

Cette surface de Riemann $ {\goth S\/}\/$ contient, ${\bf S\/}=\pi\,({\bf R\/})\/$, sous-surface de Riemann
 iso\-mor\-phe, par $\pi\/$ \`a ${\bf R\/}\/$. L'ouvert compl\'ementaire ${\goth S\/}\!\setminus\!{\bf S\/}\/$
 est standard car ses composantes sont isomorphes \`a $D\/$, par les
$(\pi\circ(c_{i}\times\imath))^{-1} : T_i=\pi\,(\{i\}\times P_1\,({\Bbb C\/})\!\setminus\!\overline{D})\rightarrow D \/$.\ffindem
\Thc Cocorollaire D| Si la surface de Riemann compacte \`a bord ${\bf R\/}\/$ est d'int\'erieur
 $U={\bf R\/}\!\setminus\!\partial {\bf R\/}$ planaire alors la surface de Riemann
compacte ${\goth S\/}$ construite par le {\petcap Corollaire D\/} est planaire.
\finc
\Dem Soit ${T}_i,\ i\in I\/$ les composantes de ${\goth S\/}\!\setminus\!{\bf S\/}\/$.
Pour tout choix de voisi\-nages $\omega_i\/$ dans $({\goth S\/}\setminus{\bf S\/})\cup V\/$
de points $p_i\in\overline{T}_i\/$
il y a un hom\'eomorphisme $H\/$ analytique par morceaux$^{14}\/$ de ${\goth S\/}$
\`a support dans $({\goth S\/}\!\setminus\!{\bf S\/})\cup V\/$ avec
$H\,({\goth S\/}\setminus(\cup\omega_i))\subset{\bf S\/}\setminus\partial{\bf S\/}\/$.

Ainsi tout $\omega\/$-cycle compact $\gamma\/$ de ${\goth S\/}\/$
s\'epare ${\goth S\/}\/$, puisque, 
en choisissant les $\omega_i\/$ disjoints de $\gamma\/$, le $\omega\/$-cycle compact
$\beta=H\,(\gamma)\/$ s\'epare l'ouvert planaire $U\/$.
\ffindem
\vfill\eject
%
%
{\parindent=0pt\par\vskip .3cm
\vskip 0mm plus -20mm minus 1,5mm\penalty-50
{\bf 4\ \/}%
{\bf D\'emonstration des th\'eor\`emes\/}%
\nobreak\parindent=20pt}%
{\parindent=0pt\par{\sl D\'emonstration du {\petcap Th\'eor\`eme 1\/}\/}\pointir\parindent=20pt}%
La surface de Riemann ${\goth P\/}\/$ \'etant compacte,
un nombre
fini d'isomorphismes $\psi_i: D\rightarrow U_i,\, i=0,\,\ldots,\, N\/$
de $D\/$ sont d'images recouvrant ${\goth P\/}\/$. Soit $r<1\/$ tel que les $V_i=\psi_i\,(D_r)\/$
couvrent encore ${\goth P\/}\/$ et $r_n\!=\!\frac{r+n}{1+n}\/$.

Comme ${\goth P\/}\/$ est connexe les $V_i\/$
se renum\'erotent dans $I=\{1,\,\ldots,\,N\}\/$ de sorte que si $0< i\leq N\/$ l'intersection
$V_{i}\cap(\cup_{j< i}V_j) \ne\emptyset\/$ est non vide.

Si $i\in I\/$ et $n\/$ entier,
la cellule $E_i^n\!=\!\psi_i\,(\overline{D}_{r_n})\/$
contient $V_i\/$ et $V_{i}^{n}\!=\!\psi_i\,({D}_{r_n})\/$. Ainsi, d'apr\`es {\bf 3.1\/}{\it (iv)\/}
 pour $p\!\leq\!N\/$, l'union
$F_p\!=\!\cup_{i<N-p}E_i^p\/$ est un $\omega\/$-polygone ferm\'e d'int\'erieur connexe
et la fermeture
$\overline{\Omega}_p\/$ de l'\/$\omega\/$-polygone compl\'ementaire $\Omega_p\!=\!{\goth P\/}\!\setminus\!F_p\/$
est  incluse dans l'union de $V_{N-p}^{p}\/$ et $\Omega_{p-1}\/$
(o\`u par convention $\Omega_{\,-1}\!=\!V_N\/$).

Les $V_i^{p}\/$ \'etant standards,
une application
inductive du {\petcap Lemme C\/} r\'ev\`ele comme standards tous les $\Omega_p\/$. 
En particulier la surface de Riemann ${\goth P\/}=\Omega_N\/$ est standard,
\'etant connexe et compacte,
elle est isomorphe \`a $P_1\,({\Bbb C\/})\/$.\ffindem
{\parindent=0pt\par{\sl D\'emonstration du {\petcap Th\'eor\`eme 2\/}\/}\pointir\parindent=20pt}%
Tout compact de ${\goth Q\/}\/$ se couvre de l'int\'erieur d'une union finie de cellules,
  un $\omega\/$-polygone (par {\bf 3.1\/}{\it (iv)\/}) relativement compact.

Ainsi la famille
${\cal URC\/}=(U_i)_{i\in I}\/$ des $\omega\/$-polygones relativement compacts de ${\goth Q\/}\/$
 est filtrante
et recouvre ${\goth Q\/}\/$.

Chaque $U_i\/$ de ${\cal URC\/}\/$ est planaire$^{15}\!\!\!\/$, isomorphe \`a l'int\'erieur de son arrondie,
il est, d'apr\`es le {\petcap Cocorollaire D\/}
isomorphe \`a un ouvert $V_i\/$ d'une surface de Riemann compacte planaire
${\goth S\/}_i$ dont chaque composante
${\goth S\/}_{i,\,j}\/$ de ${\goth S\/}_i\/$ est,
d'apr\`es le {\petcap Th\'eor\`eme 1\/}, isomorphe \`a $P_1\,({\Bbb C\/})\/$.

Si ${\goth Q\/}\/$ est compact le {\petcap Th\'eor\`eme 1\/} conclut.

Sinon 
aucune des ${\goth S\/}_{i,\,j},\ j=1,\ldots,n_i\/$
n'est incluse dans l'ouvert $V_i\/$.

Isomorphe \`a une union disjointe de $n_i\/$ ouverts de disques,
compl\'ementaires dans ${\goth S\/}_{i,\,j}\simeq P_1\,({\Bbb C\/})\/$ d'un disque ferm\'e disjoint
de l'ouvert $V_i\/$, 
chaque ouvert $U_i\/$ est donc une carte de ${\goth Q\/}\/$.
Ainsi, d'apr\`es le crit\`ere de Montel {\bf 4.1\/},
la surface de Riemann connexe ${\goth Q\/}\/$ est isomorphe \`a un ouvert de $P_1\,({\Bbb C\/})\/$.
\ffindem
\vfill\eject
%
%
%
\fancyhead[LO]{A. M. \quad Uniformisation des surfaces de Riemann}
  \fancyhead[RO]{\thepage}
  \fancyhead[RE]{\small\quad\quad Appendice 1\hfill Arcs analytiques,   
duplicatas et s\'eparation par les $\omega\/$-graphes\hfill\quad}
  \fancyhead[LE]{\thepage}
\centerline{\Large Appendices}
{\parindent=0pt\par\vskip .3cm
\vskip 2mm plus -20mm minus 1,5mm\penalty-50
{\bf 1\ \/}%
{\bf Arcs analytiques,  duplicatas et s\'eparation par les $\omega\/$-graphes\/}
\nobreak\parindent=20pt}%

Une {\it carte\/} (en un point $p\in{\goth S\/}\/$) d'une surface de Riemann ${\goth S\/}\/$ est
un ouvert $C\!\subset\!{\goth S\/}\/$ isomorphe par $\zeta : C\rightarrow V\/$, dit
{\it coordonn\'ee (en $p\/$, centr\'ee {\rm si\/} $\zeta\,(p)=0\/$) de $C\/$\/},
\`a un ouvert $V\/$ de ${\Bbb C\/}\/$. L'inverse $\gamma=\zeta^{-1}\/$ de $\zeta\/$ est 
l'\/{\it e\'ennodrooque (en $p\/$ associ\`ee \`a $\zeta\/$)\/}.

%
Un {\it arc analytique  {\rm [ou\/} $\omega\/$-arc\/}] (dans ${\goth S\/}\/$) 
{\it param\'etr\'e par $\gamma\/$, d'extr\'emit\'es, initiale {\rm ou\/} origine $p\/$ et finale $q\/$)\/},
est l'image  $\beta\/$ de l'intervalle $I\!=\![0,\,1]\/$ par une e\'ennodrooque $\gamma : V\!\rightarrow\!C\/$ d\'efinie sur
un voisinage $V\/$ de $I\/$ (avec $\gamma\,(0)\!=\!p\/$ et $\gamma\,(1)\!=\!q\/$), et  $\breve{I}\!=]0, 1[$%
.

L'\/{\it int\'erieur (formel) {\rm et le\/} bord\/}  de l'\/$\omega\/$-arc
$\beta=\gamma\,(I)\/$ sont ${\breve{\beta}} = \gamma\,({\breve{I}})\/$
et $\partial\beta=\beta\!\setminus\!{\breve{\beta}}\/$.
%
%

Un $\omega\/$-{\it graphe\/} de ${\goth S\/}\/$ est une famille ${\cal B\/}=(\beta_\lambda)_{\lambda\in \Lambda}\/$
localement finie,
de $\omega\/$-arcs dans ${\goth S\/}\/$, ses {\it ar\^etes\/},  
 avec, pour $\lambda\ne\mu,\ \beta_\lambda\!\cap\!\beta_\mu\/$
vide ou une extr\'emit\'e
commune.
%
%

La {\it valence\/} de $p\!\in\!{\goth S\/}\/$
dans ${\cal B\/}\/$  est 
$v_{\cal B\/}\,(p)={\rm card\/}\,\{\lambda\in\Lambda;\, p\ \hbox{\small extr\'emit\'e de \/} \beta_\lambda\}\/$.
%
%

L'union
 des ar\^etes $\beta_\lambda\/$ de ${\cal B\/}\/$, son {\it support\/} $\Gamma_{\cal B\/}\/$,
est partition\'e par les int\'e\-rieurs d'ar\^etes
${\breve{\beta}}_\lambda\/$ 
et son, discret, {\it ensemble $\Delta_{\cal B\/}=\{p\in \Gamma_{\cal B\/};\,v_p\ne 0\}\/$ de sommets\/}.

Une $\omega\/$-{\it cha{\^{\i}}ne\/} (resp. un $\omega\/$-{\it cycle\/}) de ${\goth S\/}\/$
est une union localement finie,
de $\omega\/$-arcs dans ${\goth S\/}\/$, 
(resp. o\`u tout point est extr\'emit\'e d'un nombre pair de ces arcs). Ainsi

{\sl L'union de deux $\omega\/$-cha{\^{\i}}nes (resp. cycles) est une $\omega\/$-cha{\^{\i}}ne (resp.
un $\omega\/$-cycle).\/}
%
%

Un $\omega\/$-graphe  est {\it ferm\'e\/} si tout $p\!\in\!{\goth S\/}\/$ est de valence $v_{\cal B\/}\,(p)\/$
 paire, $\Gamma_{\cal B\/}\/$  est alors un cycle
et (trivialement si $\Gamma=\emptyset\/$, par$^{16}$ sinon) r\'eciproquement :

{\sl Une cha{\^{\i}}ne 
$\Gamma\/$
 est support d'un graphe, qui est ferm\'e si $\Gamma\/$ est un cycle.\/}
%
%

Soit $\zeta\/$ une coordonn\'ee  centr\'ee
au sommet $s\!\in\!\Delta_{\cal B\/}\/$ et,
par locale finitude de ${\cal B\/}\/$, un voisinage $U_r\/$ de $s\/$
ne rencontrant que les ar\^etes 
$\beta_k,\, k=1,\ldots,\,v_{\cal B\/}\,(s)\/$ 
dont une extr\'e\-mit\'e est en $s\/$, n'en contenant aucune et avec
 $\zeta\,(U_r)\!=\!D_{r}\/$.

Le semi-arc $\alpha_k\!=\!\zeta\,(\beta'_k)\/$ image de la composante  $\beta'_k\/$
de $s\/$ dans $U_r\!\cap\!\beta_k\/$
est para\-m\'etr\'e par $f_k:[0,1[\rightarrow \alpha_k\/$ avec $f_k\,(0)=0\/$ et $f'_k\,(0)\ne0\/$. Ainsi,
pour $\epsilon>0\/$ petit, l'inter\-section $S_\epsilon\cap\zeta\,(\beta_k\cap U_r)\/$ est  transverse, r\'eduite
\`a  un unique point $z_k\!\in\!{\rm Im\/}\,f_k\/$ :
\Thnc {\bf 1.0\ \/}|D\'efemme|Soit $\Gamma\/$ une $\omega\/$-cha{\^{\i}}ne d'une surface de Riemann ${\goth S\/}\/$.

Alors tout $s\!\in\!{\goth S\/}\/$  est centre d'une e\'ennodrooque $\gamma\/$ dite
{\it $\Gamma\/$-normale\/},
{\it c.a.d.\/}\break d\'efinie pr\`es de ${\overline{D}}\/$ avec,
pour
$\beta'_{\gamma,k}$ des arcs initiaux  des  ar\^etes $\beta_{\gamma, k}$ issues de $s\/$
d'un graphe ${\cal B\/}\/$ dont $\Gamma\/$ est support,
$\gamma\,({\overline{D}})\cap\Gamma\!=\!
\bigcup_{k=1}^{v_{\cal B\/}\,(p)}\beta'_{\gamma,k}
\/$
 (en particulier $\gamma\,(\overline{D})\cap\Gamma=\emptyset\/$ si
$v_{\cal B\/}\,(s)=0\/$) et ${\rm card\/}\,(\beta_{\gamma,k}\cap\gamma\,(S_\rho))=1\/$
pour tout $\rho  \in]0,\,1]\/$ et $1\leq k\leq v_{\cal B\/}\,(p)\/$. 
\ffindem
\finnc

Si   ${\cal B\/}=(\beta_\lambda)_{\lambda\in\Lambda}\/$ est un $\omega\/$-graphe il y a $h_\lambda>0\/$
et $N_\lambda\/$ entier pair tel qu'une param\'etrisation $\gamma_\lambda\/$ de $\beta_\lambda\/$
est d\'efinie sur la perle $P_\lambda=P_{N_\lambda,\, h_\lambda}\/$ et les
images ${\goth p\/}_\lambda=\gamma_\lambda\,(P_\lambda)\/$ sont deux \`a deux disjointes$^{17}\/$ et couvrent
un voisinage ${\goth P\/}\/$ de
$\Gamma_{\cal B\/}\!\setminus\!\Delta_{\cal B\/}\/$.

\vfill\eject
\null\vfill

\TrimTop{-7pct}
\TrimBottom{-3pct}
\centerline{\BoxedEPSF{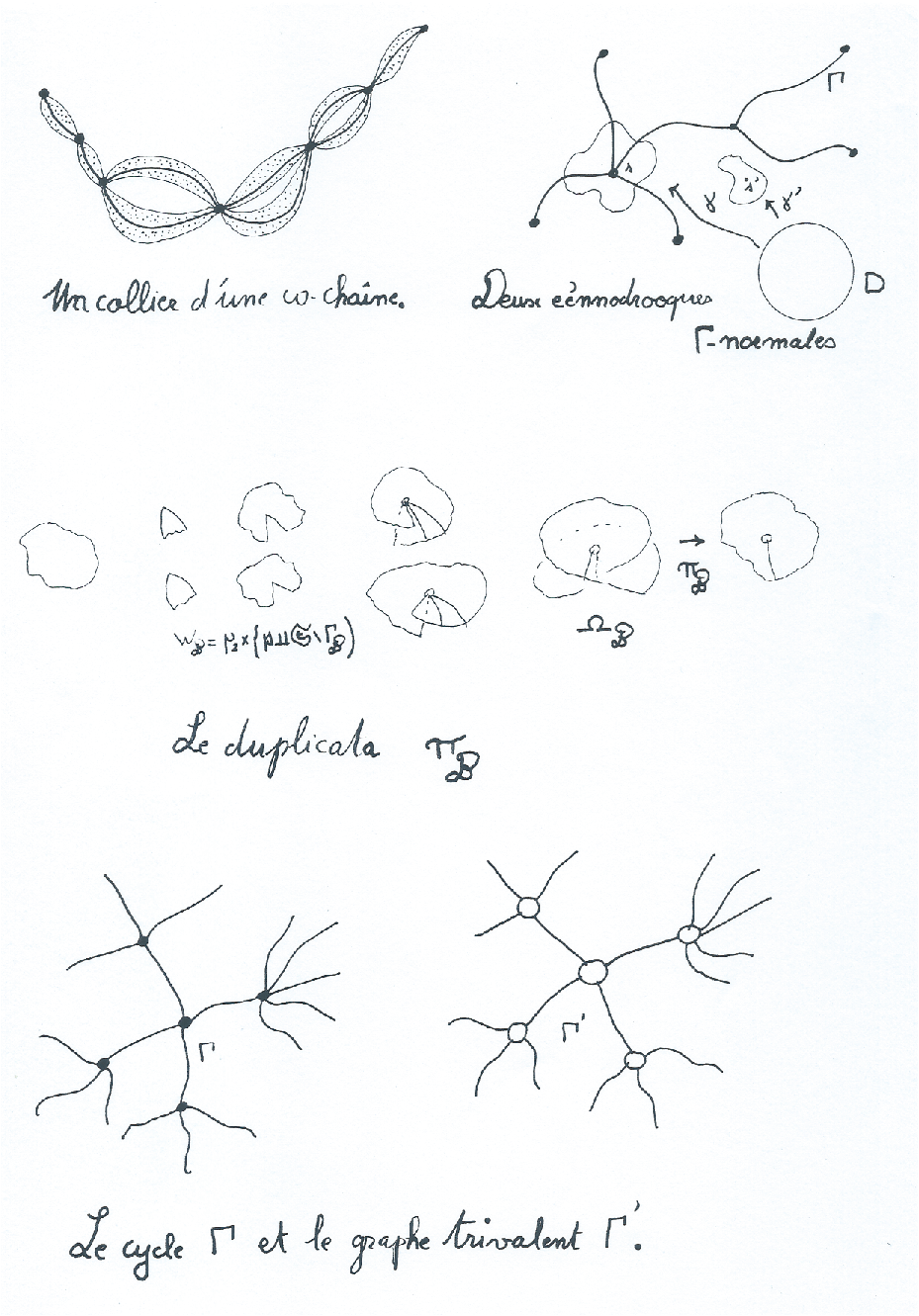 scaled 500} }
\centerline{Planche f}
\vskip2mm
\centerline{Colliers , e\'ennodroques, duplicata,\dots}
\vfill
\null\vfill\eject

Un {\it collier [param\'etr\'e] de\/}
la cha{\^{\i}}ne $\Gamma_{\cal B\/}\/$ 
de ${\cal B\/}$ (resp. -e {\it perle [-e] sur\/} l'ar\^ete $\beta_\lambda\/$)
est un tel voisinage ${\goth P\/}\/$
[$\gamma_{\goth P\/}\!=\!%
\amalg{\gamma_\lambda}_{|P_\lambda}\{(\lambda, z)\!\in\!\Lambda\!\times\!{\Bbb C\/};\,z\!\in\!P_\lambda\}\!%
\rightarrow\!{\goth P\/}\/$] (resp. ${\goth p\/}_\lambda\/$
[${\gamma_\lambda}_{|C_\lambda}\/$]).

%
%

L'application
$({\rm incl\/}_{\goth p\/}^{\goth S\/}\amalg{\rm incl\/}_{{\goth S\/}\!\setminus\!\Gamma}^{\goth S\/})\circ{\rm pr\/}_2\/$
induit, de $W_{\cal B\/}=\mu_2\!\times\!({\goth p\/}\amalg{\goth S\/}\!\setminus\!\Gamma_{\cal B\/})\/$,
sur
$$\Omega_{\cal B\/}=W_{\cal B\/}%
/\{(\eta,\,z)\sim (\eta\,\eta',\,z')\ 
\hbox{\small si \/} {\rm incl\/}_{\goth p\/}^{\goth S\/}\,(z)=%
{\rm incl\/}_{{\goth S\/}\!\setminus\!\Gamma_{\cal B\/}}^{\goth S\/}\,(z')
\ \hbox{\small et \/} \eta'\,\Im\,{\rm pr\/}_2\,(\gamma_{\goth p}^{-1}\,(z))>0
\}$$
en un {\it duplicata\/}, {\it c.a.d.\/} rev\^etement double,
$\pi_{\cal B\/} : \Omega_{\cal B\/}\rightarrow{\goth S\/}\!\setminus\!\Delta_{\cal B\/}\/$,   
trivial sur
${\goth S\/}\!\setminus\!\Gamma_{\cal B\/}\/$. 
%

\Thnc {\bf 1.1\ \/}|D\'efemme| L'ensemble {\it bord\/}
$\partial\Gamma\!=\!\{s\in\Gamma_{\cal B\/};\, v_{\cal B\/}\,(s) \equiv 1\,{\rm mod\/}\,2
\}\/$
des som\-mets pr\`es des\-quels $\pi_{\cal B\/}\/$ est non trivial,
ne d\'epend que de $\Gamma=\Gamma_{\cal B\/}\/$ et $\pi_{\cal B\/}\/$ est \'etendu par le
{\it duplicata\/} (ramifi\'e sur $\partial \Gamma\/$
si $\partial \Gamma\/$ est non vide) $\pi_\Gamma : {\goth S}_\Gamma\rightarrow {\goth S\/}\/$ {\it de\/} 
${\goth S\/}\/$ {\it sur\/} $\Gamma\/$.
\finnc
\Dem Soit $\gamma=\zeta^{-1}\/$  une  e\'ennodrooque $\Gamma\/$-normale $\gamma\/$ centr\'ee en
$s\in{\goth S\/}\/$.

Qui d\'ecrit 
$S^1\/$  
va $v_{\cal B\/}\,(s)\/$ fois, par $\zeta\,({\goth p\/}_\lambda)\/$ d'une composante d'un 
$\zeta\,({\goth p\/}_\lambda\!\setminus\!\beta_{\lambda}))\/$ \`a l'autre et $\zeta\circ\pi_{\cal B\/}\/$
induit un duplicata trivial de  $\,S^1\/$ seulement si
$\,v_{\cal B\/}\,(s)\/$ est paire.

En ce cas 
 $(\zeta\!\circ\!\pi_{\cal B\/})^{-1}\,(\overline{D}^\ast)\/$ se recolle avec
$\mu_2\!\times\!(P_1\,({\Bbb C\/})\!\setminus\!D)\/$ en un duplicata de
$P_1\,({\Bbb C\/})\!\setminus\!\{0\}\simeq{\Bbb C\/}\/$. Ce dernier trivial, $\pi_{\cal B\/}\/$
l'est pr\`es de $s\/$,
d'o\`u puisque%
, ensemble des $s\/$ ayant dans $\Gamma\/$ une base 
de voisinages $(U_n)_{n\in{\Bbb N\/}}\/$ avec $c\,(U_n)\equiv c\,(U_n^\ast)\,{\rm mod\/}\,2\/$\break (ou
 $U_n^\ast=U_n\!\setminus\!\{s\}\/$),
le bord de $\partial\Gamma\/$ d\'epend que de $\Gamma\/$,
la premi\`ere partie.

La seconde suit de la classification$^{19}\/$ des duplicatas de
$D^\ast\/$
\ffindem
Une section de ${\pi_\Gamma}_{|{\goth S\/}\!\setminus\Gamma}\/$ 
correspond \`a un $\eta\in{\goth E\/}\,({\goth S\/}\!\setminus\!\Gamma_{\cal B\/})\/$ 
et  $\pi_\Gamma\/$ est isomorphe \nobreak au\break duplicata $\pi_{\Gamma^\eta}\/$ sur
le {\it $\omega\/$-cycle $\Gamma^\eta\/$ du signe $\eta\/$\/},
union des $\beta_\lambda\/$ avec $\eta\/$ constant sur ${\goth c\/}_\lambda\!\setminus\!\beta_\lambda\/$.

Un tel  $\Gamma^{\eta_m}\!\/$ est minimal pour la co{\"{\i}}nclusion si et seulement si
il ne s\'epare aucune
composante connexe de ${\goth S\/}\/$. La co{\"{\i}}nclusion \'etant bon ordre sur les $\omega\/$-cycles de signe,
$\Gamma_{\cal B\/}\/$ contient
un $\omega\/$-cycle de signe minimal et ce d'une fa\c con relative :

\Thnc {\bf 1.2\ \/}|D\'ecofaire|
{ (i)\/} Soit $\Gamma_U\!\subset\!\overline{U}\cap\Gamma_{\cal B\/}\/$
le support d'un sous-graphe de ${\cal B\/}\/$
inclus dans la fermeture d'un  $\omega\/$-polygone $U\/$ de ${\goth S\/}\/$
avec  $\pi_{\Gamma_{\cal B\/}}\/$ trivial sur $U\!\setminus\!\Gamma_U\/$.
Alors ${\pi_{\Gamma_{\cal B\/}}}\/$ est isomorphe \`a $\pi_{\Gamma^{\eta_m}_{U}}\/$ o\`u 
le $(U,\,\Gamma_U)\/$-{\it support r\'eduit\/}
$\Gamma^{\eta_m}_U\/$ de $\Gamma_{\cal B\/}\/$  est minimal et v\'erifie
$${\Gamma_{\cal B\/}\supset\Gamma^{\eta_m}_U\supset\Gamma^{\eta_m}_U\cap{U}\subset\Gamma_U\cap U\/}$$
la derni\`ere inclusion \'etant une \'egalit\'e si $U\!\setminus\!\Gamma_U\/$ est connexe.

{ (ii)\/} Un (et tout) tel 
$\Gamma^{\eta_m}_U=\emptyset\/$
est vide et seule\-ment si 
$\pi_\Gamma\/$ est trivial.\ffindem
\finnc

Une {\it courbe simple (ferm\'ee)\/} est le support d'un $\omega\/$-graphe connexe (compact)
de
sommets de valence $2\/$, un cycle ne contenenant pas de sous-cycle strict.

Tout cycle
d'une surface de Riemann ${\goth S\/}\/$ contient une courbe simple, d'o\`u le
\Thnc {\bf 1.3\ \/}|Corollaire|{\it (i)\/} Une surface de Riemann connexe mais sans duplicata connexe
est s\'epar\'ee par toute courbe simple
 en  deux composantes connexes. 

{\it (ii)\/} Une surface de Riemann ${\goth S\/}\/$ est planaire connexe si et seulement toute courbe simple ferm\'ee
$\gamma\/$ de ${\goth S\/}\/$ a son compl\'ementaire de connexit\'e $c\,({\goth S\/}\!\setminus\!\gamma)=2\/$.\ffindem
\finnc

Une surface de Riemann, v\'erifiant la conclusion [resp. l'hypoth\`ese] de  {\bf 1.3\ \/} {\it (i)\/},
est {\it Jordan-simplement connexe\/} [resp. {\it acyclique\/}], l'\'enonc\'e de Jordan$^{3'}\/$ cit\'e
$\hbox{\rm p.\/}\,\hbox{\rm i\/}\/$ est
\Thnc {\bf 1.4\ \/}|Lemme|$\!\!\!\/$
Un $\omega\/$-polygone $U\/$ simple non co-compact de ${\Bbb C\/}$ est acyclique.
\finnc
\vskip-1mm
\Dem Soit  ${\goth p\/}_\lambda,\ {\lambda\in\Lambda}$ et
$\gamma_s,\,s\in\Delta=\Delta_U\/$
des familles de perles d'un $\omega\/$-graphe ${\cal B\/}\/$ 
de support $\Gamma_{\cal B\/}\!=\!{\rm Fr\/}\,U\/$ et d'e\'ennodrooques $\Gamma_{\cal B\/}\/$-normales
centr\'ees aux sommets de $U\/$
d'images disjointes et 
$p:\Omega\!\rightarrow\!U\/$ un duplicata de $U\/$.
 
La composante ${\goth p\/}^\ast_{\lambda}\cap U\/$ 
 de  ${\goth p\/}^\ast_{\lambda}={\goth p\/}_{\lambda}\!\setminus\!\beta_\lambda\/$ est ({\it Weyl-\/}) simplement
connexe [{\it c.a.d.\/} tous ses rev\^etements
sont triviaux], donc un duplicata
$\overline{p} : E\rightarrow\overline{U}\!\setminus\!\Delta_U\/$ \'etend $p\/$.

Ainsi $\overline{p}\/$ induit, pour $0<\epsilon<1$, un duplicata $\overline{p}^\epsilon\/$ de
la fermeture $\overline{U}_\epsilon\/$ de l'ouvert
${U}_\epsilon={U}\!\setminus\!\cup_{s\in\Delta_U\/}\gamma_s\,(D_\epsilon)\/$,
un  $\omega\/$-polygone simple de ${\Bbb C\/}\/$.

\vskip-2mm
\Thc Affirmation| Ce duplicata $\overline{p}^\epsilon\/$
 est restriction d'un duplicata $\pi^\epsilon\/$ de ${\Bbb C\/}$.
\finc
\vskip-2mm

 Si, sur ${\rm Fr\/}\,U_\epsilon\/$,$\overline{p}^\epsilon\/$ 
 est trivial, $\pi^\epsilon\/$ est recoll\'e de
${\rm pr\/}_2:\mu_2\!\times\!{\Bbb C\/}\!\setminus\!U_\epsilon\!\rightarrow\!{\Bbb C\/}\!\setminus\!U_\epsilon\/$
\`a $\overline{p}^\epsilon\/$.\hfill\break
Sinon ${\rm Fr\/}\,U_\epsilon\/$ est une courbe simple ferm\'ee que 
connecte \`a $\infty\/$, car ${\Bbb C\/}\!\setminus\!U_\epsilon\/$ non compact,
 une demi-droite $d\!\subset\!{\Bbb C\/}\!\setminus\!\overline{U}_\epsilon\/$.
Pour une topologie de sa source, $\overline{p}^\epsilon\amalg\pi_{d}=\pi^\epsilon\/$.\fcarre

 Duplicata d'un simplement connexe, $\pi^\epsilon\/$ a deux sections sur le connexe $U_\epsilon\/$.
Elles se prolongent \`a
$U=\cup U_\epsilon\/$, en deux sections trivialisant le duplicata $p\/$.\ffincdem

\Thnc {\bf 1.5\ \/}|Proposition|L'implication {\it (i)\/} de {\bf 1.3\/} est en fait une \'equivalence.
\finnc
\vskip-1mm
\Dem Soit $(\gamma_j)_{j\in J}\/$ une famille d'e\'ennodrooques normales 
  d'images  couvrant
${\goth S\/}\/$ et, pour une partie 
 $\phi\!\subset\!J\/$ finie, le $\omega\/$-polygone relati\-vement com\-pact
$U_\phi\!=\!\cup_{j\in \phi}\gamma_j\,(D)\/$,   
 et $\Gamma_\phi\!=\!\cup_{j\in \phi}\overline{\gamma}_j\,(S^1)\subset\overline{U}_\phi\/$, 
support d'un  $\omega\/$-graphe fini
${\cal B\/}_\phi\/$.

L'ensemble des parties finies non vides $\phi\subset J\/$ avec $U_\phi\/$  connexe est not\'e $\Phi\/$.\break
Un duplicata $p\/$ de ${\goth S\/}\/$, \'etant trivial sur chaque $\gamma_j\,(D)\/$, a sa restriction 
\`a $U_\phi\/$, pour $\phi\in\Phi\/$ isomorphe au duplicata sur le support 
d'un sous $\omega\/$-graphe
ferm\'e ${\cal B\/}^p_{\!\phi}\/$ de ${\cal B\/}_\phi\/$.

Selon {\bf 1.2\/} et le lemme des mariages$^{20,\,20'}\/$  il y a un $\omega\/$-graphe  ${\cal B\/}\/$ de
${\goth S\/}\/$ de support  
 l'union de r\'eduits de $\Gamma_{{\cal B\/}^p_\phi}\/$ et dont le duplicata $\pi_{\Gamma_{\cal B\/}}\/$
est isomorphe \`a $p\/$.

Soient $\gamma_s\/$ des e\'ennodrooques $\Gamma_{\cal B\/}\/$-normales centr\'ees en les sommets $s\in\Delta\/$ de ${\cal B\/}\/$ avec
les $\overline{\gamma}_s\,(\overline{D})\/$ disjointes. D'apr\`es {\bf 1.2\ \/} le
duplicata $p\/$ est isomorphe \`a $\pi_\Gamma\/$, celui
 d'un
support r\'eduit $\Gamma\/$  de l'\/$\omega\/$-graphe
$\Gamma'=(\Gamma_{\cal B\/}\!\setminus\!\cup_{s\in\Delta}\gamma_s\,(D))\cup_{s\in\Delta}\gamma_s\,(S^1)\/$.

Les sommets de $\Gamma'\/$ \'etant de valence $3\/$ et le support $\Gamma\/$ un cycle d'apr\`es {\bf 1.1\/},\break ce dernier,
s'il est non vide,
a ses sommets de valence $2\/$  donc est
 une courbe.

En ce cas ${\goth S\/}\/$ non s\'epar\'ee par le r\'eduit $\Gamma\/$ donc par chaque composante, une
courbe simple, la surface de Riemann ${\cal S\/}\/$ n'est pas Jordan-simplement
connexe.

Sinon $\Gamma\!=\!\emptyset\/$ et, isomorphe \`a
$\pi_{\emptyset}={\rm pr\/}_2:\mu_2\!\times\!{\goth S\/}\rightarrow{\goth S\/}\/\/$,
le duplicata	 $p\/$ est trivial. Le d\'ebut de preuve ayant pris ce dernier arbitraire, la r\'eciproque
de {\it (i)\/} est \'etablie.
\vskip-2mm\ffindem
  
\vfill\eject
%
\fancyhead[LO]{A. M. \quad Uniformisation des surfaces de Riemann}
  \fancyhead[RO]{\thepage}
  \fancyhead[RE]{\small\quad\quad Appendice 2\hfill Surfaces de Riemann \`a bord et double de Klein usuel\hfill\quad}
  \fancyhead[LE]{\thepage}
{\parindent=0pt\par\vskip .1cm
\vskip 0mm plus -20mm minus 1,5mm\penalty-50
{\bf 2\ \/}%
{\bf Surfaces de Riemann \`a bord et double de Klein usuel\/}
\nobreak\parindent=20pt}%

Une {\it surface de Riemann \`a bord\/} ${\bf R\/}\/$ est une vari\'et\'e  s\'epar\'ee
model\'ee sur les
ouverts de $D_+\/$ et
de changements de cartes des restrictions
d'ap\-pli\-ca\-tions holomorphes injectives d\'efinies sur des ouverts de $D\/$.
 
Si $p\in C\!\subset\!{\bf R\/}\/$ est un  point d'une carte,
 de coordonn\'ee $\zeta: C\rightarrow V\!\subset\!D_+\/$,
d'une surface de Riemann \`a bord ${\bf R\/}\/$
est d'image
$\zeta\,(p)=t\in{\Bbb R\/}\/$ r\'eelle, par le th\'eor\`eme$^{6}\/$ de l'app\-lication
ouverte, 
 il en sera de m\^eme pour toute coordonn\'ee en $p\/$.

L'ensemble des tels points est le
{\it bord\/} $\partial\,{\bf R\/}\/$ de la surface de Riemann ${\bf R\/}\/$.

Par prolongement analytique$^{18}\/$,
une extension holomorphe $f\/$ d'un changement de carte sur  un disque,
centr\'e en l'image $t=\zeta\,(p)\/$ d'un point du bord par une coordonn\'ee
v\'erifie la relation de commutation $f={\rm conj\/}\circ f\circ{\rm conj\/}\/$. 

Ainsi le {\it double\/} de ${\bf R\/}\/$, quotient de ${\bf R\/}\times\mu_2\/$
par la relation d'\'equivalence identifiant $(p,\, -1)\/$ \`a $(p,\,1)\/$ si $p\in \partial\,{\bf R\/}\/$,
muni des coordonn\'ees d\'efinies sur les doubles des cartes
$C_i, i\in I\/$ d'un atlas holomorphe
${(\zeta_i : C_i\rightarrow V_i\subset D_+)}_{i\in I}\/$ de ${\bf R\/}\/$
par la coordonn\'ee correspondante $\zeta_i\circ{\rm pr\/}_1\/$ sur
$C_i\times\{1\}\/$ et
${\rm conj\/}\circ\zeta_i\circ{\rm pr\/}_1\/$ sur $C_i\times\{-1\}\/$
est une surface de Riemann ${\goth D\/}\,{\bf R\/}\/$, la \/{\it surface de Riemann double\/} de la surface
de Riemann \`a bord ${\bf R\/}\/$, qui, si 
$\pi: {\bf R\/}\times\mu_2\rightarrow{\goth D\/}\,{\bf R\/}\/$ d\'esigne l'application quotient :
 
Contient la copie $\pi\,({\bf R\/}\times\{1\})\/$ de ${\bf R\/}\/$, encore not\'ee
${\bf R\/}\subset{\goth D\/}\,{\bf R\/}\/$.

\vskip-1mm
Et est munie de l'h\'et\'erolution $\sigma_{\bf R\/}\/$
d\'efinie par $\sigma_{\bf R\/}\,(\pi(q, \epsilon))=\pi\,(q,\,-\epsilon)\/$.

Cette {\it conjugaison (de double)\/} $\sigma_{\bf R\/}\/$
a ${\bf R\/}\/$, comme domaine fondamental, et la courbe analytique
r\'eelle $\partial {\bf R\/}={\rm Fix\/}\,\sigma_{\bf R\/}\/$ de ${\goth D\/}\,{\bf R\/}\/$, comme
ensemble de points fixe.

Le {\it double de Klein\/}
de  ${\bf R\/}\/$ est le triplet $({\goth D\/}\,{\bf R\/},\,{\bf R\/},\,\sigma_{\bf R\/}\/)\/$.

Une carte de ${\bf R\/}\/$ est {\it \'el\`ementaire\/} si elle est 
de double \'el\'ementaire. En ce cas elle a$^{13'}\/$ elle a un coordonn\'ee, dite {\it \'el\'ementaire\/}, d'image $D_+\/$.
Les bords de paires de cartes \'el\'ementaires ont quatre configuations possibles:
\Thnc {\bf 2.1\/}\ |Lemme|Soit $X=\partial C_-\cap\partial C_+\/$ l'intersection des bords de deux
cartes d'une surface de Riemann \`a bord ${\bf R\/}\/$,
sources, pour $\eta\!\in\!\mu_2\/$, des coor\-don\-n\'ees \'el\'ementaires $\zeta_\eta\!:C_\eta\rightarrow D_+\/$
et $h={\zeta_+}{|X}\circ({\zeta_-}{|X})^{-1} : \zeta_-\,(X)\rightarrow \zeta_+\,(X)\/$. Alors $X\/$ est soit

{\it i)\/} vide.
\vskip-.5mm
{\it ii)\/} inclus dans le bord de l'une des deux cartes.
\vskip-.5mm
{\it iii)\/} connexe \`a compl\'ementaire, dans chaque
$\partial C_\eta\/$, connexe et non vide.

\vskip-1mm
Quitte \`a renum\'eroter, il y a $a_-,\, b_+ \in ]-1,\,1[\/$
tels que la fonction $h\/$ est un hom\'eomorphisme
croissant de $]a_-,\, 1[=\zeta_-\,(X)\/$ sur $]-1,\,b_+[=\zeta_+\,(X)\/$.
\vskip-.5mm
{\it iv)\/} d'image par chacune, $\zeta_\eta,\, \eta\in\mu_2\/$ des deux coordonn\'ees, 
$]-1,\,1[\setminus [b_\eta,\, a_\eta]\/$ le compl\'ementaire
 d'un intervalle compact
de $]-1,\, 1[\/$.

\vskip-1mm
En ce cas, qui
arrive
si et seulement si l'union $C_1\cup C_2\/$ des deux cartes a pour bord une
composante compacte de $\partial\,{\bf R\/}\/$, la fonction
$h\/$
induit deux hom\'eo\-morphismes croissants
de $]-1, b_-[\/$ et $]a_-,\, 1[\/$ sur $]a_+,\,1[\/$ et $]-1,\,b_+[\/$
respec\-tivement. 
\finnc
\Dem Supposons que le cas {\it i)\/} ne se produit pas. Soit $Y\/$
une composante connexe de $X=\partial C_-\cap\partial C_+\/$, d'image par $\zeta_\eta\/$,
 pour $\eta\in\mu_2\/$,
un intervalle $I_\eta=]a_\eta,\, b_\eta[\/$ de $I=]-1,\, 1[\/$. Comme 
la fonction $h\/$ est restriction d'une application
holomorphe injective envoyant l'intervalle r\'eel $I_-\/$
sur l'intervalle r\`eel $I_+\/$ et des points proches dans
$D_+\!\setminus \! I\/$
dans des points de $D_+\!\setminus\! I\/$, sa d\'eriv\'ee est$^{21}\/$, sur $I_-\/$ non nulle
 et r\'eelle positive, ainsi $h\/$ est une bijection croisante de $I_-\/$
sur $I_+\/$.

Si le cas {\it ii)\/} n'a pas lieu non plus, l'une des extr\'emit\'es
$c_-\!\in\!\{a_-,\, b_-\}\/$ de $I_-\/$ est un point de $I\/$. Puisque $h\/$
est croissante born\'ee $h\,(t)\/$ a,  quand $t\/$
tend vers $c_-\/$, une limite $c_+\in\overline{I}=I\cup\{-1,\,1\}\/$.
Ce point $c_+\not\in I\/$ n'est pas dans l'intervalle $I\/$
car, sinon, les points $x_\eta={\zeta_{\eta}}_{|X}^{-1}\,(c_\eta)\/$
pour $\eta\in\mu_2\/$
seraient deux points de ${\bf R\/}\/$ distincts (car situ\'es respectivement
dans les ensembles disjoints $C_-\setminus C_+\/$ et $C_+\setminus C_-\/$)
et sans voisinages disjoints, contredisant la s\'eparation de ${\bf R\/}\/$. Ainsi,
comme $h\/$ est continue strictement croissante et
$I_+=h\,(I_-)\/$ strictement inclus dans $I\/$, car le cas
{\it ii)\/}\break est suppos\'e exclu, l'autre
extr\'emit\'e $d_-\in\{a_-,\, b_-\}\!\setminus\!\{c_-\}\/$ de $I_-\/$
n'est pas dans
l'inter\-valle $I\/$ et $\zeta_1\,(Y)\/$ est une des deux composantes
connexes de $I\setminus\{c_-\}\/$.

Comme trois telles sections d'un intervalle
ont une intersection non vide et la fonction $h\/$ est strict\-ement
croissante sur chaque composante connexe,  soit $X=Y\/$ est connexe, c'est le cas
{\it iii)\/},
sinon $X\/$ a deux composantes connexes : le cas {\it iv)\/}.

\ffindem

\vfill\eject%
%
\fancyhead[LO]{A. M. \quad Uniformisation des surfaces de Riemann}
  \fancyhead[RO]{\thepage}
  \fancyhead[RE]{\small\quad\quad Appendice 3\hfill  Bords et double de Klein des $\omega\/$-polygones\hfill\quad}
  \fancyhead[LE]{\thepage}
{\parindent=0pt\par\vskip 0cm
\vskip 0mm plus -20mm minus 1,5mm\penalty-50
{\bf 3\ \/}%
{\bf  bord et double de Klein des $\omega\/$-polygones\/}
\nobreak\parindent=20pt}%

\Thnc {3.1\ \/}\ | Lemme| Soit ${\goth S\/}\/$ une surface de Riemann alors 
{\it (i)\/} Un ouvert $U\/$  de ${\goth S\/}\/$ est un $\omega\/$-polygone ouvert de ${\goth S\/}\/$
si et seulement si ${\goth S\/}\!\setminus\!U\/$ est un $\omega\/$-polygone ferm\'e de ${\goth S\/}\/$.

{\it (ii)\/} Si $U\/$ est ouvert int\'erieur de sa fermeture et de fronti\`ere
une $\omega\/$-cha{\^{\i}}ne $\Gamma\/$ chacune de ses composantes connexes $U_i\/$ et $U\/$   sont
des $\omega\/$-polygones ouverts.

{\it (iii)\/} L'ensemble des $\omega\/$-polygones ouverts de ${\goth S\/}\/$ est stable par intersection finie.

{\it (iv)\/} L'ensemble des $\omega\/$-polygones ferm\'es de ${\goth S\/}\/$ est stable par r\'eunion finie.

{\it (v)\/} Le compl\'ementaire d'un $\omega\/$-polygone ferm\'e dans un $\omega\/$-polygone ouvert est un
$\omega\/$-polygone ouvert.
\finnc
\Dem   {\it (i)\/} Soit $V\!=\!{\rm Int\/}\,({\goth S\/}\!\setminus\!U)\/$
alors la condition portant
sur l'int\'erieur $U\!=\!{\rm Int\/}\,\overline{U}\/$ (resp. $V\!=\!{\rm Int\/}\,\overline{V}\/$) 
pour que $U\/$ (resp. $V\/$) soit un $\omega\/$-polygone ouvert, 
est (ou n'est pas) simultan\'ement
remplie par $U\/$ et $V\/$ car 
\'equivalente \`a ${\rm Fr\/}\,U\!=\!{\rm Fr\/}\,V\!\!\/$, qui remplie, donne aussi la sym\'etrie
de l'autre condition : la fronti\`ere
est un $\omega\/$-cycle.

{\it (ii)\/} Soit ${\cal B\/}\/$ un $\omega\/$-graphe de support $\Gamma\/$ et
${\goth p\/}\/$ un collier de ${\cal B\/}\/$. Si $\beta_\lambda\/$ est une ar\^ete de ${\cal B\/}\/$, la perle
perc\'ee ${\goth p\/}_\lambda\!\setminus\!\beta_\lambda\/$ est disjointe de ${\rm Fr\/}\,U\/$ ainsi, comme
${\rm Int\/}\,\overline{U}=U\/$, une de ses composantes est disjointe de $U\/$ et l'autre incluse dans $U\/$ donc
dans une composante $U_{i_\lambda}\/$ de $U\/$.
Chaque $\Gamma_i={\rm Fr\/}\,U_i\/$ est donc support d'un sous-$\omega\/$-graphe
${\cal B\/}^i\/$ de $\partial U\/$.
Les signes
$\eta\in{\goth E\/}\,({\goth S\/}\!\setminus\!\Gamma)\/$ et $\eta_i\in{\goth E\/}\,({\goth S\/}\!\setminus\!\Gamma_i)\/$,
 d\'efinis par $\eta^{-1}\,(1)=U\/$ et $\eta_i^{-1}\,(1)=U_i\/$, scindant $\pi_\Gamma\/$ et $\pi_{\Gamma_i}\/$,
ces duplicatas sont non ramifi\'es et, selon {\bf 1.1\/},
${\rm Fr\/}\,U=\partial\Gamma=\emptyset=\partial\Gamma_i={\rm Fr\/}\,U_i\/$.
Ainsi $\Gamma\/$ et $\Gamma_j\/$ sont des $\omega\/$-cycles.
D'autre part comme
$$U_i\subset {\rm Int\/}\,{\overline{U}}_i\subset{\rm Int\/}\,{\overline{U}}\cap{\overline{U_i}}\subset U\cap{\overline{U_i}}
=(\amalg_{j\in J} U_j)\cap{\overline{U_i}}=U_i\/$$
la premi\`ere condition $U_i={\overline{U_i}}\/$ est aussi v\'erifi\'ee et $U_i\/$ est un $\omega\/$-polygone ouvert.

{\it (v)\/} suit de {\it (i)\/} et {\it (iii)\/}, ce dernier, d'apr\`es {\it (i)\/}, \'equivalent \`a {\it (iv)\/},
 il suffit  d'\'etablir que l'inter\-section
de  deux
$\omega\/$-polygones ouverts $U\/$ et $V\/$  est un $\omega\/$-polygone ouvert :
$$U\cap V\!\subset\!{\rm Int\/}\,({\overline{U\cap V}})\!\subset\!{\rm Int\/}\,(\overline{U}\cap\overline{V})\!\subset\!%
{\rm Int\/}\,\overline{U}\cap{\rm Int\/}\,\overline{V}\!\subset\!U\cap V \hbox{\rm et \/}
U\cap V\!=\!{\rm Int\/}\,(\overline{{U\cap V}})$$
Les $\omega\/$-arcs de
$\partial U\/$ et $\partial V\/$ se coupant deux \`a deux en 
un, deux $\omega\/$-arcs,  ou un ensemble fini, la fronti\`ere
 ${\rm Fr\/}\,(U\cap V)\subset{\rm Fr\/}\,(U)\cup {\rm Fr\/}\,(V)\/$ est une
$\omega\/$-cha{\^{\i}}ne $\Gamma\/$ et {\it (ii)\/} conclut.
\vskip-3mm\ffindem
\Thnc{\bf 3.2\/}\ |D\'ecofaire|
Soit $U\/$ un $\omega\/$-polygone
ouvert, ${\cal B\/}\/$ un $\omega\/$-graphe de support $\partial U\/$ et
$\gamma\/$ une e\'ennodrooque $\partial U\/$-normale en un point $p\in\partial U\/$,
alors pour tout $\epsilon\in]0,\,1[\/$ :

{\it (i)\/} La {\it $(\epsilon,\,\gamma)\/$-\'etoile\/} de $U\/$ en $p\/$
est le $\omega\/$-polygone $U\cap\gamma\,(D_\epsilon)=U_{\epsilon}(=U_{\epsilon,\,\gamma})\/$.

{\it (ii)\/} 
L'\/$(\epsilon,\,\gamma)\/$-\'etoile en $p\/$ a $\frac{v_{\cal B\/}\,(p)}{2}\/$ composantes connexes
$U_{\epsilon,\,\gamma,\,i},\, 1\!\leq\!i\!\leq\!\frac{v_{\cal B\/}\,(p)}{2}$,
$\omega\/$-polygones simples, dits
{\it ($(\epsilon,\,\gamma)\/$-)secteurs de $U\/$ en $p\/$\/},  avec 
$U_{\epsilon',\,\gamma,\,i}\!\subset\! U_{\epsilon,\,\gamma,\,i}\/$  si $\epsilon'\!<\!\epsilon\/$.
\finnc
\Dem $U\/$ et $\gamma\,(D_\epsilon)\/$ \'etant
des $\omega\/$-polygones ouverts, {\it (i)\/} suit de {\bf 3.1\/}{\it (iii)\/}.
Chaque $U\cap\gamma\,(S_{\eta})$, pour $0\!<\!\eta\!<\!\epsilon\/$ est, par normalit\'e de $\gamma\/$ union de 
$\frac{v_{\cal B\/}\,(p)}{2}\/$ arcs ouverts deux \`a deux disjoints d'extr\'emit\'es sur les $v_{\cal B\/}\,(p)$
arr\^etes de $\partial U\cap U_\epsilon\/$  d'o\`u {\it (ii)\/}.\ffindem
\vfill\eject
\null\vfill

\TrimTop{-7pct}
\TrimBottom{-3pct}
\centerline{\BoxedEPSF{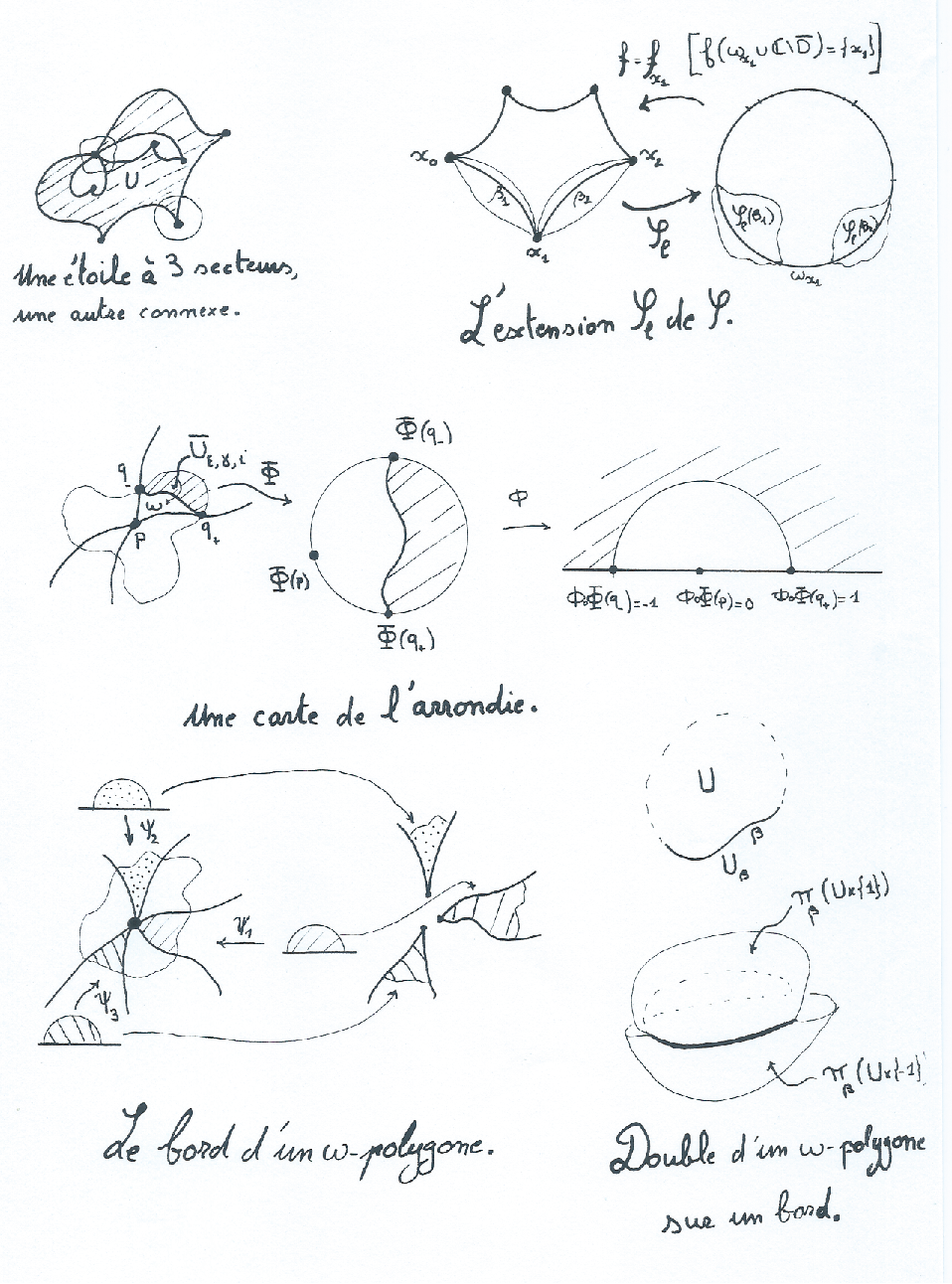 scaled 500} }
\centerline{Planche g}
\vskip2mm
\centerline{bord et double de Klein des $\omega-$polygones}
\vfill
\null\vfill\eject

\Thnc{\bf 3.3\/}\ |Proposition|Soit   $K\/$ une cellule d'int\'erieur $U\/$,
alors tout isomorphisme $\varphi : U\rightarrow D\/$ 
s'\'etend en un
hom\'eomorphisme $\Phi : K\rightarrow \overline{D}\/$.
\finnc
\Dem Soit $X\subset\partial K\/$ l'ensemble fini des sommets du
bord. D'apr\`es
le principe
de symm\'etrie de Schwarz$^{7"}\/$, $\varphi\/$ a
un prolongement $\varphi_l\/$ holomorphe injectif
d\'efini sur un voisinage de l'ensemble $K_l=K\!\setminus\!X\/$
des points lisses de $K\/$.

D'apr\`es {\bf 3.2\/}{\it (ii)\/} (ici les valences valent $2\/$) l'ensemble d'accumulation de
$\varphi\/$ sur les filtres des voisinages d'un sommet $x\!\in\!X\/$
est connexe, donc un arc de cercle $\omega_x,\, x\in X\/$. Comme
$\varphi_l\/$ est injective sur $\partial K_l\/$ ces arcs sont disjoints
et, la fermeture de l'image 
de $\varphi\/$ \'etant
$\overline{D}\supset S^1\/$, ils recouvrent le  compl\'ementaire
$S^1\!\setminus\!\varphi_l\,(\partial K_l)\/$.

Si l'int\'erieur $\alpha\/$ dans $S^1\/$
de l'un, $\omega_x\/$,
 de ces arcs est non vide alors
$U={\Bbb C\/}\!\setminus\!(S^1\!\setminus\!\alpha)\/$ est un ouvert connexe de
${\Bbb C\/}\/$ et l'application
$f : U\rightarrow S\/$ d\'efinie
par $f\,(z)=\varphi^{-1}\,(z)\/$, si $z\/$ est dans $D\/$
et $f\,(z)=x\/$, sinon est continue et holomorphe sur le compl\'ementaire
de l'arc $\alpha\/$. Le Corollaire de Morera$^{22'}\/$ assure qu'elle est holomorphe
sur $U\/$,
contredisant le
 prolongement analytique$^{18}\/$ puisque $f\/$ est
non constante sur l'ouvert connexe $U\/$, mais sa restriction
\`a l'ouvert ${\Bbb C\/}\setminus\overline{D}\/$ l'est.

Ainsi,
l'ensemble d'accumulation $\omega_x=\{\theta_x\}\/$ est r\'eduit \`a un
point et $x\mapsto \theta_x\/$ \'etend $\varphi_l\/$
en l'hom\'eomorphisme $\Phi$ cherch\'e.\ffindem
\Thnc{\bf 3.4\/}| Corollaire| Soit $U_{\epsilon,\,\gamma,\,i}\/$ un secteur en un point $p\!\in\!{\rm Fr\/}\,F\/$ 
fronti\`ere
d'un $\omega\/$-polygone ferm\'e $F\/$  d'une
surface de Riemann ${\cal S\/}\/$.

Il y  a un voisinage ouvert
$\omega\/$ de $p\/$ dans la fermeture $\overline{U}_{\epsilon,\,\gamma,\,i}\/$
hom\'eomorphe \`a
$D_+\/$
par 
$\Psi : (\omega,\,p)\rightarrow (D_+,\,0)\/$,
holomorphe sur ${U}_{\epsilon,\,\gamma,\,i}\/$ et telle que
$\Psi\,(\omega\cap {\rm Fr\/}\, F)=]-1, 1[\/$.
\finnc

\Dem Le secteur ferm\'e\ \ $\overline{U}_{\epsilon,\,\gamma,\,i}\/$ est, par {\bf 3.2\/}{\it (ii)\/}
et {\petcap Corollaire 0\/}, 
une cellule. 
Soit 
$\Phi : \overline{U}_{\epsilon,\,\gamma,\,i}\rightarrow\overline{D}\/$ l'hom\'eomorphisme, holomorphe sur
${U}_{\epsilon,\,\gamma,\,i}\/$
 de {\bf 3.3\/}.

Si $q_-,\,q_+\in \overline{U}_{\epsilon,\,\gamma,\,i}\cap\gamma\,(S_\epsilon)\/$ sont 
les sommets de $\overline{U}_{\epsilon,\,\gamma,\,i}\/$ distinct de $p\/$, l'hom\'eo\-mor\-phisme $\Psi\/$
est la restriction \`a
$V=\Phi^{-1}\,(D_+)\/$
de $\phi\circ\Phi\/$
 o\`u $\phi\/$ est l'homographie telle que $\Im \phi\,(0)>0\/$,
$\phi\,(\Phi\,(p))=0\/$ et $\{\phi\,(\Phi\,(q_+)),\, \phi\,(\Phi\,(q_-))\}=\{-1, 1\}\/$ .\ffindem

Soit $\Upsilon\/$ un ensemble indi\c cant tous les hom\'eomorphismes $\Psi\/$ que {\bf 3.4\/} donne.\break
Sur le quotient $B\/$
du produit $\Upsilon\times]-1,\,1[\/$ par l'\'equivalence  identifiant $(\upsilon,\, t)\/$
\`a $(\upsilon',\, t)\/$ si et seulement si $\Psi_\upsilon\,(t)=\Psi_{\upsilon'}\,(t')\/$ est d\'efinie 
 \/{\it le bord\/} du $\omega\/$-polygone $U\/$ :
$$\beta=\partial U : B\rightarrow {\rm Fr\/}\,U,\ {\overline{(\upsilon,\,t)}}\mapsto \Psi_\upsilon\,(t)\/$$

Si $\Upsilon_\beta\!=\!\{\upsilon\in\Upsilon; \omega_\upsilon\cap {\rm Fr\/}\, F\!\subset\!\beta\,(B)\}\/$
o\`u $\beta : B\rightarrow {\rm Fr\/}\,U\/$ est un bord d'un $\omega\/$-poly\-gone ouvert $U\/$,
le quotient de l'union disjointe ${\cal U\/}_\beta=U\amalg_{\upsilon\in\Upsilon_\beta}\omega_\upsilon\/$
par l'\'equi\-valence identifiant $v\in\omega_\upsilon\/$ \`a $u\in U\/$ si
${\rm incl\/}_{\omega_\upsilon}^{\goth S\/}\,(v)={\rm incl\/}_{U}^{\goth S\/}\,(u)\/$ est not\'e
$\rho_\beta: {\cal U\/}_\beta\rightarrow U_\beta\/$,\break une surface topologique s\'epar\'ee de bord 
$B=U_\beta\!\setminus\!\rho_\beta\,(U)\/$ et contenant $\Omega_\upsilon=\rho_\beta\,(\omega_\upsilon)\/$.
 
Sur le produit $U_\beta\times\mu_2\/$ l'\'equivalence $\sim\/$ identifiant
chaque point $(x, \epsilon)\in B\!\times\!\{\epsilon\}\/$ de la courbe $B\/$ de la 
copie d'indice $\epsilon\in\mu_2$ \`a son homologue $(x,-\epsilon)\in B\times\{-\epsilon\}\/$
dans l'autre copie et
$\pi_\beta : U_\beta\times\mu_2\rightarrow {\goth D\/}_\beta\,U=U_\beta\times\mu_2/_\sim\/$
l'application quotient.

Si $p=\rho_\beta\,(v)\in\Omega_\upsilon\/$ alors
$\Phi\,(\rho_\beta\,(v, \epsilon))={\rm conj\/}^{\frac{1-\epsilon}{2}}\circ\Psi_\upsilon\,(v)\/$
d\'efinit un hom\'eo\-morphisme
$\Phi_\upsilon\/$
du voisinage ${\goth D\/}_{\partial\,\Omega_\upsilon}\,\Omega_\upsilon\/$ de $\pi_\beta\,(p,\, \epsilon)\/$
dans ${\goth D\/}_\beta\,U\/$ sur le disque unit\'e.

Ces coordonn\'ees donnent
les structures de surface de Riemann transport\'ees par $\pi_\beta\/$, sur les 
$\pi\,({\rm Int\/}\,\Omega_\upsilon\times\{1\})\/$ et
$\pi\,({\rm Int\/}\,\Omega_\upsilon\times\{-1\})\/$ 
de celle de ${\goth S\/}\/$ et respectivement la
structure conjugu\'ee. D'o\`u un atlas holomorphe sur ${\goth D\/}_\beta\,U\/$
car, par le Corollaire de Morera$^{22'}\!\!\!\/$, elles sont aussi
holomorphiquement compatibles
entre elles.
   
L'application $\sigma_{\beta} : {\goth D\/}_\beta\,U\rightarrow {\goth D\/}_\beta\,U\/$ d\'efinie par
$\sigma_{\beta}\,(\pi_\beta\,(z,\, \epsilon))\!=\!\pi_\beta\,(z,\, -\epsilon)\/$ est, pour cette
structure holomorphe,
une h\'et\'erolution et a pour ensemble
de points fixe et domaine fondamental les sous-vari\'et\'es analytiques
images  de $B\/$ et $U_\beta\/$ par l'application continue injective et holomorphe
sur $\rho_\beta\,(U)\/$ :
$$i_{\beta} : U_\beta\rightarrow {\goth D\/}_\beta\,U,\ 
i_{\beta}\,(z)=\pi_\beta\,(z, 1)\/$$
Le
{\it double de Klein\/}
de l'$\omega\/$-polygone ouvert $U\/$ muni du bord $\beta\/$ est le triplet
$({\goth D\/}_\beta\,U,\,{\bf U\/}_{^a\beta},\,\sigma_{\beta}) \/$, double de Klein de
${\bf U\/}_{^a\beta}=i_{\beta}\,(U_\beta)\/$,
la {\it surface de Riemann \`a bord arrondie de $U\/$ sur $C\/$\/}. Le Corollaire
de Morera$^{22'}\/$ en donne la caract\'erisation :

\Thnc {\bf 3.5\/}|Proposition| Si ${\goth S\/}\/$
est une surface de Riemann munie de $\sigma\/$, une
h\'et\'ero\-lution  et si
$f: U_\beta\rightarrow{\goth S\/}\/$ est continue,
holomorphe sur $U\/$ et telle que
$f\,(B)\,\subset\,{\rm Fix\/}\,\sigma\/$,
alors il y a une unique application $F : {\goth D\/}_\beta\,U\rightarrow{\goth S\/}\/$
holomorphe, $(\sigma_{\beta},\, \sigma)\/$-invariante ({\it i.e.\/}
$\sigma\circ F=F\circ\sigma_{\beta}\/$) telle que $F\circ i_{\beta}=f\/$.\ffindem
\finnc

D'o\`u, puisque le compl\'ementaire dans la sph\`ere de Riemann
$P_1\,({\Bbb C\/})\setminus I\/$ d'un intervalle ferm\'e
non r\'eduit \`a un point $I\/$ 
 de $P_1\,({\Bbb R\/})\/$
 est
isomorphe au disque unit\'e : 

\Thnc {\bf 3.6\/}| Corollaire| Le double d'une carte simple r\'eguli\`ere sur un arc ouvert non dense de sa fronti\`ere
est \'equivariament isomorphe \`a $(D,\,D_+,\,{\rm conj\/})\/$,
le disque unit\'e muni du demi-disque Nord et de la conjuguaison complexe.\ffindem
\finnc
\vfill\eject
\fancyhead[LO]{A. M. \quad Uniformisation des surfaces de Riemann}
  \fancyhead[RO]{\thepage}
  \fancyhead[RE]{\small\quad\quad Appendice 4\hfill Le crit\`ere de Montel\hfill\quad}
  \fancyhead[LE]{\thepage}
{\parindent=0pt\par\vskip .3cm
\vskip 0mm plus -20mm minus 1,5mm\penalty-50
{\bf 4\ \/}%
{\bf Une caract\'erisation de carte dans les surfaces de Riemann\/}%
\nobreak\parindent=20pt}%

Une famille  $(C_i)_{i\in I}\/$ de parties de ${\goth S\/}\/$ est
{\it filtrante\/} si pour toute paire $\{C_i,\, C_j\}\/$ de ces parties, il y a
dans la famille une  $C_k\supset C_i\cup C_j\/$ 
les contenant toutes deux.
\Thnc {\bf 4.1\ \/}|Crit\`ere de Montel| Soit dans une surface de Riemann ${\goth S\/}\/$
un ouvert connexe $U\/$ union d'une famille filtrante $(C_i)_{i\in I}\/$
de cartes de ${\goth S\/}\/$ alors $U\/$ est
une carte de ${\goth S\/}\/$
qui,
si chaque $C_i\/$ est isomorphe \`a $D\/$,
 est isomorphe \`a $D\/$ ou \`a ${\Bbb C\/}\/$.
\finnc
\Dem Soit $p_0\in C_0\subset {\goth S\/}\/$ fix\'e, $\varphi : 5D\rightarrow U_1\/$ un isomorphisme sur un voisinage de $p_0\/$,
 avec $\varphi\,(0)=p_0\/$ et $U^0=U\!\setminus\!\varphi\,(\overline{D})\/$.
Comme $(C_i)_{i\in I}\/$ est filtrante, il y a $J\subset I\/$ tel que $(C_j)_{j\in J}\/$,
est filtrante, recouvre $U\/$, et chaque
$C_j
\/$ contient $\varphi(4D)=U_0\/$.

Si $j\in J\/$ soit
$\zeta'_j : C_j\rightarrow V_j\subset{\Bbb C\/}\/$ une coordonn\'ee avec
$\zeta_j\,(p_0)=0\/$.
Il y a $t>0\/$ tel que $\zeta_j=h_t\circ\zeta'_j\/$ v\'erifie
$\zeta_j\,(C_j\setminus\varphi\,(\overline{D}))\subset{\Bbb C\/}\setminus\overline{D}\/$
mais $\zeta_j\,(q_j)=\theta_j\in S^1\/$,
pour un point
$q_j=\varphi\,(\eta_j)\in \varphi\,(S^1)\/$.
D'apr\`es le th\'eor\`eme de
Montel$^{23}\/$, 
les restrictions
$${\zeta_j}_{|U^0} : U^0\rightarrow {\Bbb C\/}\!\setminus\!\overline{D}
\subset P_1\,({\Bbb C\/})\!\setminus\!\overline{D}\build{\rightarrow}_{^\sim}^{\imath}
D\/$$
des $\zeta_j\/$ \`a $U^0\/$
forment une famille normale de fonctions m\'eromorphes.

Il y a donc 
$K\subset J\/$, cofinale pour l'ordre sur $J\/$
d'inclusion des $C_j, j\in J\/$, telle que 
$(\zeta_k)_{k\in K}\/$
converge suivant le filtre des sections de $K\/$, uniform\'ement sur
tout compact de $U^0\/$, vers 
$\zeta : U^0\rightarrow\{\infty\}\cup{\Bbb C\/}\setminus\overline{D}\/$. 
Limite uniforme sur tout compact de fonctions holomorphes,
$\zeta\/$ est m\'eromorphe$^{8'''}\/$ et n'est pas une constante $c_v\/$ car :

i) si $v\ne\infty\/$ la formule de Cauchy$^{24}\/$, exprimant
 $f_j=\zeta_j\circ\varphi_{|D_2}:D_2\rightarrow{\Bbb C\/}$
par des int\'egrales sur $S_3\/$,
donnerait que les restrictions
${\zeta_j}_{|\varphi\,(D_2)}\/$ convergent uniform\'ement sur tout compact
vers $c_v\/$, contredisant $\zeta_j\,(p_0)=0\/$ et
$|\zeta_j\,(q_j)|=1\/$. 

ii) si $v=\infty\/$ le principe
de l'argument$^{25}\/$,  donne $j\in J\/$ tel que
$g_j\,(D_2)\!\supset\!g_0\,(\overline{D}_2)\/$, o\`u $g_j=\zeta_j\circ\varphi\/$.
Alors
$f=
h_{\frac{1}{2\,\eta_j}}\circ g_j^{-1}\circ h_{\frac{\theta_j}{\theta_0}}\circ g_0\circ h_{2\,\eta_0}\/$,
holomorphe non sujective de $D\/$ dans $D\/$,
violerait le lemme
de Schwarz$^{26}\/$, puisque 
$f\,(0)=0\/$ et $f\,(\frac{1}{2})=\frac{1}{2}\/$.

L'infini n'\'etant dans l'image des $\zeta_j\/$, le principe de l'argument
assure que $\infty\/$ n'est pas valeur de $\zeta\/$.
La formule de Cauchy donne la convergence uniforme sur tout compact de $\varphi\,(D_2)\/$.
Ainsi les $\zeta_j, j\in J\/$, holomorphes injectives,
 convergent
uniform\'ement sur tout compact de $U\/$ vers une fonction
holomorphe $\zeta\/$ qui n'\'etant pas constante est$^{27}\/$ un isomorphisme
holomorphe de $U\/$ sur un ouvert de ${\Bbb C\/}$.\fcarre

La derni\`ere assertion suit du {\petcap Th\'eor\`eme 0\/}, car une union $V\/$ filtrante d'ouverts $V_i\/$
de ${\Bbb C\/}$ v\'erifiant les hypoth\`eses du {\petcap Th\'eor\`eme 0\/}
en v\'erifie$^{28}\/$ la derni\`ere.\ffincdem
\vskip-2mm

\rmc Remarque|	
Th\'eor\`eme 0 s'\'evite  :
Choisir $V_i=D_{R_i}\/$, le Lemme de Schwarz$^{26}\/$ assure $R_i\leq R_j\/$ si $C_i\subset C_j\/$.
Ainsi $(D_{R_k})_{k\in K}\/$ est
filtrante 
d'union
$\cup D_{R_i}=D_{R}\/$, o\`u
$R={\rm sup\/}\,R_{i}\in]0,\,\infty]\/$ et  on convient $D_\infty={\Bbb C\/}$.
Partie de l'argument pr\'ec\'edent
donne $L\subset K\/$ cofinale telle que les
$\chi_{l}=\zeta\circ\zeta_{l}^{-1} : D_{R_{l}}\rightarrow {\Bbb C\/}\/$ convergent suivant le filtre
des sections de $L\/$ vers 
$\chi : D_{R}\rightarrow {\Bbb C\/}\/$. Par prolongement analytique$^{18}\/$, 
\'etant l'inclusion au voisinage de $0\/$, $\chi\/$ est $D_R\hookrightarrow{\Bbb C\/}\/$ donc $\ \cdots\/$\fcarre
\finc
\vfill\eject

\null\vskip-5mm

\pagenumbering{Roman}
\setcounter{page}{1}
\def\folio{$\emptyset\/$}
\fancyhead[LO]{A.M.\quad\quad  Uniformisation des surfaces de Riemann }
  \fancyhead[RO]{\thepage}
  \fancyhead[RE]{\quad\quad \small Epilogue \hfill Bibliographie comment\'ee\hfill}
  \fancyhead[LE]{\thepage}
%
\centerline{\Large  Epilogue}
\centerline{\bf Commentaires sur litt\'erature et tradition orale}
\vskip2mm
\parc

A l'automne 1989, J.-P. Demailly {\bf [D\^y]\/} racontait ainsi l'uniformisation :

Une surface de Riemann simplement connexe
non compacte ${\goth S\/}$
se couvre des composantes relativement compactes $U_\gamma\/$
du compl\'ementaire de courbes simples ferm\'ees $\gamma\/$ analytiques r\'eelles.
La surface de Riemann compacte simplement connexe \`a bord $U_\gamma\cup\gamma\/$,
a pour double une surface de Riemann
compacte simplement connexe isomorphe (par Riemann-Roch) \`a la droite projective, ainsi $U_\gamma\/$
est isomorphe au disque et gr\^ace \`a Montel ${\goth S\/}\/$ est isomorphe au disque ou au plan.

L'analyse \'etait donc \'elimin\'ee du th\'eor\`eme d'uniformisation puisque Riemann Roch ne n\'ecessite rien de plus
que le th\'eor\`eme de compacit\'e de Montel.  

Des id\'ees issues de la topologie g\'eom\'etrique ({\bf [Dy]\/} et {\bf [Ss]\/})
 donn\`erent le  Lemme A, l'extension dans la cat\'egorie des surfaces de Riemann
de la caract\'erisation topologique,
attibu\'ee \`a Brown, mais qu'on ne trouve que dans {\bf [Dy]\/}, de la sph\`ere comme seule vari\'et\'e compacte
couverte par deux disques permettant 
d'obtenir, sans Riemann-Roch ni g\'eom\'etrie analytique, le cas compact de l'uniformisation ({\bf [Mn]\/}).

La preuve  
utilisait un engoufrement par sym\'etrie de Schwarz dans une triangulation de la surface, proc\'ed\'e
 qui, bien auparavant (Cf. {\bf [Wn]\/}, {\bf [Cy]\/}) fourni \`a B.L. van der Waerden
le th\'eor\`eme d'uniformisation, en supposant la triangulabilit\'e, avec un minimum d'analyse, pr\'efigurant
l'\'enonc\'e moderne ({\bf [Me]\/}, {\bf [Be-Sn]\/}) :\hfill\break
{\sl Les surfaces de Riemann sans composante compacte sont des vari\'et\'es de Stein\/}\hfill\break
mais faisait dispara{\^{\i}}tre la belle id\'ee de double
gard\'ee dans le pr\'esent texte.

La th\'eorie d'Eilenberg dans le plan
donne que si $\varphi : C\rightarrow V\subset{\Bbb C\/}\subset P_1\,({\Bbb C\/})\/$
est une carte connexe dans une surface de Riemann planaire connexe ${\goth P\/}\/$ alors les composantes connexes
$E_i,\, i\in I\/$ de ${\goth P\/}\!\setminus\!C\/$ ont dans ${\goth P\/}\/$ des voisinages $U_i\/$
deux \`a deux disjoints
tels que pour chaque $i\!\in\!I\/$ il y a une composante $F_i\/$ de $P_1\,({\Bbb C\/})\!\setminus\!V$ avec 
${\cup}_{i\in I}^{} F_i\!=\!P_1\,({\Bbb C\/})\!\setminus\!V\/$ et
$\varphi\,(C\cap U_i)\cup F_i\/$ est un voisinage de $F_i\/$ dans $P_1\,({\Bbb C\/})\/$. Ceci permit d'obtenir par une
induction n\'ecessitant un ordre de s\'eparation sur les composantes du compl\'ementaire des composantes
d'une carte l'\'enonc\'e :\break
\vskip-2mm
{\sl Une surface de Riemann planaire compacte couverte de deux cartes est isomorphe \`a $P_1\,({\Bbb C\/})\/$.}

Le double d'un $\omega\/$-polygone simple relativement compact dans une surface de Riemann planaire
${\goth P\/}\/$ v\'erifiant les hypoth\`eses de cet \'enonc\'e, le  Th\'eor\`eme 1
en d\'ecoule, par une modification de la preuve donn\'ee p.10. Le th\'eor\`eme d'uniformisation vient alors
par la r\'eduction de Demailly.

Ceci (cont\'e en 1999) 
\'evite  Lemmes B et C mais exige
 plus de topologie :
{\sl Une des com\-posantes
du compl\'ementaire d'une courbe simple ferm\'ee dans une surface de Riemann simplement connexe est relativement compacte\/} et
{\sl le double d'un $\omega\/$-polygone simple dans une surface de Riemann planaire
est planaire\/}. Bien qu'il n'en soit rien dit dans {\bf [Mn]} ces deux derni\`eres
\og assertions \'evidentes\fg sont non triviales, surtout dans une surface non encore sue paracompacte.

Implicites dans toutes ces approches est une arithm\'etique de la s\'eparation dans les surfaces,
pr\'ecis\'ement \'enonc\'ee et \'etablie par Jordan ({\bf [Jn]\/})
et dont l'Appendice 1 donne une version pour les surfaces
de Riemann qui avec des modifications mineures s'applique
aux surfaces analytiques non n\'ecessairement paracompactes$^{29}\!\!\!\/$. Myst\'erieusement ce
\og calcul de Jordan\fg est bien moins cit\'e
que le \og Th\'eor\`eme de Jordan\fg prouv\'e dans {\bf [Jn]\/}
mais  
n'y ayant pas le statut d'\'enonc\'e.

L'alternative choisie dans le texte r\'eduit la topologie n\'ecessaire \`a la partie facile de cette 
\og topologie de Jordan\fg
(l'appendice 1 sans  {\bf 1.5\/}) et retrouve le Satz 1 de {\bf [Be-Sn]} :

{\sl Un ouvert rela\-tivement compact d'une surface de Riemann
est isomorphe \`a un ouvert \`a compl\'ementaire standard dans une surface de Riemann compacte.\/}

{\parindent=0pt o\`u il y avait une id\'ee proche de celle du double (et la base des constructions
 de l'appendice 1)  :
\nobreak\parindent=20pt}%

Pour  inciser le long d'ar\^etes et
\og recoller en croix des l\`evres de la plaie\fg
sans besoin d'applications quasi-conformes assurant l'inexistence de singularit\'e  aux sommets,
 Benke et Stein se limitent aux surfaces de Riemann \'etal\'ees sur $P_1\,({\Bbb C\/})\/$
et aux triangles pr\'eimages de
triangles  de c\^ot\'e arcs de cercle-droites : les recollements, sinon
l'identit\'e comme dans le cas du double, sont des homographies et l'uniformisation d'un voisinage \'epoint\'e
du sommet se donne  explicitement.

\finc
\vfill\eject
\def\folio{\uppercase\expandafter{\romannumeral \pageno}}
\centerline{\bf R\'ef\'erences}
\vskip 5mm
\hangindent=1cm\hangafter=1\noindent{\bf [Be-Sn]\/}\
{\petcap Benke H., Stein K}\pointir 
{\sl Entwicklung analytischer Funktionen auf Riemannschen Fl\" achen.},
Math. Annalen, {\bf 120\/}, 430-461, (1948).

\hangindent=1cm\hangafter=1\noindent{\bf [Cy]}\
{\petcap Caratheodory C}\pointir 
{\sl Conformal representation.},
Cambridge tracts in Mathematics n$^\circ$28,
Cambridge University Press, (1941).

\hangindent=1cm\hangafter=1\noindent{\bf [Cy']}\
{\petcap Caratheodory C}\pointir 
{\sl Theory of functions of a complex variable}, Vol 1  et 2, 
Chelsea, (1954).

\hangindent=1cm\hangafter=1\noindent{\bf [D\^y]}\
{\petcap Demailly J.-P}\pointir 
{\sl Comment r\'eduire le th\'eor\`eme d'uniformisation au cas compact}, 
Conversation de comptoir, Gerland, Automne 1989.

\hangindent=1cm\hangafter=1\noindent{\bf [Dy]}\
{\petcap Douady A}\pointir 
{\sl Plongements de sph\`eres (d'apr\`es Mazur et Brown)}, 
 Expos\'e N$^\circ$ 205 du S\'eminaire Bourbaki, 1-6, W.A. Benjamin (1960-1961).

\hangindent=1cm\hangafter=1\noindent{\bf [Dy']}\
{\petcap Douady A}\pointir 
{\sl Arrondissement des ar\^etes\/}, 
Expos\'e N$^\circ$ 3 du S\'eminaire H. Cartan, (1961-1962).

\hangindent=1cm\hangafter=1\noindent{\bf [Fr]}\
{\petcap Forster O}\pointir 
{\sl Riemannsche Fl\" achen}, 
Heidelberger Taschenb\" ucher 184, Springer-Verlag, (1977).

\hangindent=1cm\hangafter=1\noindent{\bf [Jn]}\
{\petcap Jordan C}\pointir 
{\sl Cours d'analyse}, 
Gauthier-Villars, (1894).

\hangindent=1cm\hangafter=1\noindent{\bf [Me]}\
{\petcap Malgrange B}\pointir 
{\sl  Existence et approximation des solutions des \'equations aux d\'eriv\'ees partielles et des \'equations
de convolution.}, Ann. Inst. Fourier {\bf 6\/}, 271-355, (1956).

\hangindent=1cm\hangafter=1\noindent{\bf [Mn]}\
{\petcap Marin A}\pointir 
{\sl  Le th\'eor\`eme d'uniformisation de Riemann d'apr\`es Demailly et
Stallings.}, Sem. di Geom. reale, 
Univ. Pisa, Dip. di Mat.,
Sez. di geom. e alg. 1.41,
{\bf 560\/}, 35-36, (1990).

\hangindent=1cm\hangafter=1\noindent{\bf [Na]}\
{\petcap Nevanlinna R}\pointir 
{\sl  Uniformisierung}, Grundlehren der Mathematischen Wissenschaften
LXIV, Springer-Verlag (1953).

\hangindent=1cm\hangafter=1\noindent{\bf [Ro]}\
{\petcap  Rad\'o T}\pointir
{\sl \"  Uber den Begriff der Riemannschen Fl\" ache.},
Acta Szeged {\bf 2\/}, 101-121, (1925).

\hangindent=1cm\hangafter=1\noindent{\bf [Rn]}\
{\petcap Rudin W}\pointir 
{\sl  Real and Complex Analysis}, McGraw-Hill (1966), 2$^{\rm nd\/}$ edition (1974).

\hangindent=1cm\hangafter=1\noindent{\bf [Sr]}\
{\petcap Springer G}\pointir 
{\sl  Introduction to Riemann surfaces}, Addison-Wesley (1957).

\hangindent=1cm\hangafter=1\noindent{\bf [Ss]}\
{\petcap  Stallings J.R}\pointir
{\sl The piecewise linear structure of Euclidean space.},
Proc. camb. Phil. Soc. {\bf 58\/}, 481-488, (1962).

\hangindent=1cm\hangafter=1\noindent{\bf [Wn]}\
{\petcap van der Waerden B.L}\pointir
{\sl Topologie und Uniformisierung der Riemann\-schen Fl\" achen.},
Ber. s{\H a}chs Akad. Wiss. Leipzig,  math.-phys. Kl. {\bf 93}, 147-160, (1941).

\hangindent=1cm\hangafter=1\noindent{\bf [Wl]}\
{\petcap  Weyl H}\pointir
{\sl Die Idee der Riemannschen Fl\" ache},
Teubner (1913)  et $2^{\rm nd\/}$ed. (1955).
\vfill\eject
%
%
\fancyhead[LO]{A. M. \quad Uniformisation des surfaces de Riemann}
  \fancyhead[RO]{\thepage}
  \fancyhead[RE]{\quad\quad\small Epilogue\hfill Notes \hfill}
  \fancyhead[LE]{\thepage}
\centerline{\bf Notes\/}
\vskip .2 cm
{\parindent=0pt\small
$^{0}\/$
et pas seulement voir {\bf [D\^y]\/}$^{0'}\/$, pour les textes,
on consultera la bibliographie de
{\bf [Na]\/} ainsi que les r\'ef\'erences donn\'ee par H. Weyl dans les notes du \S 20
 de l'\'edition de 1955 {\bf [Wl]\/} (p. 137-139).

$^{0'}\/$ \`a qui revient de porter le chapeau {\bf \^\/} de l'ajout
du pr\'esent opuscule \`a la litt\'erature.

$^{1}\/$\
p.137 de l'\'edition de 1955 de {\bf [Wl]\/}.


$^{2}\/$\
Cf. \S 315-323 de {\bf [Cy']\/}, {\bf 14.8\/} de {\bf [Rn]\/}
ou la preuve de Koebe dans l'exercice 26
du chapitre {\bf 14\/} de ce dernier manuel, plus concis,
choisi (\'edition de 1974) ici
pour r\'ef\'erer \`a l'analyse comp\-lexe \'el\'ementaire n\'ecessit\'ee :    
elle est incluse dans quatre chapitres de {\bf [Rn]\/} :
{\bf 10\/}, {\bf 11\/}, {\bf 12\/} et {\bf 14\/}.

$^{3}\/$\
{\sl Un domaine de Jordan n'a pas de rev\^etement double non trivial.\/}
\vskip-1mm
Pour la preuve sous cette forme et l'\'equivalence du sens \'etymologique donn\'e par Jordan \`a
l'ex\-pres\-sion \og simple connexit\'e\fg \`a celui qu'elle
 a depuis
{\bf [Wl]\/} voir l'appendice 1. L'original est :

$^{3'}\/$\
Lemme 511 (559) du t. II de la 2$^{\hbox{\rm i\`eme}}$  (3$^{\hbox{\rm i\`eme}}$ ) \`edition   1894 (tirage de 1959) de {\bf [Jn]}.

$^{4}\/$\
pour une preuve sans uniformisation
voir \S 23 de {\bf [Fr]\/} ou IV \S 3
de {\bf [Na]\/}. \'Etablie par Rad\'o ({\bf [Ro]\/}) en 1925 comme
corollaire de l'uniformisation des surfaces des Riemann s\'eparables, la paracompacit\'e (\'equivalente \`a la s\'eparabilit\'e de
chaque composante connexe)  faisait,  auparavant, partie de la d\'efinition des surfaces de Riemann
(voir {\bf [Wl]\/}, et aussi {\bf 2.18\/} et {\bf 6.16\/}
 de {\bf [Na]\/}).

$^{5}\/$\
sans \^etre exhaustif on peut citer {\bf [Wl]\/}, {\bf [Na]\/}, {\bf [Sr]\/}, 
{\bf [Fr]\/}. Le m\'emoire {\bf [Cy]\/} \'evoqu\'e dans les com\-men\-taires fait
 exception, mais utilise triangulations et topologie combinatoire.

$^{6}\/$\
par le th\'eor\`eme de l'application ouverte {\bf 10.32\/} de 
{\bf [Rn]\/}.

$^{7}\/$\
voir \S 341-342 de {\bf [Cy']\/}, \S 141-145 de {\bf [Cy]\/}. Les renvois prim\'es
adaptent {\bf 11.17\/} de {\bf [Rn]\/} :

$^{7'}\/$\
L'appli\-cation
$\Lambda : {\goth S\/}_{r^2}\!\setminus\!\varphi\,(S_r)\rightarrow{\Bbb C\/},\
\Lambda\,(z')=\zeta_r\,(z')\/$ si $z'\in {\goth S\/}_r\/$ et
$\Lambda\,(z)=
\sigma_1\circ\zeta_{r}\circ%
\varphi\circ\sigma_{r^2}\circ\varphi^{-1}\,(z)\/$ si
$z\in {\goth S\/}_{r^2}\!\setminus\!{\overline{\goth S\/}}_r\/$, est holomorphe  injective et le logarithme
de son
module se prolonge en
une fonction $\ell\/$ continue sur ${\goth S\/}_{r^2}\/$ telle que, si
$\psi\/$ et $\lambda=\log\,|\psi|\/$
sont d\'efinie sur
un voisinage $V\/$ de ${\Bbb R\/}\/$, convexe et invariant par conjugaison
comme compos\'ees de $\Lambda\/$ et $\ell\/$ avec
$\varphi\circ h_r\circ\exp\circ h_i\/$ alors $\lambda\/$ a la propri\'et\'e de
la moyenne, puisqu'harmonique hors de ${\Bbb R\/}$ et v\'erifiant
$\lambda\circ{\rm conj\/}=-\lambda\/$, et est donc harmonique
({\bf 11.16\/} de {\bf [Rn]\/}).
Ainsi 
$\lambda\/$  est partie imaginaire
d'une fonction holomorphe ({\bf 11.13\/} de {\bf [Rn]\/})
$h\/$, que l'on peut choisir de sorte que $g=\exp\circ h\/$ soit
\'egale \`a $\psi\/$ en un point de $V\!\cap\!H\/$, donc sur ce 
convexe et aussi sur $V\!\setminus\!{\Bbb R\/}\/$,
puisque, pour $\phi=g,\,\psi\/$ on a $\phi\circ{\rm conj\/}=\sigma_1\circ\phi\/$.
Ainsi $\exp\circ\ell\/$ est module d'une fonction holomorphe
$\xi : {\goth S\/}_{r^2}\rightarrow {\Bbb C\/}\/$ 
telle que
$\xi\,\circ\varphi\circ\sigma_{r^2}\circ\varphi^{-1}=\sigma_1\circ\xi\/$ et prolongeant $\Lambda\/$. Cette fonction $\xi\/$
est
injective, puisque$^{6}\/$, sur ${\goth S\/}_{r^2}\!\setminus\!\varphi^{-1}(S_r)\/$, elle l'est.

$^{7"}\/$\
Soit $\alpha : D\rightarrow {\goth S\/}\/$ une application holomorphe injective telle que
$\alpha\,(D_+)\/$ est un voisinage d'un point de $\partial K\!\setminus\!X\/$.
L'argument de la note pr\'ec\'edente $^{7'}\/$ se d\'eroule sur l'application
$\psi : D\!\setminus\!]-1,\,1[\rightarrow{\Bbb C\/}\/$ d\'efinie
par $\psi\,(z)=\varphi\circ\alpha\,(z)\/$, si $\Im z>0\/$ et
$\psi\,(z)=\sigma_1\circ\psi\circ{\rm conj\/}\/$, sinon.

$^{8}\/$\
{\bf 10.20\/} de {\bf [Rn]\/}. Dans les renvois prim\'es, l'appliquer aux fonctions holomorphes born\'ees :

$^{8'}\/$
la compos\'ee
$\varphi^{-1}\circ\xi^{-1}\circ\imath : {\Bbb C\/}^{\ast}\cap\imath\,(\xi\,(U))\rightarrow D\/$.

$^{8"}\/$
les compos\'ees $\imath\circ\psi\circ\varphi : D^\ast\rightarrow{\Bbb C\/}\/$
et $\varphi^{-1}\circ\psi^{-1}\circ\imath : {\Bbb C\/}^{\ast}\cap\imath\,(\psi\,(U))\rightarrow D\/$.

$^{8'''}\/$\
la compos\'ee
$\imath\circ\zeta :U^0\rightarrow P_1\,({\Bbb C\/}\!\setminus\!\overline{D})
\rightarrow D\/$.

$^{9}\/$\
une fonction holomorphe $f\/$
ne s'annulant pas sur $\xi\,({\goth S\/}_{r^2})\/$ a une racine carr\'ee holomorphe
$g_0\/$
sur un disque, voisinage de $\overline{D}\/$ dans $\xi\,({\goth S\/}_{r^2})\/$,
donc sur
$U_0=\overline{D}\cup\sigma_1\circ\zeta_r\circ
\varphi\,(D_{r}\!\setminus\!D_{r'})\/$
pour un $r'<r\/$, elle a aussi une racine carr\'ee holomorphe $g_\epsilon\/$
sur 
$U_\epsilon=\sigma_1\circ\zeta_r\!\circ%
\varphi\,(D\!\setminus\!(D_{r^2}\cup{\Bbb R\/}_\epsilon))\/$ pour
$\epsilon\!\in\!\mu_2\/$, puisque ces ouverts sont \'el\'ementaires.
 Les intersections $U_0\cap U_\epsilon\/$ pour $\epsilon\in\mu_2\/$
 \'etant connexes, il y a $\eta_\epsilon\in\mu_2\/$ tel que
$\eta_\epsilon g_\epsilon=g_0\/$ sur
$U_0\cap U_\epsilon\/$. Comme les deux composantes connexes
de $U_-\cap U_+\/$ rencontrent $U_0\cap U_1\/$,
ces $\eta_\epsilon g_\epsilon\/$ se recollent avec $g_0\/$ en une racine
carr\'ee holomorphe de $f\/$.
\vfill\eject
$^{10}\/$\
{\sl Une h\'et\'erolution $\sigma\/$ de $P_1\,({\Bbb C\/})\/$  est conjugu\'ee
\`a\/ $\sigma_{-1}\/$ si elle est libre et \`a\/ $\sigma_1\/$ sinon.\/}
\vskip-1mm

{\parindent=20pt
En effet, $c_\infty\/$\- et ${\rm Id\/}\/$
 n'\'etant des h\'et\'erolutions,  il y a 
$z_1\in{\Bbb C\/}\/$ avec $z_1\ne\sigma\,(z_1)=z_2\in{\Bbb C\/}$.
Pour $\lambda\in{\Bbb C\/}^\ast\/$, soit $\varphi : z\mapsto\frac{z-z_1}{z-z_2},\ %
\varphi_\lambda=h_\lambda\circ\varphi\/$
et
$f_{\epsilon,\, \lambda}=
\sigma_\epsilon\circ\varphi_\lambda\circ\sigma\circ\varphi_\lambda^{-1}\/$,
isomorphisme de $P_1\,({\Bbb C\/})\/$ fixant $0\/$ et $\infty\/$,
a en $0\/$ pour d\'eriv\'ee
$f'_{\epsilon,\,\lambda}\,(0)=\epsilon\,|\lambda|^{-2}f'_{1,\,1}\,(0)\/$, et il y a
$f=f_{\epsilon,\,\lambda}\/$ avec $f'\,(0)=1\/$.
\vskip-1mm

Alors $z\mapsto\frac{\sigma_1\circ f\circ\sigma_1^{-1}\,(z)}{z}\/$,
est, par le th\'eor\`eme$^{8}\/$ des singularit\'es inexistantes, restriction
 d'une fonc\-tion enti\`ere tendant vers $1\/$ quand $z\/$ tend vers l'infini donc $c_1\/$,
selon le th\'eor\`eme$^{11}\/$ de Liouville. Ainsi $f={\rm Id\/}\/$ et
$\sigma=
\varphi_\lambda^{-1}\circ\sigma_\epsilon\circ\varphi_\lambda\/$
est conjugu\'ee \`a $\sigma_\epsilon\/$.\ffindem
\nobreak\parindent=0pt}

$^{11}\/$\
{\bf 10.23\/} de {\bf [Rn]\/}.

$^{11'}\/$\
de fermeture incluse, si $1<t'<t\/$, dans
${\goth G\/}^{t'}\!\!\!\/$, ce domaine ${\goth G\/}^t\/$
est isomorphe \`a
un ouvert du compl\'ementaire dans ${\Bbb C\/}\/$ d'un disque $\{|z-z_0|\leq\epsilon\}$,
sur lequel $\imath\circ\tau_{-z_0}\/$ est holomorphe born\'ee non constante.
Il ne peut donc, par le th\'eor\`eme de Liouville$^{11}\/$
\^etre isomorphe \`a ${\Bbb C\/}\/$.

$^{12}\/$\
Deux points distincts $\pi\,(x_1,\, n_1)\ne\pi\,(x_2,\, n_2)\/$
ont pour voisinages disjoints,
$p^{-1}\,(\omega_i)\/$ pour $i=1,\,2\/$ o\`u $\omega_i\/$  sont
des voisinages disjoints des $x_i\/$  si $x_1\ne x_2\/$
et, si $x_1=x_2\/$, les images par $\pi\/$ des ouverts satur\'es
$((U_{{\overline{n_i-1}}}\!\cap\!U_{{\overline n}_i})\times\{n_i-1\})\cup %
(U_{{\overline n}_i}\times\{n_i\})
\cup((U_{{\overline n}_i}\!\cap\!U_{\overline{n_i+1}})\times\{n_i+1\})\/$ 
qui sont disjoints car, comme $|n_1-n_2|>1\/$ la condition {\it (ii)\/} et
le choix de $N\/$ donnent $|n_1-n_2|\geq N>2\/$.

$^{13}\/$\
les isomorphismes
de $D\/$ sont ({\bf 12.6\/} de {\bf [R]\/}) les $h_\lambda\circ\varphi_{\alpha}\/$, pour $|\lambda|=1\/$
(et $|\alpha|<1\/$).

$^{13'}\/$ {\sl Une h\'et\'erolution $\varsigma\/$  de $D\/$ est conjugu\'ee \`a ${\rm conj\/}\/$.\/}
\vskip-1mm

{\parindent=20pt
En effet$^{13}\/$, $\varsigma\,({\bar z})=h_\lambda\circ\varphi_{\alpha}\,(z)\/$,
soit $t=\frac{1-\sqrt{1-|\alpha|^2}}{|\alpha|^2}\/$ et $\alpha'=t{\bar \alpha}\/$.
Ainsi, pour $t\in[-1,\,1]\/$ :
\vskip-1mm
$$\varsigma\,(\alpha')=\varsigma\,(t\,{\bar \alpha})=%
\frac{t-1}{1-t\,|\alpha|^2}\,{{\lambda\,\alpha}}=%
\frac{1-t}{1-t\,|\alpha|^2}\,{{\varsigma\,(\varsigma\,({\bar \alpha}))}}\,%
=\frac{1-t}{1-t\,|\alpha|^2}\,{\bar \alpha}=\alpha'\/$$
ainsi $\alpha'\/$ est point fixe de
$\varsigma\/$ et
$\varphi_{\alpha'}\circ\varsigma\circ\varphi_{\alpha'}^{-1}=\lambda'\,{\bar z}\/$,
conjugu\'ee \`a ${\rm conj\/}\/$ (par $h_\mu\,z\/$ si $\lambda'\,\mu^2=1\/$).\ffindem 
\nobreak\parindent=0pt}

$^{13''}\/$\
par $^{13}\/$, $\theta=h_\lambda\circ\varphi_{\alpha}\/$ o\`u $\alpha\in]-1,\,1[\/$ et $\lambda\in\mu_2\/$,
 par \'equivariance
$\theta\circ{\rm conj\/}={\rm conj\/}\circ\theta\/$.
 Donc $\alpha\in]-1,\,1[\!\setminus\!\{0\}\/$ car $\theta\/$ libre sur $D\/$.  
Alors $-\frac{\alpha}{|\alpha|}\,\vartheta_+\/$ conjugue
$\theta\/$ \`a ${h_t}_{|L_+} : L_+\rightarrow L_+\/$ o\`u $t=\frac{1+|\alpha|}{1-|\alpha|}>1\/$.

$^{14}\/$\
Si $p_i=\pi\,(i,\,z_i)\/$ l'inverse $G=H^{-1}\/$ peut \^etre donn\'e par : si
$x=\pi(i,\,\imath\,(u))\/$, o\`u $|u|>\sqrt{r_i}\/$ alors
$G\,(x)=
\pi\,(i,\, \imath\,((\frac{\sqrt{r_i}}{|u|})^{N_i-1}\,\imath\,(u)+%
(1-(\frac{\sqrt{r_i}}{|u|})^{N_i})\,\imath\,(z_i))\/$
pour un $N_i>0\/$  grand et $G\,(x)=x\/$ sinon.

$^{15}\/$\
Si $\Gamma\/$, compact dans $U_i\/$, ouvert  connexe de ${\goth Q\/}\/$, connexe et localement
connexe et 
$C=U_i\!\setminus\!\Gamma\/$ est connexe,
alors
${\rm Fr\/}\,C\supset{\rm Fr\/}\,U_i={\rm Fr\/}\,({\goth Q\/}\!\setminus\!U_i)\ (\ne\emptyset\/$  si $U_i\ne{\goth Q\/}\/$)
et
${\goth Q\/}\setminus\Gamma=C\cup({\goth Q\/}\!\setminus\!U_i)\/$ est connexe.

$^{16}\/$\
Soit $\Gamma\!=\!\build{\cup}_{\lambda'\in \Lambda'}^{}\!\!\beta'_{\lambda'}$ une cha{\^{\i}}ne et
$\beta'_{\lambda'}\!\cap\!(\partial\beta'_{\lambda'}\cup\!\!\!\build{\cup}_{\lambda"\ne \lambda'}^{}\beta'_{\lambda"})\!=\!%
\{\gamma_{\lambda'}\,(t_{k}^{\lambda'});\, 0\!=\!t_{0}^{\lambda'}\!<\!\cdots\!<t_{l_{\lambda'}}^{\lambda'}\!=\!1\}\/$.
Pour
$\lambda=(n,\,\lambda')\in\Lambda=\{(n,\,\lambda')\in{\Bbb N\/}\times \Lambda';\, 0\leq n< l_{\lambda'}\}\/$,
soit $\gamma_\lambda=\gamma'_{\lambda'}\circ\tau_{t_{n}^{\lambda'}}\circ h_{t_{n+1}^{\lambda'}-t_{n}^{\lambda'}}\/$
et $\beta_\lambda=\gamma_\lambda\,(I)\/$.
La famille ${\cal B\/}=(\beta_\lambda)_{\lambda\in \Lambda}\/$ est un graphe 
d'union $\Gamma\/$, ferm\`e si $\Gamma\/$ est un $\omega\/$-cycle.

$^{17}\/$\
Soit $\beta_\mu\/$  avec $\emptyset\ne\beta_\mu\cap\beta_\lambda\/$, il y a pour $\nu\in\{\lambda,\,\mu\}$,
quitte \`a changer certains $\gamma_\nu\/$ en $\gamma_\nu\circ\tau_1\circ h_{-1}\/$, une param\'etrisation
 $\gamma_\nu\/$ de $\beta_\nu\/$ avec $\{\gamma_\nu\,(0)\}\,=\,\{p\}\,=\,\beta_\mu\cap\beta_\lambda\/$.
Soit $\,\sum_{m=1}^{\infty}\, a_m\,z^m\/$
(o\`u $\,a_1\ne0\/$)\break le d\'eveloppement en $0\/$ de l'isomorphisme local $\,\zeta_\lambda\circ\gamma_\mu\/$.
Puisque $\beta_\lambda\ne\beta_\mu\/$,  il y a, par
analyticit\'e$^{18}\/$ un entier
$m_{\lambda,\,\mu}\/$ avec $\,m_{\lambda,\,\mu}=0\,\/$ si $\,a_1<0\/$,
et sinon $\,\Im\,a_{m_{\lambda,\,\mu}}\ne0\/$
mais $\,\Im\,a_m=0\/$ pour $\,m<m_{\lambda,\,\mu}\/$.
Soit $N_\lambda\!=\!2\,\max \{m_{\lambda,\,\mu};\,\beta_\mu\cap\beta_\lambda\!\ne\!\emptyset\}\/$ et
$h'_\lambda>0\/$ grand avec $\gamma_\lambda\/$ d\'efini sur
$P'_\lambda\!=\!P_{N_\lambda,\,h'_\lambda}\/$, il y a un
voisinage $U\/$ de $\Delta_{\cal B\/}\/$ disjoint, pour tout $\lambda\ne \mu\/$, de
$\gamma_\lambda\,(P'_\lambda)\cap\gamma_\mu\,(P'_\mu)\/$. La famille ${\cal B\/}\/$ \'etant localement finie et, 
les $\gamma_\lambda\,(P_{N_\lambda,\,h})\/$ d\'ecroissant en $h\/$ d'intersection les arcs ouverts deux \`a deux disjoints
$\beta_\lambda\!\setminus\!\partial\beta_\lambda\/$, il y a des $h_\lambda\geq h'_\lambda\/$
tels que les ${\goth p\/}_\lambda=\gamma_\lambda\,(P_{N_\lambda,\,h_\lambda})\/$
sont deux \`a deux disjoints.

$^{18}\/$\
{\bf 10.18\/} de {\bf [Rn]\/}.

\vfill\eject
$^{19}\/$\
{\sl Un duplicata non trivial $p:X\!\rightarrow\!D^\ast\/$
 est isomorphe au {\it carr\'e\/} 
$q : D^\ast\!\rightarrow\!D^\ast,\, z\mapsto z^2\/$.\/}
\vskip-1mm

{\parindent =20pt
De source $L_-\/$ simplement connexe,
$\tilde{p}=\exp_{|L_-}=p\circ\pi : L_-\rightarrow X\rightarrow D^\ast\/$ se factorise par $p\/$, ce dernier
 \'etant un duplicata,
$\pi\circ\tau_{2\,\pi\,i}=\tau\circ\pi\/$ pour un isomorphisme $\tau:X\rightarrow X\/$
de $p\/$.
\vskip-1mm

Ainsi
$\pi\circ\tau_{2\,2\,\pi\,i}=\pi\circ\tau_{2\,\pi\,i}\circ\tau_{2\,\pi\,i}=\tau^2\circ\pi=\pi\/$,
et $\pi\circ h_{2}=\phi\circ\exp\/$ pour $\phi:D^\ast\rightarrow X\/$ holomorphe v\'erifiant
$p\circ\phi\!=\!q\/$ et $p^{-1}\,(q\,(z))\!=\!\{\phi\,(z),\,\tau\,(\phi\,(z))\!=\!\phi\,(-z)\}\/$,
si $z\!\in\!D^\ast\!\/$ 
car  $\tau\/$ est libre, puisque  $D^\ast\!\/$ est connexe et $p\/$ non trivial.
Bijectif, l'holomorphisme $\phi\/$ est un isomorphisme.\nobreak\ffindem
\nobreak\parindent=0pt}

$^{20}\/$\
Une famille,
index\'ee par les couples comparables $\phi<\phi'\in\Phi\/$ d'un ensemble bien ordonn\'e filtrant,
de \/{\it d\'esirs\/} $\pi_{\phi}^{\phi'} : F_{\phi}^{\phi'}\rightarrow F_{\phi}\/$
surjectifs sur l'ensemble fini $F_{\phi}\/$ d\'efini sur l'ensemble
$F_{\phi}^{\phi'}\!\subset\!F_{\phi'}\/$ 
{\it (le) pr\'etendant (dans $F_{\phi'}\/$)\/} est {\it compatible\/} (ou d\'efinit un {\it d\'esir syst\'ematique\/})  
si :

\vskip-1mm
{\parindent =20pt
Pour tout $\phi<\phi'<\phi"\/$, alors
$\pi_{\phi}^{\phi"}=\pi_{\phi}^{\phi'}\circ\pi_{\phi'}^{\phi"}\/$.
Sa {\it limite projective\/} $\limproj_{\, \phi\in\Phi}(F_\phi,\pi_{\phi}^{\phi'})\/$,
ou ensemble de ses {\it r\'ealisations\/} ${\cal R\/}\/$
  est l'ensemble des $(f_\phi)_{\phi\in\Phi}\/$ tel que pour tout $\phi<\phi'\/$ la $\phi'^{\hbox{\small i\`eme\/}}\/$
composante $f_{\phi'}\in F_{\phi}^{\phi'}\/$ est un pr\'etendant de d\'esir $\pi_{\phi}^{\phi'}\,(f_{\phi'})=f_\phi\/$.
\vskip-1mm
\thc Lemme des mariages|Un d\'esir syst\'ematique a des r\'ealisations :
${\displaystyle{\cal R\/}=\limproj_{\, \phi\in\Phi}(F_\phi,\pi_{\phi}^{\phi'})\ne\emptyset}\/$.\finc
\vskip-2mm
Car ${\cal R\/}\/$ 
est intersection ferm\'ee, avec propri\'et\'e de l'intersection finie,
du compact  $\prod_{\phi\in \Phi\/}F_\phi\/$.
\vskip-2mm
\ffindem
\nobreak\parindent=0pt}

$^{20'}\/$\
Appliqu\'e \`a la famille des
$F_\phi=\{\Gamma^{\eta_m}; \hbox{\small supports r\'eduits de \/} \Gamma_{{\cal B\/}^p_\phi}\},\, \phi\in\Phi\/$
que {\bf 1.2\/} munit des d\'esirs compatibles
$\pi_{\phi}^{\phi'}, \Gamma^\eta\mapsto\Gamma^\eta\cap\overline{U}_\phi\/$
sur les pr\'etendants
$ F_{\phi}^{\phi'}=\{\Gamma^{\eta}\in F_{\phi'}; \Gamma^{\eta}\cap\overline{U}_\phi\in F_\phi\}\/$.

$^{21}\/$\
{\bf 10.33\/} de {\bf [Rn]\/}.

$^{22}\/$\
{\bf 10.17\/} de {\bf [Rn]\/}, on en tire le \og Petit principe$^{7}$ de sym\'etrie\fg, dit aussi \og Corollaire de Morera\fg :

$^{22'}\/$\
{\sl Une fonction continue, holomorphe hors d'une courbe analytique r\'eelle, est holomorphe\/}.

$^{23}\/$\
{\bf 14.6\/} de {\bf [Rn]\/}.

$^{24}\/$\
{\bf 10.15\/} de {\bf [Rn]\/}.

$^{25}\/$\
Tout
$w\in g_0\,(D_2)\/$ est, car
$g_j\rightarrow\infty$ uniform\'ement sur $S_2\/$, dans la composante de
$0\!=\!g_j\,(0)\/$
de ${\Bbb C\/}\!\setminus\!g_j\,(S_2)\/$. Par {\bf 10.43\/} de {\bf [Rn]\/}, 
$f_w\!=\!\tau_{-w}\circ g_j\/$ a, dans $D_2\/$ autant de z\'eros que $f_0\/$,
donc $w\in g_j\,(D_2)\/$.

$^{26}\/$\
{\bf 12.2\/} de {\bf [Rn]\/}.

$^{27}\/$\
voir le dernier paragraphe de la preuve de {\bf 14.8\/} dans {\bf [Rn]\/}.

$^{28}\/$\
Si $g_{i_j}^2=f_{|V_{i_j}},\, i=1,\, 2\/$ pour $f:V\rightarrow{\Bbb C\/}^\ast\/$ et $g_{i_j}\/$
holomorphes avec
$g_{i_1}\,(z_0)=g_{i_2}\,(z_0)\/$, alors, sur $V_{i_1}\cap V_{i_2}\/$, $g_{i_1}=g_{i_2}\/$ 
 car, sur le connexe $V_{i_j}\/$,
$g_{i_j}\/$ est restriction de la racine carr\'ee $g_{k}\/$ de $f_{|V_{k}}\/$ o\`u $V_{k}\supset V_{i_1}\cup V_{i_2}\/$,
v\'erifiant
$g_{k}\,(z_0)=g_{i_1}\,(z_0)=g_{i_2}\,(z_0)\/$ : ces $g_i\/$ se recollent
en une racine de $f\/$.

$^{29}\/$\
Car si $\alpha : ]-1,\,1[\subset U\rightarrow{\Bbb R\/}^2\/$ est une courbe analytique r\'eelle d\`efinie sur un voisinage
simplement connexe de $]-1,\,1[\/$, la primitive 
co\"{\i}ncidant avec $\alpha\/$ en $0\/$, de l'application holomorphe
$x+i\,y=z\mapsto \frac{\partial}{\partial x}\,\alpha\,(z)-\frac{\partial}{\partial y}\,\alpha\,(z)\/$,
\'etend la restriction de $\alpha\/$ \`a $]-1,\,1[\/$
et, si $\alpha_{|]-1,\,1[}\/$ est injective, est injective sur un voisnage de $]-1,\,1[\/$.
}

\vfill\eject
\null\vfill

\TrimTop{-7pct}
\TrimBottom{-3pct}
\centerline{\BoxedEPSF{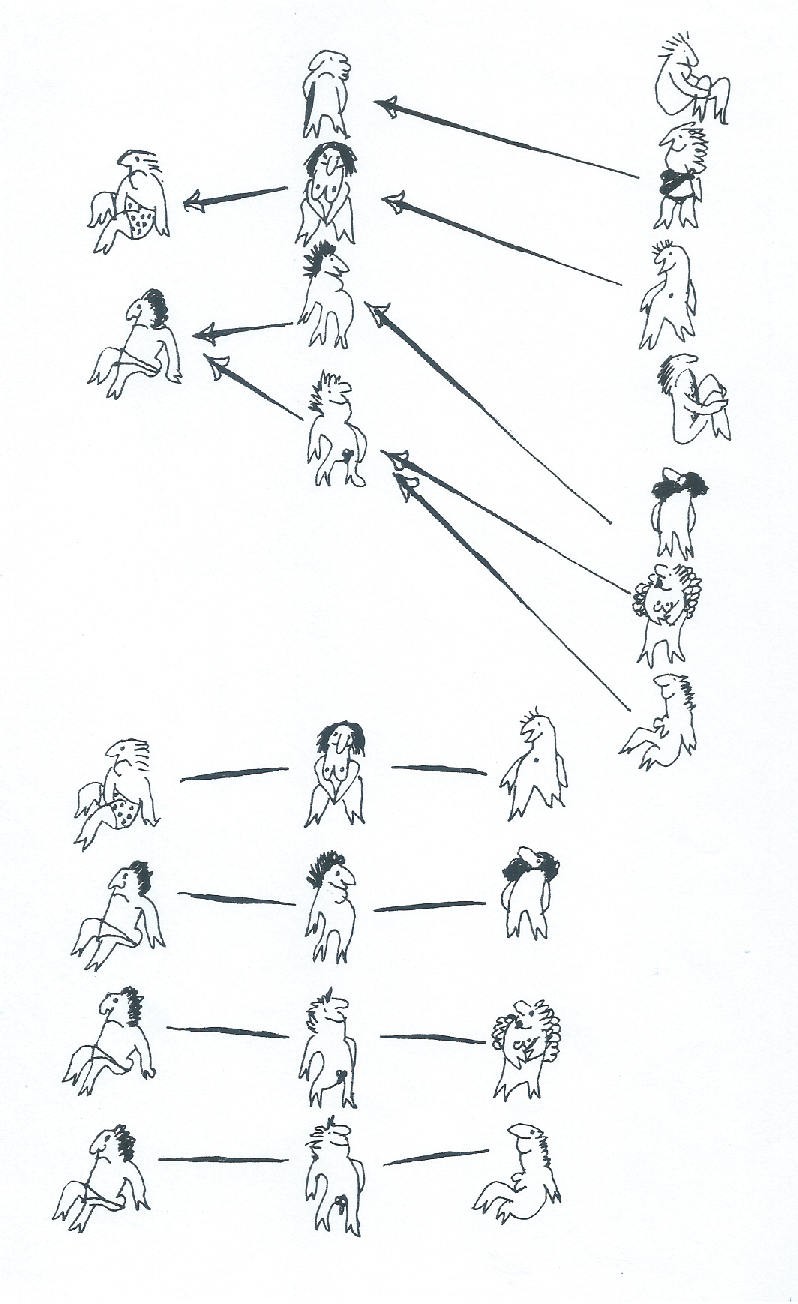 scaled 500} }
\centerline{Planche h}
\vskip2mm
\centerline{Le lemme des mariages}
\vfill
\null\vfill\eject
%
\fancyhead[LO]{A. M. \quad\quad Uniformisation des surfaces de Riemann}
  \fancyhead[RO]{\thepage}
  \fancyhead[RE]{\small \quad\quad\quad Epilogue\hfill Glossaire \hfill}
  \fancyhead[LE]{\thepage}
%
\def\page #1{\leaders\hbox to 2mm{\hfil.\hfil}\hfill
\rlap{\hbox to 11mm{\hfill #1}}\par}
\rightskip=10mm
{\vsize=20,5cm
\null\vskip-5mm
\centerline{\bf Glossaire}
\parc
Sans reprendre les d\'enominations des objets et morphismes standards (p.{\romannumeral 4} et {\romannumeral 5}),
cette liste donne la localisation des d\'efinitions des \og mots et symboles propres\fg du texte.
\finc
\parcs\it{
\centerline{Lexique}
acyclique \page {14}
arc analytique \page{12}
arrondie \page{{\romannumeral 3} {\small et\/} 21}
ar\^ete (\/{\rm d'un $\omega\/$-graphe\/}) \page{12}
bord (\/{\rm d'une surface de Riemann\/}) \page {16}
bord (\/{\rm d'un $\omega\/$-arc\/}) \page {11}
bord (\/{\rm d'un $\omega\/$-polygone ouvert ou ferm\'e\/}) \page {{\romannumeral 3} {\small et\/}  19}
carte \page{12}
cellule \page {4}
collier (param\'etr\'e) (\/{\rm d'une $\omega\/$-cha{\^{\i}}ne\/})\page {14}
coin \page {{\romannumeral 3}}
conjugaison (\/{\rm de double\/}) \page {15}
connexit\'e \page {\romannumeral 5}
courbe simple (ferm\'ee) \page {14}
coordonn\'ee (\/{\rm en un point d'une surface de Riemann\/}) \page {12}
cycle analytique\page {{\romannumeral 2} {\rm et\/} 12}
d\'esir (syst\'ematique)\page {\uppercase\expandafter{\romannumeral 4}}
domaine de Jordan \page {\romannumeral 1}
double (de Klein) (\/{\rm d'une surface de Riemann\/}) \page {{\romannumeral 3} {\small et\/} 16}
double (de Klein) (\/{\rm d'un $\omega\/$-polygone\/} sur un bord $\beta\/$) \page {{\romannumeral 3} {\small et\/} 19}
duplicata \page {14}
e\'ennodrooque \page {12}
\'el\'ementaire (\/{\small ouvert dans une surface de Riemann\/}\/) \page {\romannumeral 2}
\'el\'ementaire (\/{\small carte ou coordonn\'ee dans une surface de Riemann \`a bord\/}\/) \page {16}
ferm\'e (\/{\small $\omega\/$-graphe\/}\/) \page {12}
filtrante (\/{\rm famille (de parties)\/}) \page{22}
graphie  \page {\romannumeral 5}
h\'et\'erographie  \page {\romannumeral 5}
h\'et\'erolution \page {\romannumeral 3}
homographie  \page {\romannumeral 5}
int\'erieur formel (\/{\small d'un $\omega\/$-arc\/}\/)\page {12}
isomorphisme \page {\romannumeral 1}
Jordan-simplement connexe \page {14}
limite projective (\/{\rm d'un d\'esir syst\'ematque\/})\page {\uppercase\expandafter{\romannumeral 5}}
normale (\/{\small e\'ennodrooque  $\Gamma\/$-normale\/}\/)\page {11}
p\'eriph\'erique (\/{\rm cellule dans un $\omega\/$-polygone ferm\'e\/}) \page {3}
perle  \page {\romannumeral 4}
perle sur une arr\^ete \page {14}
planaire (\/{\small surface de Riemann\/}\/) \page {\romannumeral 2}
polygone analytique (\/{\small ouvert \/}\/)\page {{\romannumeral 2}}
polygone analytique (\/{\small ferm\'e \/}\/)\page {4}
pr\'etendant \page {\uppercase\expandafter{\romannumeral 5}}
r\'ealisation (\/{\rm d'un d\'esir syst\'ematique\/})\page {\uppercase\expandafter{\romannumeral 5}}
r\'eduit (\/{\rm support r\'eduit d'un support de graphe\/})\page {14
}
signe locaux (\/{\rm groupe des\/})\page {\romannumeral 5}
simple (\/{\rm $\omega\/$-polygone\/})\page {4}
simlement connexe \page {\romannumeral 1}
simplement connexe (\/{\small aux sens de Jordan\/}\/) \page {14}
sommet  (\/{\small d'un $\omega\/$-polygone et d'un $\omega\/$-graphe\/}\/) \page {{\romannumeral 3} {\rm et\/} 12}
standard (\/{\small ouvert dans une surface de Riemann\/}\/) \page {\romannumeral 2}
support (\/{\small d'un $\omega\/$-graphe\/}\/) \page {12}
surface de Riemann \page {\romannumeral 1}
surface de Riemann \`a bord \page {16}
valence (\/{\rm d'un point dans un $\omega\/$-graphe d'une surface de Riemann \/}) \page{11}
$\omega\/$-arc \page {12}
$\omega\/$-cha{\^{\i}}ne \page {12}
$\omega\/$-cycle \page {{\romannumeral 2} {\rm et\/} 12}
$\omega\/$-cycle d'un signe \page {14}
$\omega\/$-graphe \page {11}
$\omega\/$-polygone (\/{\small ouvert et ferm\'e \/}\/)\page {{\romannumeral 2} {\rm et\/} 3}
\vskip3mm
\centerline{symblolique}
$\partial\,{\bf R\/},\, \partial\,U\/$ \page {16}
$i_\beta\/$
$({\cal D\/}_\beta\,U,\,{\bf U\/}_{^a\beta},\,\sigma_\beta)\/$ \page {\romannumeral 3}
${\cal B\/}\/$ \page {12}
$\Gamma_{\cal B\/}\/$ \page {12}
$\Delta_{\cal B\/}\/$ \page {12}
${\goth p\/}\/$, ${\goth p\/}_\lambda\/$ \page {14}
$\partial\Gamma\/$ \page {14}
$\pi_{\cal B\/}\/$ \page {14}
$\pi_\Gamma\/$ \page {14}
${\goth S\/}_\Gamma\/$ \page {14}
$\Gamma^\eta\/$,  $\Gamma^\eta_U\/$ \page {14}
$U_{\epsilon,\,  \gamma,\, i}\/$ \page {18}
}\fincs
}
\vfill\eject
%
\fancyhead[LO]{A. M. \quad\quad Uniformisation des surfaces de Riemann}
  \fancyhead[RO]{\thepage}
  \fancyhead[RE]{\small\quad\quad\quad Epilogue\hfill G\'en\'erique de fin\hfill}
  \fancyhead[LE]{\thepage}
\centerline{\bf Table des mati\`eres}
\vskip20mm
\centerline{\bf Pr\'eambule : Aper\c cu sur quatre paragraphes\dots}
\vskip2mm
Introduction \page {\romannumeral 1}
La preuve et ses sept \'enonc\'es, le cas compact \page {\romannumeral 2}
Arrondissement des brisures du bord et le cas g\'en\'eral \page {\romannumeral 3}
Objets et morphismes standards dans le plan et la sph\`ere \page {\romannumeral 4}
\vskip5mm
\centerline{\bf D\'eroulement de la D\'emonstration\/}
\vskip2mm
 \S 1 D\'emonstration du Lemme A  \page {2}
\S 2 Le Lemme B et d\'emonstration du Lemme C  \page {4}
\S 3 D\'emonstration du Lemme D et Corollaire D  \page {8}
\S 4 D\'emonstration des Th\'eor\`emes \page {11}
\vskip5mm
\centerline{\bf Appendices}
\vskip2mm
{\bf 1\/} Arcs analytiques, duplicatas et s\'eparation par les $\omega\/$-graphes \page {12}
{\bf 2\/} Surfaces de Riemann \`a bord et double de Klein usuel \page {16}
{\bf 3\/} Bord et double de Klein des $\omega\/$-polygones \page {18}
{\bf 4\/} Une caract\'erisation de cartes dans les surfaces de Riemann \page {22}
\vskip5mm
\centerline{\bf Epilogue}
\vskip2mm
Commentaires sur lit\'erature et tradition orale \page{\uppercase\expandafter{\romannumeral 1}}
R\'ef\'erences \page {\uppercase\expandafter{\romannumeral 2}}
Notes \page {\uppercase\expandafter{\romannumeral 3}}
Glossaire \page {\uppercase\expandafter{\romannumeral 7}}
Table des mati\`eres \page {\uppercase\expandafter{\romannumeral 9}}

G\'en\'erique de fin \page {\uppercase\expandafter{\romannumeral 10}}
Contacts et affiliation \page {\uppercase\expandafter{\romannumeral 11}}
\vfill\eject
\null\vskip5mm

{\leftskip-8mm\rightskip-2mm\parindent=0pt
{\small\it\og Klein, Poincar\'e und Koebe
ist es vor allem zu verdanken,
wenn heute die Theorie
der Uniformisierung, welche innerhalb der komplexen Funktiontheorie eine zentrale
Stellung beanspruchen darf, als ein mathematisches Geb\" aude von bezonderer Harmonie
und Gro\stz z\" ugigkeit 
vor uns steht.\fg\/}
\vskip2mm
\hfill{\small  \hfill H. Weyl, ({\small\sl Die Idee der Riemannschen Fl\" ache} (1913), p. 142)}
\par}
\nobreak\parindent=20pt%
\vskip 25mm

\setcounter{footnote}{0}
Les quatres paragraphes et autant d'appendices de ce texte n'existe\-raient pas sans
les \oe uvres explicitement ou implicitement cit\'ees,
diverses demandes d'une r\'ealisation de leur r\'edaction, ni
la patience de L. Guillou pour nos balbutiements et son
\'erudition communicative de cette \og litt\'erature math\'ematique\fg dont l'auteur,
avant d'entreprendre sous sa direction ce T.E.R.%
\footnote%
{Travail d'\'etude et de recherche, exercice de style de nature initiatique, introduit et exig\'e
chez des tribus extratropicales (bien au nord du Cancer et tr\`es \`a l'est de l'Amazonie) d\`es la fin
du vingti\`eme si\`ecle, pour  faire passer les licencieuces au statut de ma{\^{\i}}tresses
de Math\'ematique.}
d'histoire des math\'ematiques, n'imaginaient pas l'existence.

L'auteur manifeste aussi sa sympathie aux doigts et \`a l'\oe il de monsieur Carreira, 
qui  effectu\`erent coifures, costumes, maquillages,
et firent passer de  nos l\`evres \`a ce \og manuscript frapp\'e\fg une discr\`ete couche de \TeX.

\vfill
{\leftskip-8mm\rightskip-2mm\parindent=0pt
{\small\it \og J'ai toujours regrett\'e de ne pas l'avoir connu en pleine jeunesse quand, brun et
basan\'e \`a l'image d'un conquistador et tout fr\'e\-mis\-sant des perspectives scientifiques
qu'ouvraient la psychologie du XIX$^{\hbox{\it e\/}}\/$ si\`ecle, il \'etait parti \`a la conqu\^ete
spirituelle du nouveau monde\fg\/}
\vskip-2mm
\hfill{\small \hfill C. L\`evi-Strauss,  ({\small\sl Tristes tropiques\/})}
\vskip6mm
\centerline{ (et de ses) non encore uniformis\`ees}
\vskip2mm
\centerline{${\goth ante-Koebiennes\/}\ {\goth Surfaces\/}\ {\goth de\/}\  {\goth Riemann}\/$}
\vfill\eject
%
\null\vfill\vfill\vfill\vfill\vfill\vfill\vfill\vfill\vfill\vfill\vfill
\parc
\noindent{
Metteur en sc\`ene et secr\'etaire~: Alexis Marin Bozonat\hfill\break
\null\ courriel~: alexis.charles.marin@gmail.com}

\noindent{
Univ. Grenoble Alpes, CNRS, IF, 38100, Grenoble, France}
\finc

\vfill

\enddocument